\documentclass[reqno]{amsart}
\usepackage{amssymb}
\usepackage{amsmath}
\usepackage{url}

\def\scaledpicture#1by#2(#3scaled#4){{
\dimen0=#1  \dimen1=#2
\divide\dimen0 by 1000 \multiply\dimen0 by #4
\divide\dimen1 by 1000 \multiply\dimen1 by #4
\picture \dimen0 by \dimen1 (#3 scaled #4)}}
\def\dfigure#1by#2(#3scaled#4offset#5:#6)
    {\medskip
     \vglue 2mm minus 2mm
     $$
       \hbox{
         \hglue#5
         {\scaledpicture #1 by #2 (#3 scaled #4)}
       }
     $$
     \par\goodbreak
     \vglue 2mm minus 2mm
     \medskip}

\def\qmod#1#2{{\hbox{}^{\displaystyle{#1}}}\!\big/\!\hbox{}_{
\displaystyle{#2}}}

\def\resto#1#2{{
#1\hskip 0.4ex\vline_{\hskip 0.2ex\raisebox{-0,2ex}
{{${\scriptstyle #2}$}}}}}

\def\A{{\mathbb A}}
\def\B{{\mathbb B}}
\def\C{{\mathbb C}}

\def\E{{\mathbb E}}
\def\F{{\mathbb F}}

\def\H{{\mathbb H}}

\def\N{{\mathbb N}}

\def\P{{\mathbb P}}
\def\Q{{\mathbb Q}}
\def\R{{\mathbb R}}
\def\Z{{\mathbb Z}}

\def\cringle{\mathaccent23}
\def\union{\mathop{\bigcup}}
\def\qed {\hfill\vrule height6pt width6pt depth0pt \smallskip}

\def\map{\longrightarrow}
\def\textmap#1{\mathop{\vbox{\ialign{
                                  ##\crcr
      ${\scriptstyle\hfil\;\;#1\;\;\hfil}$\crcr
      \noalign{\kern 1pt\nointerlineskip}
      \rightarrowfill\crcr}}\;}}
\def\bigtextmap#1{\mathop{\vbox{\ialign{
                                  ##\crcr
      ${\hfil\;\;#1\;\;\hfil}$\crcr
      \noalign{\kern 1pt\nointerlineskip}
      \rightarrowfill\crcr}}\;}}
      
\newcommand{\cal}{\mathcal}
\def\textlmap#1{\mathop{\vbox{\ialign{
                                  ##\crcr
      ${\scriptstyle\hfil\;\;#1\;\;\hfil}$\crcr
      \noalign{\kern-1pt\nointerlineskip}
      \leftarrowfill\crcr}}\;}}
 
\def\ag{{\mathfrak a}}
\def\bg{{\mathfrak b}}
\def\cg{{\mathfrak c}}

\def\g{{\mathfrak g}}

\def\tg{{\mathfrak t}}

\def\Ag{{\mathfrak A}}

\def\Rg{{\mathfrak R}}

\newtheorem{sz}{Satz}[section]
\newtheorem{thry}[sz]{Theorem}
\newtheorem{pr}[sz]{Proposition}
\newtheorem{re}[sz]{Remark}
\newtheorem{co}[sz]{Corollary}
\newtheorem{dt}[sz]{Definition}
\newtheorem{lm}[sz]{Lemma}

\begin{document}
\def\Pr{{\rm Pr}}
\def\tr{{\rm Tr}}
\def\End{{\rm End}}
\def\Aut{{\rm Aut}}
\def\Spin{{\rm Spin}}
\def\U{{\rm U}}
\def\SU{{\rm SU}}
\def\SO{{\rm SO}}
\def\PU{{\rm PU}}
\def\GL{{\rm GL}}
\def\spin{{\rm spin}}
\def\u{{\rm u}}
\def\su{{\rm su}}
\def\so{{\rm so}}
\def\ub{\underbar}
\def\pu{{\rm pu}}
\def\Pic{{\rm Pic}}
\def\Iso{{\rm Iso}}
\def\NS{{\rm NS}}
\def\deg{{\rm deg}}
\def\Hom{{\rm Hom}}
\def\Aut{{\rm Aut}}
\def\h{{\mathfrak h}}
\def\Herm{{\rm Herm}}
\def\Vol{{\rm Vol}}
\def\pf{{\bf Proof: }}
\def\id{{\rm id}}
\def\i{{\mathfrak i}}
\def\im{{\rm im}}
\def\rk{{\rm rk}}
\def\ad{{\rm ad}}
\def\h{{\bf H}}
\def\coker{{\rm coker}}
\def\dbar{\bar{\partial}}
\def\Lo{{\Lambda_g}}
\def\niq{=\kern-.18cm /\kern.08cm}
\def\Ad{{\rm Ad}}
\def\RSU{\R SU}
\def\ad{{\rm ad}}
\def\dva{\bar\partial_A}
\def\da{\partial_A}
\def\p{{\rm p}}
\def\sp{\Sigma^{+}}
\def\sm{\Sigma^{-}}
\def\spm{\Sigma^{\pm}}
\def\smp{\Sigma^{\mp}}
\def\oo{{\scriptstyle{\cal O}}}
\def\ooo{{\scriptscriptstyle{\cal O}}}
\def\sw{Seiberg-Witten }
\def\pa{\partial_A\bar\partial_A}
\def\Dr{{\raisebox{0.17ex}{$\not$}}{\hskip -1pt {D}}}
\def\gr{{\scriptscriptstyle|}\hskip -4pt{\g}}
\def\subsetint{{\  {\subset}\hskip -2.45mm{\raisebox{.28ex}
{$\scriptscriptstyle\subset$}}\ }}
\def\nr{\parallel}
\def\ra{\rightarrow}
\def\Ob{{\rm Ob}}
\def\t{\mathrm t}

\def\varLim@#1#2{%
\vtop{\m@th\ialign{##\cr
\hfil$#1\operator@font Lim$\hfil\cr
\noalign{\nointerlineskip\kern1.5\ex@}#2\cr
\noalign{\nointerlineskip\kern-\ex@}\cr}}%
}
\def\varinjLim{%
\mathop{\mathpalette\varLim@{\rightarrowfill@\textstyle}}\nmlimits@
}

\title[Cohomotopy invariants]
{Cohomotopy invariants and the universal cohomotopy invariant jump formula}\author{Christian Okonek,   Andrei Teleman}

\begin{abstract}  Starting from ideas of Furuta, we develop a general formalism for the construction of cohomotopy invariants associated with a certain class of $S^1$-equivariant non-linear  maps between Hilbert bundles.   Applied to the Seiberg-Witten map, this formalism yields a new class of cohomotopy Seiberg-Witten invariants which have clear functorial properties with respect to diffeomorphisms of 4-manifolds.   Our invariants and the Bauer-Furuta classes are directly comparable for 4-manifolds with $b_1=0$; they are equivalent when $b_1=0$ and $b_+>1$, but are finer in the case $b_1=0$, $b_+=1$ (they detect the wall-crossing phenomena).

 We study fundamental properties of the new invariants in a very general framework. In particular we prove a universal cohomotopy invariant jump formula and a multiplicative property. The formalism applies to other gauge theoretical problems, e.g. to the theory of gauge theoretical (Hamiltonian) Gromov-Witten invariants.
\end{abstract}
\maketitle

\tableofcontents

\section{Introduction}\label{intro}

\subsection{Motivation}\label{mot} The goal of this article is to develop a general formalism for the construction of cohomotopy invariants associated with a certain class of non-linear  maps between Hilbert bundles.   The main example we have in mind is the Seiberg-Witten map, but the formalism applies to other interesting classes of maps related to gauge theoretical problems as well. 

The first stable-homotopy Seiberg-Witten invariants have been introduced independently by  M. Furuta and S. Bauer. Furuta first used ``finite dimensional approximations" of the monopole map in his work on the 11/8 conjecture \cite{Fu1}, and then introduced a class of refined Seiberg-Witten invariants (called   ``stable homotopy version of the Seiberg-Witten invariants") in a geometric, non-formalized way in \cite{Fu2}.  In this preprint  Furuta acknowledges  independent work by Bauer \cite{B3}.
According to Furuta, the new invariants belong  to a certain inductive limit of sets of homotopy classes  of maps associated with ``finite dimensional approximations" of the Seiberg-Witten map. The structure  and the functorial properties of this inductive limit (with respect to diffeomorphisms between 4-manifolds) have not been worked out in this article.  A precise version of the new invariants has been introduced later by Bauer-Furuta in \cite{BF}: the Bauer-Furuta classes belong to  certain stable cohomotopy groups associated with a presentation $(E,F)$ of the K-theory element  ${\rm ind}(\Dr)$ defined by a fixed $Spin^c$-structure.  This element   ${\rm ind}(\Dr)$  belongs to the K-theory group $K(B)$, where $B=H^1(X;\R)/H^1(X;\Z)$ is the Picard group  of the base manifold $X$. 

In this article we propose a different construction of cohomotopy invariants which has the following advantages: Our construction yields a larger class of invariants, which are well defined in all cases, are always finer than the classical invariants, and have clear functorial properties. In order to explain  the advantages of the new formalism in a non-technical  way, we consider again the Seiberg-Witten case.

It is well known that the Seiberg-Witten map $\mu$ can be regarded as an $S^1$-equivariant bundle map between Hilbert bundles over the torus $B$ (see  \cite{BF} and section \ref{swmap} of this article).  We first  choose the perturbing form
in the second Seiberg-Witten equation in the ``bad way", i.e. such that the equations have reducible solutions (solutions with trivial spinor component); we make this ``bad choice" even in the case $b_+(X)>1$! In ``classical" Seiberg-Witten theory one perturbs the second Seiberg-Witten equation using a nontrivial self-dual harmonic form $\kappa\in i\H^+\setminus\{0\}$, and gets a new map $\mu_\kappa$ which defines a moduli space which does not contain   reductions.  Instead of a constant perturbation $\kappa$, we consider a map $\kappa:B\to i\H^+\setminus\{0\}$, and perturb the Seiberg-Witten map $\mu$ (regarded as bundle map over $B$) using this map.  The associated invariant will depend  on the homotopy class $[\kappa]\in [B, S(i\H^+)]$.  This leads to the following questions:
\begin{enumerate}
\item Does one obtain new invariants in this way?
\item If so, does one have a {\it universal cohomotopy invariant  jump formula}, i.e. a formula which describes the jump of the cohomotopy  invariant when one passes from one homotopy class to another?
\item Use again constant perturbation forms $\kappa$, but let $\kappa$ vary in the sphere $S(i\H^+)$ and regard the obtained map $\tilde \mu$ as an $S^1$-equivariant  bundle map  over the larger basis $B\times S(i\H^+)$. Does this universal perturbation $\tilde\mu$ yield more differential topological information than the individual perturbations $\mu_\kappa$? If not, express the cohomotopy invariant associated with $\tilde \mu$ in terms of the invariant associated with  $\mu_\kappa$ and topological invariants of   $X$.
 \end{enumerate}
These questions are interesting as soon as $b_1\geq b_+-1$ (even for $b_+>1$!)  and they are also interesting {\it for the classical invariant}, because for non-constant perturbations $\kappa$ one gets new   Seiberg-Witten type moduli spaces. The universal wall-crossing formula \cite{LL}, \cite{OO}, \cite{OT} for the  {\it full Seiberg-Witten invariant}
\footnote{The full Seiberg-Witten invariant \cite{OT} is an element in $\Lambda^*(H^1(X;\Z))$. The numerical Seiberg-Witten invariant (the original invariant introduced by Witten) is the degree 0 term of the full invariant.}
    should be a formal consequence of a  universal  cohomotopy invariant jump formula. These questions  will be completely answered in this article. 
  
  Another motivation for proposing a new formalism was the need to have {\it  well defined} invariants,  with clear {\it functorial properties}. Recall that the classical full Seiberg-Witten  invariant can be regarded as an element of $[\wedge^* H^1(X,\Z)]^{Spin^c(X)}$, where $Spin^c(X)$ denotes the  torsor of equivalence classes of $Spin^c$-structures.  Therefore this invariant  belongs to a group which is obviously functorial with respect to pairs $(h,\theta)$ consisting of an orientation preserving homotopy equivalence $h:X\to X'$, and a bijection $\theta:Spin^c(X')\to Spin^c(X)$ which is compatible with the Chern class maps $Spin^c(X)\to H^2(X,\Z)$, $Spin^c(X')\to H^2(X',\Z)$ and the $H^2(X,\Z)$, $H^2(X',\Z)$-actions  on the two sets.  Such a pair will be called a $Spin^c$-homotopy equivalence.  We will say that an assignement $X\mapsto G(X)\in{\cal A}b$  is {\it topologically functorial}  on the category of smooth 4-manifolds if it is functorial with respect to  $Spin^c$-homotopy equivalences. It is natural to require that the refined Seiberg-Witten invariant   belongs to a  group  $G(\cdot)$  which  is topologically functorial, as it is the case for the classical invariant. In other words, we want the group to which the invariant belongs to have much stronger functorial properties than the invariant itself. This  is  important for practical reasons; for instance, if one wants to classify the $Spin^c$-homotopy equivalences $X\to X$ which are realized by diffeomorphisms, one will essentially need the {\it topological} functoriality of the group to which the invariant belongs.

The definition of the stable cohomotopy group  used in \cite{BF} depends on the choice of a presentation $(E,B\times\C^n)$   of   ${\rm ind}(\Dr)\in K(B)$ (see \cite{BF} p. 8-9). Since in general such a presentation has homotopically non-trivial automorphisms, the obtained cohomotopy groups cannot be regarded as invariants of the K-theory element ${\rm ind}(\Dr)\in K(B)$. This makes it   difficult to control the  functorial properties of the Bauer-Furuta stable cohomotopy groups as defined in \cite{BF} with respect to homeomorphisms (or even  diffeomorphisms) of 4-manifolds, and to understand in which sense  the constructed class is well defined.

Using Segal cocycles instead of finite rank presentations (\cite{BF} p. 7-8) does not remove the problem, because of monodromy phenomena in the space of Segal cocycles\footnote{Contrary to what is  often stated in the literature,  the space of possible Dirac operators associated with a fixed equivalence class of $Spin^c$-structures is not contractible (see section \ref{swmap}). So even if one considers only $Spin^c$-Dirac operators, one does not get   a contractible space  of Segal cocycles. }. 
A similar difficulty   concerns the concept   ``Thom spectrum of a virtual bundle", used by Bauer-Furuta (see \cite{BF} p. 8) and other authors in order to give a geometric  interpretation of the Bauer-Furuta classes.  One can indeed associate a Thom spectrum to a {\it fixed presentation} $(E,B\times\C^n)$  of a K-theory element $x\in K(B)$, but unfortunately {\it not to $x$ itself.}  
For 4-manifolds with $b_1=0$, the Bauer-Furuta class  gives a well defined invariant, which can easily  be   identified with the image of  our invariant under a boundary morphism of cohomotopy groups.   The two invariants are equivalent when $b_1=0$, $b_+>1$. 

Note that an elegant  construction described by Furuta \cite{Fu3} leads to a well-defined invariant belonging to a group which is  ${\cal C}^\infty$-functorial for  arbitrary $b_1$ and  $b_+>1$. Furuta uses the {\it universal family of Dirac operators} associated with a metric and a class  of $Spin^c$-structures $\cg\in Spin^c(X)$  in order to remove the ambiguity in the choice of a $Spin^c$-structure $\tau\in\cg$ and get a well-defined Segal cocyle.   We will explain this  formalism in section \ref{swmap}.\\

 Our new point of view has the following advantages:
\begin{enumerate}
\item The new cohomotopy  Seiberg-Witten invariants  are finer than the full classical Seiberg-Witten invariants {\it in all cases}, including the case of manifolds with $b_1\geq b_+-1$ and including  the invariants associated with non-constant perturbations $\kappa:B\to  i\H^2_+\setminus\{0\}$. In the case $b_1\geq b_+-1$ we prove a universal cohomotopy invariant jump formula; the  universal wall-crossing formula for the classical invariant is a formal consequence of this result.
\item  Our invariant belongs to a cohomotopy group which is intrinsically associated with the base 4-manifold, and  is topologically  functorial in the sense explained above.
\end{enumerate}
\vspace{1mm}
{\bf Remark:} An interesting development in cohomotopy Seiberg-Witten theory concerns invariants defined for families of 4-manifolds parameterized by a compact space  \cite{Fu2}. Most parts of our construction generalize immediately to this situation; note however that in the family case the map $\kappa$ above has to be replaced by a section of a certain sphere bundle.

\subsection{Summary of results}
  In the first section we construct a graded cohomotopy group  associated with a  K-theory element  $x\in K(B)$.   To every representative $(E,F)\in x$ we associate the graded  group ${_{S^1}\alpha}^*_B(S(E)_{+B},F^+_B)$, where  ${_{S^1}\alpha}^*_B(\cdot,\cdot)$ stands for  the   $S^1$-equivariant stable cohomotopy group functor on the category of pointed $S^1$-spaces over $B$; it is obtained by stabilizing  with spaces of the form $(\eta\oplus\xi_0)^+_B$, where $\eta$ is a a complex, and $\xi_0$ a real bundle.  Note that we do not use all  characters of $S^1$ in the stabilizing process; for this reason we do not use the standard notation  ${_{S^1}\omega}^*_B$ found in the literature \cite{CJ}.
We define $\alpha^*(x)$ to be the inductive limit of the graded groups ${_{S^1}\alpha}^*_B(S(E)_{+B},F^+_B)$ with respect to the category ${\cal T}(x)$ of representatives $(E,F)$ of $x$. Since ${\cal T}(x)$ is not a small filtering category (see \cite{AM}, and  section \ref{limits} below),  this  limit cannot be obtained using the classical construction. It will be constructed in two steps: First we stabilize the graded group ${_{S^1}\alpha}^*_B(S(E)_{+B},F^+_B)$ with respect to  standard representative enlargements  $(E,F)\mapsto (E\oplus U,F\oplus U)$, and we obtain a new graded group   $\hat \alpha^*( E,F)$, which still depends on the fixed pair $(E,F)$. The groups $\hat \alpha^*( E,F)$, $\hat \alpha^*( E',F')$ defined by two representatives  $(E,F)$, $(E',F')$ of $x$ are {\it non-canonically} isomorphic. The group $\alpha^*(x)$ will be the quotient of $\hat \alpha^*( E,F)$ by the equivalence relation generated by the inductive limit  of the automorphism groups $\Aut(E\oplus U)\times \Aut(F\oplus U)$.  We give an explicit description of $\alpha^*(x)$ as a quotient of the group $\hat \alpha^*( E,F)$ associated with any representative $(E,F)$ by the action of the image of the $J$-homomorphism ${_{S^1}J}:K^{-1}(B)\to  {_{S^1}\alpha}^0(B)^\times$ in the group of units ${_{S^1} \alpha}^0(B)^\times$ of the ground ring ${_{S^1}\alpha}^0(B): ={_{S^1}\alpha}^0_B(B_{+B},B_{+B})$. In other words, we are able to control the effect of bundle automorphisms in our inductive limit and we obtain a graded group which is intrinsically associated with the K-theory element $x$. We believe that this construction is of independent   interest from the point of view of homotopy theory.

A way to understand the role of the graded group $\alpha^*(x)$ is the following: Because of  the presence of homotopically non-trivial bundle automorphisms, one cannot define the projectivization $\P(x)$ of a K-theory element $x\in K(B)$ (neither in the category of topological spaces nor in the category of spectra). {\it The graded group $\alpha^*(x)$ plays the role of  what should be the  cohomotopy group of a formal projectivization of the K-theory element $x$.}   

In the second section we first introduce a distinguished class of non-linear  maps $\mu$ between Hilbert bundles over a compact base $B$. The $\C$-linear part of the linearization  of such a map $\mu$ at the zero section is a linear Fredholm operator, so it defines a K-theory  element $x\in K(B)$. The goal of the section is the construction of an invariant $\{\mu\}\in \alpha^*(x)$. This invariant is constructed using {\it finite dimensional approximations} of the map $\mu$. In order to get these approximations we make   use  of the retractions $\rho_A:{\cal A}^+\setminus S(A^\bot)\to A^+$ associated with finite dimensional subspaces $A$ of a Hilbert space ${\cal A}$,  as in \cite{BF}. This method to construct finite dimensional approximations   applies to a very large class of non-linear maps, whereas   Furuta's original method based on $L^2$-orthogonal projections on direct sums of eigenspaces (see \cite{Fu2})  is limited  to maps whose linearizations are elliptic differential operators.

The main difference between our definition and the construction of the Bauer-Furuta classes  given in \cite{BF}, is that 
\begin{enumerate}
\item our construction  uses only spaces fibered over the base $B$. In particular we avoid using Thom spaces,
\item  we  treat the real and the complex summands in our finite dimensional approximations separately. 
\end{enumerate}
Therefore, from this point of view, our construction is {\it closer to the original ideas of Furuta} \cite{Fu2}.  Having the finite dimensional approximation, a representative of the invariant is  an element in a group of the form ${_{S^1}\alpha}^*_B(S(E)_{+B},F^+_B)$ obtained by a  simple geometric construction,  which we call {\it the cylinder construction}\footnote{After completing the first version of this article we found out 
about the results \cite{C} on Leray-Schauder index theory.  Under the assumption that the group-action is free on a neighborhood of the  zero-set of the vector field, one can define a refinement  of the usual Poincar\'e-Hopf vector field index, which is probably related to our refinement of the Bauer-Furuta class.}. The  class  obtained in this way carries more information than the one defined in \cite{BF}.   In    \ref{swmap} we show that the Seiberg-Witten map associated with a Spin$^c$-structure $\tau$ on a Riemannian 4-manifold $M$ with $b_+(M)>0$ yields a non-linear Fredholm map $sw_\kappa$ (depending on a  twisting map   $\kappa:B=H^1(X;\R)/H^1(X;\Z)\to i\H^2_+\setminus\{0\}$) which belongs to our distinguished class of maps. Hence the general theory applies and a yields a cohomotopy Seiberg-Witten invariant $\{sw_\kappa\}\in\alpha^{b_+(M)-1}({\rm ind}(\Dr))$, which only depends of the homotopy class of $\kappa$. Our construction of the bundle map $sw_\kappa$ is  different from  and somewhat simpler   than  the one given in \cite{BF}.   In this section we also explain why the space of Dirac operators associated with a fixed equivalence class of $Spin^c$-structures is not contractible, so there is no way to distinguish a contractible class of Segal cocycles defining the $K$-theory element ${\rm ind}(\Dr)$. This makes clear why the cohomotopy groups defined in \cite{BF} (which depend essentially on the choice of a Segal cocycle) cannot be regarded as intrinsic (functorial) invariants of the base manifold.

In the third section we prove several fundamental properties of the invariant $\{\mu\}\in \alpha^*(x)$ in our general, abstract framework: 
\begin{enumerate}
\item We study the image of our invariant under the Hurewicz morphism,  and we prove that the Poincar\'e dual of this image can be identified with  the virtual fundamental class of the vanishing locus. In other words, the full homology invariant associated with the virtual  fundamental class of the ``moduli space" (i.e.  the $S^1$-quotient of the vanishing locus of   $\mu$) can be identified with the Hurewicz image of the cohomotopy invariant.  Moreover, the Hurewicz morphism is an isomorphism when the ``expected dimension" vanishes.
\item We prove a formal  universal cohomotopy invariant jump formula for our refined cohomotopy invariants.
\item We prove a general product formula for the invariant $\{\mu_1\times\mu_2\}$ associated with  a product of maps; in this formula we allow one of the factors to have zeros on the $S^1$-fixed point locus.  When both factors are nowhere vanishing on their fixed point loci, we prove a vanishing result   for the Hurewicz image of the invariant.
\end{enumerate}

Specialized  to the Seiberg-Witten map, these properties automatically yield important results for the new cohomotopy Seiberg-Witten invariants.  

The first result shows that the cohomotopy Seiberg-Witten invariant is a refinement of the classical {\it full} Seiberg-Witten invariants in all cases. Combined with the second property, this also yields a universal invariant jump formula for the full classical Seiberg-Witten invariant in the case $b_1(X)\geq b_+(X)-1$. 

The vanishing result in (2) reproves the classical vanishing theorem for the Seiberg-Witten invariant of a direct sum in the case where both summands $X_i$ have $b_+(X_i)>0$.   

The third  result provides the topological formalism for proving a formula for the cohomotopy  invariant of a connected sum of two 4-manifolds, even in the case when one term of the sum has $b_+=0$. The analytic techniques required for this general gluing formula are not discussed in this article.  
\\
\\
{\bf Acknowledgment:} The first author would like to thank  Tammo tom Dieck for useful and interesting discussions at the beginning of this project.   The second author is indebted to Mikio Furuta for  useful discussions and for explaining the construction of the ``universal Seiberg-Witten map" described in section  \ref{swmap}.

\section{Cohomotopy groups associated with elements in $K(B)$}\label{cohomotopy}

\subsection{Definition of ${_{S^1}\alpha}^*_B(X,Y)$}\label{definition}
Let $B$ be a compact topological space endowed with the trivial $S^1$-action. Let ${\cal C}_B$  be the category defined in the following way: the objects of ${\cal C}_B$ are vector bundles over $B$ of the form
$$\xi=\eta\oplus \xi_0 \ ,
$$
where $\eta$ is a complex vector bundle endowed with the standard $S^1$-action and $\xi_0$ is a real vector bundle  endowed with the trivial $S^1$-action; for two objects $\xi=\eta\oplus \xi_0$, $\xi=\eta'\oplus \xi_0'$ a morphism $\nu:\xi\to \xi'$ is a pair $(i,\zeta)$ consisting of an $S^1$-equivariant bundle embedding $i=\iota\oplus i_0:\xi\to \xi'$ and a complement $\zeta=\nu\oplus \zeta_0$ of $i(\xi)$ in $\xi'$.  Composition of morphisms is defined in a natural way. A morphism $u=(i,\zeta):\xi\to\xi'$ defines a push-forward morphism $A(u):A(\xi)\to A(\xi')$, where $A(\xi):=A(\eta)\times A(\xi_0)$ is the automorphism group of $\xi$. We obtain in this way a functor $A:{\cal C}_B\to{\cal G}r$. In the terminology of section \ref{limits}, the pair $({\cal C}_B,A)$ is a  category with automorphism push-forward.

 Let  $X\to B$, $Y\to B$ be two pointed $S^1$-spaces over $B$.  The assignment
 $$\xi\mapsto {_{S^1}\pi}^0_B(X\wedge_B\xi^+_B,Y\wedge_B\xi^+_B)
 $$
(where ${_{S^1}\pi}^0_B(X,Y)$ stands for the set of  homotopy classes  of $S^1$-equivariant base point preserving maps over $B$) is functorial with respect to morphisms in ${\cal C}_B$: for a morphism $u=(i,\zeta):\xi\to \xi'$, the push-forward class $u_*([f])$ is defined  using $i\circ f\circ i^{-1}$ on $i(\xi)$ and $\id_\zeta$ on its complement $\zeta$. Therefore this assignment defines   a functor ${_{S^1}\pi}^0_B(X\wedge_B\cdot ,Y\wedge_B\cdot):{\cal C}_B\to{\cal S}ets$. It is not clear at all that an inductive limit of this functor exists, because ${\cal O}b({\cal C}_B)$ is neither a filtering nor a small category (see section \ref{limits}).  
\begin{pr}\label{TSA} Let $\xi=\eta\oplus\xi_0\in {\cal O}b({\cal C}_B)$, $\ag=(\alpha,a_0)\in A(\xi)$, and $u=(i,\zeta)$ the standard morphism $\eta\oplus\xi_0=\xi\to \tilde\xi:=(\eta\oplus\eta)\oplus (\xi_0\oplus\xi_0)$ defined by  $(y,x)\mapsto ((y,0),(x,0))$. For every $[f]\in {_{S^1}\pi}^0_B(X\wedge_B\xi^+_B,Y\wedge_B\xi^+_B)$ one has
$$u_*(\ag_*([f])=u_*([f])\ .
$$
\end{pr}
\pf  Identifying $\tilde \xi$ with $\xi\oplus\xi$ one can write $u_*(\ag_*[f])=[g]$ where $g$ is the composition
$$(\id_X\wedge_B [\ag\oplus\id_\xi]^+_B)\circ (f\wedge_B\id_{\xi^+_B})\circ (\id_X\wedge_B [\ag^{-1}\oplus\id_\xi]^+_B): X\wedge_B[\xi\oplus\xi]^+_B\to Y\wedge_B[\xi\oplus\xi]^+_B. $$
Let $R_t$ be the automorphism of $\xi\oplus\xi$ defined by the matrix
$$r_t:=\left(
\begin{array}{cr}
\cos (t\frac{\pi}{2})&-\sin (t\frac{\pi}{2})\\ \sin (t\frac{\pi}{2})&\cos (t\frac{\pi}{2})
\end{array}\right)\ .
$$
For an automorphism $\bg$ of $\xi$ note that $r_t\circ (\bg\oplus \id_\xi)\circ r_t^{-1}$ defines a homotopy between $\bg\oplus \id_\xi$ and $\id_\xi\oplus \bg$. This shows that $g$ is homotopic to the map
$$g':=(\id_X\wedge_B [\id_\xi\oplus \ag]^+_B)\circ (f\wedge_B\id_{\xi^+_B})\circ (\id_X\wedge_B [\id_\xi\oplus  \ag^{-1}]^+_B)=f\wedge_B \id_{\xi^+_B}
$$
which is a representative of the class $u_*([f])$.
\qed
\vspace{2mm}

We define the stable cohomotopy group ${_{S^1}\alpha}^0_B(X,Y)$ by 

$${_{S^1}\alpha}^0_B(X,Y):=\varinjlim\hskip-22pt\raisebox{-11pt}{$\scriptstyle(n,m)\in\N^2$\ \ }\ {_{S^1}\pi}^0_B(X\wedge_B[\underline{\C}^n\oplus\underline{\R}^m]^+_B,Y\wedge_B[\underline{\C}^n\oplus\underline{\R}^m]^+_B) \ .$$
In this formula and in the rest of the paper we use the notation $\underline{V}$ for the trivial bundle $B\times V$ over the base $B$. This inductive limit has a natural Abelian group structure (see \cite{CJ} p. 168 for the non-equivariant case). 
\begin{pr}\label{change} The functor ${_{S^1}\pi}^0_B(X\wedge_B\cdot ,Y\wedge_B\cdot):{\cal C}_B\to{\cal S}ets$ admits an inductive limit, which can be identified with ${_{S^1}\alpha}^0_B(X,Y)$.
\end{pr}
\pf  Let ${\cal N}^2$ be the   small category associated with the ordered set $(\N\times\N,\leq)$ and consider the functor $\Theta:{\cal N}^2\to {\cal C}_B$ which assigns to a pair $(n,m)$ the trivial bundle $\underline{\C}^n\oplus\underline{\R}^m$ over $B$, and to an inequality $(n,m)\leq (n',m')$ the standard morphism between the corresponding trivial bundles. Using the terminology of section \ref{limits}, ${\cal N}$ is a small filtering category, and $\Theta$ is a cofinal functor from ${\cal N}$ to  the   category $({\cal C}_B,A)$, which is a category with automorphism push-forward. By definition  ${_{S^1}\alpha}^0_B(X,Y)$ is just the limit of the composed functor ${_{S^1}\pi}^0_B(X\wedge_B\cdot ,Y\wedge_B\cdot)\circ\Theta$. On the other hand, Proposition \ref{TSA} shows that the functor ${_{S^1}\pi}^0_B(X\wedge_B\cdot ,Y\wedge_B\cdot)$ satisfies the ``trivial stable actions" axioms TSA, $\Theta$SA. The result follows therefore from Proposition \ref{trivact} in section \ref{limits}.
 \qed
 
 Note that  Proposition \ref{change} implicitly yields a canonical map 
 $$c_\xi:{_{S^1}\pi}^0_B(X\wedge_B\xi^+_B ,Y\wedge_B\xi^+_B)\to {_{S^1}\alpha}^0_B(X,Y)$$
  for every $\xi\in{\cal O}({\cal C}_B)$, such that the system $(c_\xi)_{\xi\in{\cal O}({\cal C}_B)}$ satisfies the universal property of the inductive limit.\\

As in the non-equivariant case we put
$${_{S^1}\alpha}^p_B(X,Y):={_{S^1}\alpha}^0_B( X\wedge_B (\underline{\R}^N)^+_B, Y\wedge_B (\underline{\R}^{N+p})^+_B)\ (N,\ N+p\geq 0)\ .$$
Each ${_{S^1}\alpha}^p_B(X,Y)$ is a bimodule over the ring
$${_{S^1}\alpha}^0(B):={_{S^1}\alpha}^0(B_+,S^0)={_{S^1}\alpha}^0_B(B_{+B},B_{+B})\ ,
$$
and ${_{S^1}\alpha}^*_B(X,Y):=\oplus_{p\in\Z}\ {_{S^1}\alpha}^p_B(X,Y)$ is a graded bimodule over the graded ring ${_{S^1}\alpha}^*(B)=\oplus {_{S^1}\alpha}^p(B)$, where
$${_{S^1}\alpha}^p(B):={_{S^1}\alpha}^p(B_+,S^0)={_{S^1}\alpha}^0(B_+,S^p)\ .
$$
Right and left multiplication  with elements in ${_{S^1}\alpha}^0(B)$ coincide (see \cite{CJ} p. 172).   
\begin{re}  In the special case when $Y$ is of the form $Y=\zeta^+_B$ with $\zeta\in{\cal C}_B$, one has a canonical identification
$${_{S^1}\alpha}^0_B(X,\zeta^+_B)={_{S^1}\alpha}^0\left(\qmod{X\wedge_B[\zeta']^+_B}{\infty}\ ,\ V^+\right)\ , 
$$
where $\zeta\oplus\zeta'=\underline{V}$, and $V$ has the form $\C^k\oplus\R^l$. In the terminology of \cite{BF} the latter group is a stable cohomotopy group formed with respect to the universum generated by the $S^1$-representations $\C$ and $\R$.
\end{re}

\subsection{The computation of ${_{S^1}\alpha}^k(B_+,V^+)$}\label{comp}

Let $S^1\to O(V)$ be an orthogonal representation of $S^1$.  Our next goal is the computation of the group  ${_{S^1}\alpha}^k(B_+,V^+)$ for $k\geq 0$. In particular, we obtain explicit descriptions of the positive summands  ${_{S^1}\alpha}^k(B)={_{S^1}\alpha}^k(B_+,[\R^k]^+)$ of the graded ring ${_{S^1}\alpha}^*(B)$.

Replacing $V$ by $V\oplus \R^k$, we can reduce the problem to the case $k=0$. One has 
$${_{S^1}\alpha}^0(B_+,V^+)=\varinjlim\hskip-21pt\raisebox{-11pt}{$\scriptstyle(n,m)\in\N^2$}\left[B_+\wedge[\C^n\oplus\R^m]^+,V^+\wedge[\C^n\oplus\R^m]^+\right]^{S^1}_0\ ,
$$
where   $[\cdot,\cdot ]^{S^1}_0$ stands for the set of homotopy classes of $S^1$-equivariant maps between two pointed $S^1$-spaces. \\

 According to Hauschild's splitting theorem (Satz 3.4 in \cite{H}) there is a natural identification
\begin{equation}\label{split}
\left[B_+\wedge[\C^n\oplus\R^m]^+,[V\oplus\C^n\oplus\R^m]^+\right]^{S^1}_0=\ \ \ \ \ \ \ \ \ \ \ \ \ \ \ \ \ \ \ \ \ \ \ \ \ \ \ \ \ \ \ \ \ \ \ \ \ \ \ \ \ \ \ \ \ \ \ \ \ \ \ \ \ \ \ \ \ \ \ \ \ \ \ \ \ \ \end{equation}
$$  \left[B_+\wedge[\R^m]^+, [V^{S^1}]^+\wedge[\R^m]^+\right]_0\times \left[B_+\wedge\left[\qmod{[\C^n\oplus\R^m]^+}{[\R^m]^+}\right],V^+\wedge[\C^n\oplus\R^m]^+\right]^{S^1}_0
$$
where the projection on the first factor is given by restriction to the fixed point set. There exists a homeomorphism  of $S^1$-spaces
$$
\qmod{[\C^n\oplus\R^m]^+}{[\R^m]^+}\approx S(\C^n)_+\wedge S^{m+1}\ .
$$
Indeed, one has
$$
\qmod{[\C^n\oplus\R^m]^+}{[\R^m]^+}\approx \qmod{S(\C^n\oplus \R^{m+1})}{S(\R^{m+1})}\approx $$
$$\approx\qmod{S(\C^n)\times D(\R^{m+1})\cup D(\C^n)\times S(\R^{m+1})}{D(\C^n)\times S(\R^{m+1})}\approx
S(\C^n)_+\wedge S^{m+1}\ .$$

Using the natural identification
$$B_+\wedge[S(\C^n)_+\wedge S^{m+1}]\approx S(\C^n)_+\wedge[B_+\wedge S^{m+1}]\approx \qmod{S(\C^n)\times [B_+\wedge S^{m+1}]}{S(\C^n)\times\{*\}}
$$
and denoting by $\tilde V_n$ the associated bundle $S(\C^n)\times_{S^1}V$ over $\P(\C^n)$ we find
$$ \left[B_+\wedge\left[\qmod{[\C^n\oplus\R^m]^+}{[\R^m]^+}\right],V^+\wedge[\C^n\oplus\R^m]^+\right]^{S^1}_0\cong  \  \ \ \ \ \ \ \ \ \ \ \ \ \ \ \ \ \ \ \ \ \ \ \ \ \ \ \ \ \ \ \ \  \ \ \ $$
$$ \ \ \  \ \ \ \ \ \ \ \ \ \ \ \ \ \ \ \   \ \ \   \cong \left[\qmod{S(\C^n)\times [B_+\wedge S^{m+1}]}{S(\C^n)\times\{*\}},V^+\wedge[\C^n\oplus\R^m]^+\right]^{S^1}_0\cong 
$$
$$ \ \ \ \ \ \ \ \ \ \ \ \ \ \ \ \ \ \ \ \ \ \ \ \ \  \ \ \  \ \ \  \ \  \cong  {_{S^1}\pi}^0_{S(\C^n)}\left(S(\C^n)\times [B_+\wedge S^{m+1}],S(\C^n)\times[V\oplus\C^n\oplus\R^m]^+\right)\cong\ \ \ \ \ \ \ \ \ \ \ \ 
$$
$$ \ \ \ \ \ \ \ \ \ \ \ \ \ \ \  \ \ \ \ \ \  \ \ \  \ \ \  \ \ \   \cong  {\pi}^0_{\P(\C^n)}\left(\P(\C^n)\times [B_+\wedge S^{m+1}],[\tilde V_n\oplus {\cal O}_{\P(\C^n)}(1)^{\oplus n}\oplus\underline{\R}^{m}]^+_{\P(\C^n)}\right)\cong\ \ \ \ \ \ \ \ \ \ $$
$$   \cong  {\pi}^0_{\P(\C^n)}\left(\left[\P(\C^n)\times [B_+\wedge S^1]\right]\wedge_{\P(\C^n)} \underline{S}^{m}],  [ \tilde V_n\oplus{\cal O}_{\P(\C^n)}(1)^{\oplus n}]^+_{\P(\C^n)}\wedge_{\P(\C^n)}\underline{S}^m\right) .$$

The limit over $m$ of this set is 
$$\omega^0_{\P(\C^n)}\left(\P(\C^n)\times [B_+\wedge S^{1}],[\tilde V_n\oplus {\cal O}_{\P(\C^n)}(1)^{\oplus n}]^+_{\P(\C^n)}\right) \ .
$$
Now note that 
$$\tilde V_n\oplus \underline{\C}\oplus T_{\P(\C^n)}\cong  \tilde V_n\oplus {\cal O}_{\P(\C^n)}(1)^{\oplus n}\ .
$$
Therefore, applying the duality isomorphism given in Proposition 12.41 \cite{CJ} to the map $\pi:\P(\C^n)\to \{*\}$, one gets
$$\omega^0_{\P(\C^n)}\left(\P(\C^n)\times [B_+\wedge S^{1}],[\tilde V_n\oplus {\cal O}_{\P(\C^n)}(1)^{\oplus n}]^+_{\P(\C^n)}\right)\cong \omega^0(B_+\wedge S^1,\pi_*([\tilde V_n\oplus \underline{\C}]^+_{\P(\C^n)})) 
$$
$$\cong \omega^0(B_+\wedge S^1,   T({\tilde V_n\oplus \underline{\C}}) )\cong\omega^0(B_+\wedge S^1 ,  T({\tilde V}_n)\wedge S^2 )\cong \omega^0(B_+,  T({\tilde V}_n)\wedge S^1 )\ ,
$$
where $T(\cdot)$ stands for the  Thom space functor. This shows that
$$
\varinjlim\hskip-26pt\raisebox{-11pt}{$\scriptstyle(n,m)\in\N^2$}\hskip-2pt\left[B_+\wedge\left[\qmod{[\C^n\oplus\R^m]^+}{[\R^m]^+}\right],V^+\wedge[\C^n\oplus\R^m]^+\right]^{S^1}_0\hskip-8pt
\cong\omega^0(B_+,  T(ES^1\times_{S^1} V)\wedge S^1 ) $$
where $ES^1\times_{S^1} V$ is the vector bundle associated with the universal $S^1$-bundle $ES^1\to BS^1=\P^\infty$ and the  fiber $V$. Using formula  (\ref{split}) we obtain the following
\begin{pr}\label{fs} One has a natural group isomorphism
\begin{equation}\label{product}
{_{S^1}\alpha}^0(B_+,V^+)\cong \omega^0(B_+,[V^{S^1}]^+)\times \omega^0(B_+,  T(ES^1\times_{S^1} V)\wedge S^1 ) \ .
\end{equation}
where the projection on the first factor is given by restriction to the fixed point set. In particular
$${_{S^1}\alpha}^k(B)\cong \omega^k(B)\times \omega^k(B_+,  \P^\infty_+\wedge S^1)\ .
$$
\end{pr}
\begin{re} \label{ideal} The second summand in the decomposition 
$${_{S^1}\alpha}^0(B)\cong \omega^0(B)\times \omega^0(B_+,  \P^\infty_+\wedge S^1)$$
 is called  ``the free summand" in \cite{CK}. The projection ${_{S^1}\alpha}^0(B)\to \omega^0(B)$ is given by restriction to the fixed point set, hence it is a ring homomorphism. Therefore the free summand $\omega^0(B_+,\P^\infty_+\wedge S^1)$ is an ideal of ${_{S^1}\alpha}^0(B)$, and one has a natural ring isomorphism
$$\omega^0(B)\simeq\qmod{{_{S^1}\alpha}^0(B)}{\omega^0(B_+,\P^\infty_+\wedge S^1)}\ .
$$
\end{re}

\begin{co} \label{indlimN} Suppose that $B$ is a finite CW complex. Restriction to the fixed point set defines an isomorphism
$$\varinjlim\hskip-16pt\raisebox{-11pt}{$\scriptstyle N\in\N$}\ {_{S^1}\alpha}^k(B_+,[\C^N]^+)\textmap{\cong} \omega^k(B)\ .
$$
\end{co}
\pf Indeed, taking $V=\C^N\oplus\R^k$, the second summand in (\ref{product}) is:
$$\omega^0(B_+,  T(ES^1\times_{S^1} [\C^N\oplus \R^{k+1}]))=\varinjlim\hskip-16pt\raisebox{-11pt}{$\scriptstyle l\in\N$}\ \pi^0(B_+\wedge[\R^l]^+,  T(ES^1\times_{S^1} [\C^N\oplus \R^{k+1+l}]))$$

Recall that the Thom space of a real vector bundle of rank $r$ over a CW complex $X$ admits a CW decomposition consisting of a single 0-dimensional cell and cells of dimension $\geq r$. Therefore, for $N$ sufficiently large any map $B_+\wedge[\R^l]^+\to  T(ES^1\times_{S^1} [\C^N\oplus \R^{k+1+l}])$  is homotopically trivial.
\qed

\subsection{The groups ${\alpha}^*(x)$ associated with an element $x\in K(B)$}\label{alpha}

Fix an element $x\in K(B)$.  We define a category ${\cal T}(x)$ in the following way: the objects of ${\cal T}(x)$ are the presentations   of $x$. For two such presentations $(E,F)$, $(E',F')$, a morphism $\tau:(E,F)\to (E',F')$ is a system $\tau=(i,j,E_1,F_1,k)$ consisting of bundle monomorphisms $j:E\hookrightarrow E'$, $i:F\hookrightarrow F'$, complements $E_1$ and $F_1$ of $i(E)$ and $j(F)$ in $E'$ and $F'$ respectively, and an isomorphism $k:E_1\to F_1$.

With every   $(E,F)\in x$ we associate the graded group ${_{S^1}\alpha}_B^*(S(E)_{+B},F^+_B)$.  In this formula the sphere bundle $S(E)$ is defined by $S(E):=(E\setminus 0_E)/{\R_{>0}}$; alternatively one can use   an arbitrary  Hermitian metric on $E$. We claim that a morphism $\tau:(E,F)\to (E',F')$ induces a morphism
$$
\tau_*: {_{S^1}\alpha}_B^*(S(E)_{+B},F^+_B)\map {_{S^1}\alpha}_B^*(S(E')_{+B},[F']^+_B)\ .
$$

Note first that, for Euclidean or Hermitian vector spaces $V$, $W$, one has a contraction
$$S(V\oplus W)\to S(V)_+\wedge W^+
$$
induced by the map
$$c:S(V\oplus W)=[S(V)\times D(W)]\cup_{S(V)\times S(W)} [D(V)\times S(W)]\map$$
$$\map \qmod{S(V)\times D(W)}{S(V)\times S(W)}\simeq\qmod{S(V)\times W^+}{S(V)\times \infty_W}=S(V)_+\wedge W^+\ .
$$
It is useful  to have  explicit analytic formulae for the contraction map $c$. One can define $W^+$ in two equivalent ways: as the one-point compactification of $W$, and as the quotient $D(W)/S(W)$. Accordingly, the contraction maps $c$, $c':S(V\oplus W)\to S(V)_+\wedge W^+$ are given by the formulae:
\begin{equation}\label{contraction}
c(v,w)=\left\{\hskip -8pt
\begin{array}{ll}
\left(\frac{1}{\|v\|}v,\frac{1}{\sqrt{1-\|w\|^2}}w\right)&v\ne 0\\
*&v= 0 
\end{array}\right. \ ,\  c'(v,w)=\left\{\hskip -8pt
\begin{array}{ccc}
\left(\frac{1}{\|v\|}v, w\right)&v\ne0\ \\
*&v=0. 
\end{array}\right.
\end{equation}
To save on notations we will still write $c$  instead of $c'$ when the second definition  of $W^+$ is used.

In the presence of a morphism $\tau=(i,j,E_1,F_1,k):(E,F)\to (E',F')$ we choosing Hermitian metrics on $E'$ and $F'$ which make the isomorphisms $i$, $j$, $k$ isometries and the decompositions $E'=i(E)\oplus E_1$, $F'=j(F)\oplus F_1$ orthogonal. We get a  map
$$S(E')_{+B}=S(i(E)\oplus E_1)_{+B}\textmap{c} S(i(E))_{+B}\wedge_B (E_1)^+_B\ ,
$$
which is well defined up to homotopy (the section $+_B$ on the left is mapped fiberwise to the distinguished section on the right). One obtains morphisms
$${_{S^1}\alpha}_B^*(S(E)_{+B},F^+_B) \textmap{(i,j)\simeq}{_{S^1}\alpha}_B^*(S(i(E))_{+B},j(F)^+_B)=
$$
$$={_{S^1}\alpha}_B^*(S(i(E))_{+B}\wedge_B (F_1)^+_B,j(F)^+_B\wedge_B(F_1)^+_B)={_{S^1}\alpha}_B^*(S(i(E))_{+B}\wedge_B (F_1)^+_B,(F')^+_B) $$
$$\stackrel{k}{\simeq}{_{S^1}\alpha}_B^*(S(i(E))_{+B}\wedge_B (E_1)^+_B,(F')^+_B)\textmap{c^*}{_{S^1}\alpha}_B^*(S(E')_ {+B},(F')^+_B)\ .
$$

The composition of these maps will be denoted by $\tau_*$.  One checks that $\tau_*$ is a morphism of ${_{S^1}\alpha}^*(B)$-modules and that, for any two composable morphisms $\tau$, $\tau'$, one has
$$(\tau'\circ\tau)_*=\tau'_*\circ \tau_*\ .
$$
In other words, the assignment $(E,F)\mapsto{_{S^1}\alpha}_B^*(S(E)_{+B},F^+_B)$ is functorial, so it defines a functor $\ag^*_x:{\cal T}(x)\to {\cal A}b^*$, where ${\cal A}b^*$ is the category of graded Abelian groups.\\ \\
{\bf Example:} Suppose that the stable class $\varphi\in {_{S^1}\alpha}_B^0(S(E)_{+B},F^+_B)$ is represented by an $S^1$-equivariant map $f:S(E)\to F^+_B$ over $B$ (or, equivalently, by an $S^1$-equivariant map $S(E)_{+B}\to F^+_B$ of pointed spaces over $B$). Let $U$ be  a complex vector bundle over $B$ and let $\tau$ be the obvious morphism $(E,F)\to (E\oplus U, F\oplus U)$. Then $f$ defines a map
$$ \qmod{[S(E)\times_B U^+_B]}{S(E)\times_B\infty_U}\map \qmod{F^+_B\times_B  U^+_B}{F^+_B\times_B\infty_U}
$$
which, composed from the right with the contraction 
$$S(E\oplus U)\to \qmod{[S(E)\times_B U^+_B]}{S(E)\times_B\infty_U}$$
and from on left  with the contraction
$$\qmod{F^+_B\times_B  U^+_B}{F^+_B\times_B\infty_U}\to \qmod{F^+_B\times_B  U^+_B}{\left[F^+_B\times_B\infty_U\cup\infty_F\times_B U^+_B\right]}=(F\oplus U)^+_B
$$
gives an $S^1$-equivariant map $S(E\oplus U)\to (F\oplus U)^+_B$ over $B$. This map represents $\tau_*(\varphi)\in {_{S^1}\alpha}_B^0(S(E\oplus U)_{+B},(F\oplus U)^+_B)$.\\

  Let $a\in \Aut(E)$ be a unitary gauge transformation of the bundle $E$. Composing with the induced automorphisms  $S(a)$  of the sphere bundles $S(E)_{+B}$ defines  a morphism 
$${_{S^1}\alpha}_B^*(S(E)_{+B},F^+_B)\textmap{S(a)^*}  {_{S^1}\alpha}_B^*(S(E)_{+B},F^+_B)\ .
$$
On the other hand, $a$ defines an element $[a^+_B]\in {_{S^1}\pi}^0_B(E^+_B,E^+_B)$, whose stable class $\{a^+_B\}$ is a  unit in the ground ring ${_{S^1}\alpha}^0(B)$ and defines  multiplication automorphisms
$${_{S^1}\alpha}_B^*(S(E)_{+B},F^+_B)\textmap{m(a)}  {_{S^1}\alpha}_B^*(S(E)_{+B},F^+_B) \ .
$$
Clearly these automorphisms depend only on the homotopy class of $a$. 

\begin{pr} \label{gaugeaction}  Let   $\varphi\in {_{S^1}\alpha}^*(S(E)_{+B},F^+_B)$ and  $a\in\Aut(E)$.  Let  $\tau$ be  the obvious morphism $\tau:(E,F)\to (E\oplus E,F\oplus E)$. In ${_{S^1}\alpha}^*(S(E\oplus E)_{+B},[F\oplus E]^+_B)$ is holds
$$\tau_*(S(a)^*(\varphi))=\tau_*(m(a)(\varphi))\ .$$
\end{pr}
\pf For simplicity we prove the statement only in degree 0. We may assume that $a$ is a unitary automorphism with respect to a Hermitian structure on $E$. Suppose that $\varphi$ is represented by 
$$[f]\in {_{S^1}\pi}^0_B(S(E)_{+B}\wedge_B\xi^+_B, F^+_B\wedge_B \xi^+_B)\ .
$$
 We will prove  that the natural representatives 
 $$p,\ q\in {_{S^1}{\rm Map}}_B(S(E\oplus E)_{+B}\wedge_B\xi^+_B\wedge_B E^+_B, (F\oplus E)^+_B\wedge_B\xi^+_B\wedge E^+_B)$$ of  $\tau_*(S(a)^*([f]))$, $\tau_*(m(a)([f]))$ are homotopic, so they define the same element in
$${_{S^1}\pi}^0_B(S(E\oplus E)_{+B}\wedge_B\xi^+_B\wedge_B E^+_B, (F\oplus E)^+_B\wedge_B\xi^+_B\wedge_B E^+_B)\ . $$
We suppose for simplicity that $\xi$ is trivial, to save on notations. Consider the contraction map $c:S(E\oplus E)_{+B}\to S(E)_{+B}\wedge_B E^+_B$ defined by the first formula in (\ref{contraction}), and introduce the maps
$$\Psi,\ \chi: S(E)_{+B}\wedge_B E^+_B \wedge_B E^+_B\map  F^+_B\wedge_B E^+_B\wedge_B E^+_B
$$
defined by
$$\Psi:=[ f\circ S(a)]\wedge_B \id_{E^+_B}\wedge_B\id_{E^+_B}\ ,\ \chi:=f\wedge_B \id_{E^+_B}\wedge_B a^+_B\ .
$$
Using our definitions it is easy to see  that $p=\Psi\circ (c\wedge_B\id_{E^+_B})$, $q=\chi\circ (c\wedge_B\id_{E^+_B})$.  Use the same method as in the proof of Proposition \ref{TSA} (conjugation with the rotations of $E\oplus E$ defined by the matrices $r_t$) to construct a homotopy 
$$\chi=f\wedge_B(\id_E\oplus a)^+_B\simeq f\wedge_B(a\oplus \id_E)^+_B=f\wedge_B a^+_B\wedge_B \id_{E^+_B}:=\chi'\ .$$
It suffices to construct a homotopy between $\Psi\circ (c\wedge_B\id_{E^+_B})$, and  $\chi'\circ (c\wedge_B\id_{E^+_B})$, and for this it suffices to construct a homotopy between the maps $\Psi_0\circ c$ and $\chi'_0\circ c$, where
$$\Psi_0:=[ f\circ S(a)]\wedge_B \id_{E^+_B}=(f\wedge_B\id_{E^+_B})\circ (S(a)\wedge_B\id_{E^+_B})\ ,$$
$$\chi'_0:=f\wedge_B a^+_B=(f\wedge_B\id_{E^+_B})\circ(\id_{E^+_B}\wedge_B a^+_B)\ .$$
Note that  $(S(a)\wedge_B\id_{E^+_B})\circ c=c\circ S(a\oplus \id_E)$, and $(\id_{S(E)}\wedge_B a^+_B)\circ c=c\circ S(\id_E\oplus a)$. In these formulae we use the fact that $a$ is a unitary. On the other hand, using again conjugation with the rotations defined be the matrices $r_t$,  we see that $S(a\oplus \id)\simeq S(\id_E\oplus a)$. Therefore
$$\Psi_0\circ c= (f\wedge_B\id_{E^+_B})\circ c\circ S(a\oplus \id)\simeq (f\wedge_B\id_{E^+_B})\circ c\circ S(\id_E\oplus a)=$$
$$=(f\wedge_B\id_{E^+_B})\circ (\id_{S(E)}\wedge_B a^+_B)\circ c=\chi'_0\ ,
$$
which completes the proof.
\qed
\vspace{2mm}

A similar statement holds for the action of an automorphism $b\in \Aut(F)$. Denote by $[b^+_B]_*$ the automorphism  of ${_{S^1}\alpha}^*(S(E)_{+B},F^+_B)$ defined by   composition with  $b^+_B$.
 \begin{pr} \label{gaugeF}  The automorphisms   $[b^+_B]_*$,  $m(b)$  coincide on  ${_{S^1}\alpha}^*(S(E)_{+B}, F^+_B)$. %
 \end{pr}
 The proof uses similar arguments as the proof of Proposition \ref{gaugeaction} but is substantially easier.
\qed

An automorphism  $c\in\Aut(U)$ defines a automorphism $\sigma(c)$ of the graded group $\alpha^*(S(E\oplus U)_{+B},[F\oplus U]^+_B)$  defined by  $f\mapsto [\id_F\oplus c]^+_B\circ f\circ S(\id_E\oplus c)^{-1}$.
\begin{co} \label{condition} Let $\tau:(E\oplus U,F\oplus U)\to (E\oplus U\oplus E\oplus U,F\oplus U\oplus E\oplus U)$ be the natural morphism. Then for any $\varphi\in \alpha^*(S(E\oplus U)_{+B},[F\oplus U]^+_B)$ one has
$$\tau_*(\sigma(c)(\varphi))=\tau_*(\varphi)\ .
$$
\end{co}
\pf Indeed, one has
$$\tau_*\circ \{[\id_F\oplus c]^+_B\}_*=\tau_*\circ m(c)\ ,\ \tau_*\circ \{S(\id_E\oplus c)^{-1}\}^*= \tau_*\circ(m(c)^{-1})\ .$$
On the other hand the morphism $\tau_*$ is ${_{S^1}\alpha}^0(B)$-linear.
\qed

Consider now the category  ${\cal U}_B$ of all finite rank complex vector bundles over $B$.  A morphism $\nu:U\to U'$ in the category ${\cal U}_B$ is a pair $(i,U_1)$ consisting of a bundle embedding $i:U\to U'$ and a complement $U_1$ of $i(U)$ is $U'$. This category can be regarded in an obvious way as a category with automorphism push-forward (see section \ref{limits}). The assignment $U\mapsto {_{S^1}\alpha}_B^*(S(E\oplus U)_{+B},(F\oplus U)^+_B)$ is functorial with respect to morphisms in ${\cal U}_B$, so it  defines a functor $\ag_{E,F}^*:{\cal U}_B\to {\cal A}b^*$. Since ${\cal U}_B$ is not a small category, it is not clear whether this functor has an inductive limit (see sections \ref{definition}, \ref{limits}).
We put
\begin{equation}\label{alphahat}
\hat\alpha^*(E,F):=\varinjlim\hskip-16pt\raisebox{-10pt}{$\scriptstyle n\in\N$\ \ }{_{S^1}\alpha}_B^*(S(E\oplus \underline{\C}^n)_{+B},(F\oplus \underline{\C}^n)^+_B)\ .
\end{equation}
\begin{pr}\label{U} The functor $\ag_{E,F}^*$ admits an inductive limit  which can be identified with $\hat\alpha^*(E,F)$.
\end{pr}
\pf  Let ${\cal N}$ be the category associated with the ordered set $(\N,\leq)$ and let $\Theta:{\cal N}\to {\cal U}_B$ be  the cofinal functor $n\mapsto \underline{\C}^n$ (see section \ref{limits}).  By Corollary \ref{condition}, the functor $\ag_{E,F}^*$ satisfies the trivial stable action axiom $\Theta$SA. The result follows now from Proposition \ref{trivact} in section \ref{limits}.
\qed
\vspace{2mm}\\
In particular one has canonical morphisms $c_U:{_{S^1}\alpha}_B^*(S(E\oplus U)_{+B},(F\oplus U)^+_B)\to \hat\alpha^*(E,F)$ for any complex bundle $U$, and the system $(c_U)_U$ is $\ag_{E,F}^*$-compatible and satisfies the universal property of the inductive limit. Note that $\hat\alpha^*(E,F)$ has a natural structure of a graded  ${_{S^1}\alpha}^*(B)$ bimodule.    By Propositions  
\ref{gaugeaction} and \ref{gaugeF} we get:
\begin{re}\label{viaunits} The action of the gauge groups $\Aut(E\oplus U)$, $\Aut(F\oplus U)$ on $\hat\alpha^*(E,F)$ is induced by the canonical ${_{S^1}\alpha}^0(B)^\times$-action  defined by its module structure via the morphisms $\Aut(E\oplus U)\to {_{S^1}\alpha}^0(B)^\times$, $\Aut(F\oplus U)\to{_{S^1}\alpha}^0(B)^\times$ defined by $a\mapsto a^+_B$.
\end{re}

 A morphism $\tau=(i,j,E_1,F_1,k):(E,F)\to (E',F')$  between two presentations $(E,F)$, $(E',F')$ of $x$ induces a sequence a morphisms $(E\oplus\underline{\C}^n,F\oplus\underline{\C}^n)\to (E'\oplus\underline{\C}^n,F'\oplus\underline{\C}^n)$, so it induces a  morphism   $\hat\tau_*:\hat\alpha^*(E,F)\textmap{\simeq} \hat\alpha^*(E',F')$. It is easy to see that $\hat\tau_*$ is an isomorphism: it suffices to note that the exists an isomorphism $\theta: (E',F')\to (E\oplus U, F\oplus U)$ (with $U:=E_1$) such that  $\theta\circ\tau$ is the standard morphism $(E,F)\to (E\oplus U,F\oplus U)$, and to apply Proposition \ref{U}. Therefore we obtain a functor $\hat\ag^*_x:{\cal T}(x)\to{\cal A}b^*$ whose associated morphisms $\hat\ag^*_x(\tau)=\hat\tau_*$ are all isomorphisms.
 According to Proposition \ref{iso} an  inductive limit  of this functor exists and can be identified with a quotient of $\hat\alpha^*(E,F)$, for any fixed presentation $(E,F)$ of $x$. Therefore we can make
\begin{dt}  Define
$$\alpha^*(x):=\varinjlim\hskip-24pt\raisebox{-10pt}{$\scriptstyle(E,F)\in x$}\hat\alpha^*(E,F)\ .
$$
\end{dt}
\begin{re} This inductive limit is also an inductive limit of the functor $\ag^*_x$ introduced at the beginning of this section. The existence of the inductive limit of this functor is a non-trivial statement. 
\end{re}
We introduce now the notations
$$\A(E):=\varinjlim\hskip-18pt\raisebox{-10pt}{$\scriptstyle N\in\N$}\ \Aut(E\oplus \C^N),\ \A(F):=\varinjlim\hskip-18pt\raisebox{-10pt}{$\scriptstyle N\in\N$} \Aut(F\oplus \C^N)\ .$$
The two groups $\A(E)$, $\A(F)$ act on the graded group $\hat\alpha^*(E,F)$ via the group morphisms $l:\A(E)\to {_{S^1}\alpha}^0(B)^\times$, $r:\A(F)\to {_{S^1}\alpha}^0(B)^\times$ (see Remark \ref{viaunits}), so the two actions commute.  Let $\Z[\A(E)]$, $\Z[\A(F)]$ be the group rings of $\A(E)$, $\A(F)$, $I[\A(E)]$, $I[\A(F)]$ the augmentation ideals, and $\lambda:\Z[\A(E)]\to {_{S^1}\alpha}^0(B)$, $\rho:\Z[\A(F)]\to {_{S^1}\alpha}^0(B)$ the ring morphisms associated with the group morphisms $l$, $r$.  Using Proposition \ref{iso} we get
\begin{re}\label{firstdesc} For every presentation $(E,F)\in x$ there is a canonical isomorphism
$$\alpha^*(x)\textmap{\simeq}\qmod{\hat\alpha^*(E,F))}{\lambda( I[\A(E)])\hat\alpha^*(E,F)+\rho(I[\A(F)])\hat\alpha^*(E,F)}\ .
$$
\end{re}
In the next section we will see that $\A(E)$, $\A(F)$ are both isomorphic to $K^{-1}(B)$ and we will identify the images $\lambda( I[\A(E)])$,  $\rho(I[\A(F)])$ of the two ideals in ${_{S^1}\alpha}^0(B)$ with the image of  the ideal $I[K^{-1}(B)]$ under  the ring morphism $\Z[K^{-1}(B)]\to {_{S^1}\alpha}^0(B)$ induced by the $J$-map $K^{-1}(B)\to {_{S^1}\alpha}^0(B)^\times$.

\subsection{The $S^1$-equivariant $J$-map and the description of $\alpha^*(x)$}\label{J}

Let $\pi:E\to B$ be a Hermitian vector bundle over a compact basis, and let $a$, $b\in\Aut(E)$ be two unitary automorphisms. We define a map
$$\Delta_E(a,b):S(E)_{+B}\wedge_B\underline{S}^1\map E^+_B
$$
in the following way: We use the models
$$S(E)_{+B}\wedge_B\underline{S}^1\cong \qmod{S(E)\times[0,1]}{S(E)\times\{0,1\}}\ ,\ E^+_B\cong \qmod{D(E)}{_B S(E)}$$
for the two spaces, and define  
$$\Delta_E(a,b)([e,t]):=\left\{
\begin{array}{ccc}
{\ } [(1-2t)a(e)]&{\rm for}& 0\leq t\leq \frac{1}{2}\ \ \\
{\ }[(2t-1)b(e)] &{\rm for}& \frac{1}{2}\leq t\leq 1\ .
\end{array}\right.
$$
Consider the contraction map
$$\cg_E: E^+_B\map S(E)_{+B}\wedge_B \underline{S}^1 
$$
induced by $e\mapsto [(\frac{1}{\|e\|} e,\|e\|)]$. One has
\begin{equation}
\{\Delta_E(a,b)\}=\{b^+_B\}-\{a^+_B\}\ .
\end{equation}
\begin{dt} The J-homomorphism associated with a Hermitian bundle $E$ is the morphism $J_E:\pi_0(\Aut(E))\to {_{S^1}\alpha}^0_B(B)^\times$ defined by $J_E([a]):=\{a^+_B\}$.
\end{dt}

We introduce the map
$$\Theta_E: \pi_0(\Aut(E))\map {_{S^1}\alpha}^{-1}_B(S(E)_{+B},E^+_B)\ ,\ \Theta_E([a]):=\{\Delta_E(\id_E,a)\}\ .
$$
Let $\partial_E:{_{S^1}\alpha}^{-1}_B(S(E)_{+B},E^+_B)\to {_{S^1}\alpha}^0_B(E^+_B,E^+_B)$ be the connecting morphism in the long exact cohomotopy sequence: 
\begin{equation}\label{gysin}
\dots\to  {_{S^1}\alpha}^{-1}_B(S(E)_{+B},E^+_B)\stackrel{\partial_E}{\to}{_{S^1}\alpha}^0_B(E^+_B,E^+_B)\to{_{S^1}\alpha}^0_B(B_{+B},E^+_B)\to\dots
\end{equation}
associated with   $E^+_B$ and the cofiber sequence
$$S(E)_{+B}\map D(E)_{+B}\map E^+_B\ .
$$
Since $\partial_E$ acts by composition with the contraction $\cg_E$, we see that the diagram
 \begin{equation}\label{diagram}
\begin{array}{c}
\unitlength=1mm
\begin{picture}(100,30)(15,0)
\put(25,25){$\pi_0(\Aut(E))$}
\put(48,26){\vector(2,0){22}}
\put(56,27.5){$\Theta_E$}
\put(75,25){${_{S^1}\alpha}^{-1}_B(S(E)_{+B},E^+_B)$}
\put(33,22){\vector(0,-2){10}}
\put(92,22){\vector(0,-2){10}}
\put(93,17){$\partial_E$}
\put(28,15){$J_E$}
\put(48,9){$\cdot-1$}
\put(25,7){$ {_{S^1}\alpha}^0(B)^\times
$}
\put(45,8){\vector(2,0){14}}
\put(62,7){$ {_{S^1}\alpha}^0(B)$}
\put(77,7){$=$}
\put(81,7){$ {_{S^1}\alpha}^0_B(E^+_B,E^+_B)$}
\end{picture} 
\end{array}
\end{equation} 
\vspace{-8mm}\\
is commutative.
\begin{re} Let $\omega^0(B_+,\P^\infty_+\wedge S^1)\subset {_{S^1}\alpha}^0(B)$ be the free summand of the ring  ${_{S^1}\alpha}^0(B)$ (see Proposition \ref{fs}). For any $[a]\in \pi_0(\Aut(E))$ it holds
$$J_E([a])-1\in \omega^0(B_+,\P^\infty_+\wedge S^1)\ .
$$
\end{re}
Indeed, $\omega^0(B_+,\P^\infty_+\wedge S^1)$ is the kernel of the morphism $\rho:{_{S^1}\alpha}^0(B)\to \omega^0(B)$ given by restriction to the fixed point set. Therefore
$$\rho(J_E([a]))=\rho(\{a^+_B\})=\{(a^+_B)^{S^1}\}=\{\id_{B_{+B}}\}\ .
$$
\qed
\begin{pr} One has
\begin{enumerate}
\item $$\varinjlim\hskip-12pt\raisebox{-10pt}{${\scriptstyle N}\ $}\ \pi_0(\Aut(E\oplus\underline{\C}^N))=K^{-1}(B)$$
\item The system of morphisms $(\partial_{E\oplus\underline{\C}^N})_{N\in\N}$ defines an isomorphisms
$$\partial:\varinjlim\hskip-12pt\raisebox{-10pt}{${\scriptstyle N}$}\ \ {_{S^1}\alpha}^{-1}_B(S(E\oplus\underline{\C}^N)_{+B},[E\oplus\underline{\C}^N]^+_B)\map \omega^0(B_+,\P^\infty_+\wedge S^1)\ .$$
\end{enumerate}
\end{pr}
\pf   Let $\Phi$ be a complex bundle on $B$. For any automorphism $a\in\Aut(\Phi)$ we construct a bundle $\Phi_a$  over $B\times S^1$ in the following way:  we consider the bundle  $\Phi\times [0,1]$  over $B\times[0,1]$ and we  identify $\Phi\times\{0\}$ with $\Phi\times\{1\}$ via $a$. This bundle comes with an obvious identification $\resto{\Phi_a}{B\times\{0\}}\simeq\resto{\p_B^*(\Phi)}{B\times\{0\}}$, so the difference $[\Phi_a]-[\p_B^*(\Phi)]$ defines an element  $k_\Phi(a)\in K(B\times S^1,B\times\{0\})$. It is easy to see that the obtained map $k_\Phi:\Aut(\Phi)\to K(B\times S^1,B\times\{0\})=K^{-1}(B)$   is a group morphism. Taking the   limit  over $N$ of the system of morphisms $k_{E\oplus\underline{\C}^N}$ we obtain a morphism
$$\kappa_E:\varinjlim\hskip-12pt\raisebox{-10pt}{${\scriptstyle N}\ $}\ \pi_0(\Aut(E\oplus \underline{\C}^N))\to K^{-1}(B)\ .
$$

Let $E'$ be a complement of $E$ and fix an isomorphism $E'\oplus E\cong \underline{\C}^n$. The assignment $a\mapsto \id_{E'}\oplus a$ defines an {\it injective} morphism
$$i_{E'}:\varinjlim\hskip-12pt\raisebox{-10pt}{${\scriptstyle N}\ $}\ \pi_0(\Aut(E\oplus \underline{\C}^N))\to \varinjlim\hskip-12pt\raisebox{-10pt}{${\scriptstyle N}\ $}\ \pi_0(\Aut(\underline{\C}^{n+N}))\ .
$$
Similarly, we obtain an obvious {\it injective} morphism 
$$j_E:\varinjlim\hskip-12pt\raisebox{-10pt}{${\scriptstyle N}\ $}\ \pi_0(\Aut(\underline{\C}^N))\to \varinjlim\hskip-12pt\raisebox{-10pt}{${\scriptstyle N}\ $}\ \pi_0(\Aut(E\oplus \underline{\C}^N))\ .$$
 Hence we have   morphisms
$$\varinjlim\hskip-12pt\raisebox{-10pt}{${\scriptstyle N}\ $}\ \pi_0(\Aut(\underline{\C}^N))\stackrel{j_E}{\to} \varinjlim\hskip-12pt\raisebox{-10pt}{${\scriptstyle N}\ $}\ \pi_0(\Aut(E\oplus \underline{\C}^N))\stackrel{i_{E'}}{\to}\varinjlim\hskip-12pt\raisebox{-10pt}{${\scriptstyle N}\ $}\ \pi_0(\Aut(\underline{\C}^{n+N}))\textmap{\kappa_{\underline{\C}^n}}  K^{-1}(B)\ .
$$
 The composition $i_{E'}\circ j_E$ is clearly an isomorphism. Moreover, it is well-known that $\kappa_{\underline{\C}^n}$ is an isomorphism, for every $n\in\N$.  Since $i_{E'}$ is injective, we see that $\kappa_E=\kappa_{\underline{\C}^n}\circ i_{E'}$ is injective. On the other hand, $\kappa_{\underline{\C}^n}\circ i_{E'}\circ j_E=\kappa_E\circ j_E$ is an isomorphism, so $\kappa_E$ is also surjective.

For the second isomorphism, we take the direct limit over $N$ in the  cohomotopy exact sequence (\ref{gysin}) associated with $E\oplus \underline{\C}^N$. We have
$$\varinjlim\hskip-12pt\raisebox{-10pt}{${\scriptstyle N}\ $}\  {_{S^1}\alpha}^k_B([E\oplus \underline{\C}^N]^+_B,[E\oplus \underline{\C}^N]^+_B)= {_{S^1}\alpha}^k(B)\ .
$$
On the other hand, the system of morphisms defined by restriction to the fixed point set (see section \ref{comp}) defines a morphism
$$r^k_E:\varinjlim\hskip-12pt\raisebox{-10pt}{${\scriptstyle N}\ $}\  {_{S^1}\alpha}^k_B(B_{+B},[E\oplus \underline{\C}^N]^+_B)\to \omega^k(B_+,S^0)=\omega^k(B)\ .
$$
Using again a complement $E'$ of $E$ as above, we obtain morphisms
$$\varinjlim\hskip-18pt\raisebox{-10pt}{$\scriptstyle N\in\N$}\ {_{S^1}\alpha}^k(B_{+}, [\C^N]^+)=\varinjlim\hskip-18pt\raisebox{-10pt}{$\scriptstyle N\in\N$}\ {_{S^1}\alpha}^k_B(B_{+B},B\times [\C^N]^+)\to \ \ \ \ \ \ \ \ \ \ \ \ \ \ \ \ \ \ \ \ \ \ \ \ \ \ \ \ \ \ \ \ \ $$
$$
\ \ \ \ \ \ \ \ \ \ \ \to \varinjlim\hskip-12pt\raisebox{-10pt}{${\scriptstyle N}\ $}\  {_{S^1}\alpha}^k_B(B_{+B},[E\oplus \underline{\C}^N]^+_B)\to \varinjlim\hskip-12pt\raisebox{-10pt}{${\scriptstyle N}\ $}\  {_{S^1}\alpha}^k_B(B_{+B},[\underline{\C}^{n+N}]^+_B)\textmap{r^k_{\underline{\C}^n}} \omega^k(B) \ .
$$
The morphism $\varinjlim\hskip-18pt\raisebox{-10pt}{$\scriptstyle N\in\N$}\ {_{S^1}\alpha}^k_B(B_{+B},B\times [\C^N]^+)\to\varinjlim\hskip-12pt\raisebox{-10pt}{${\scriptstyle N}\ $}\  {_{S^1}\alpha}^k_B(B_{+B},[\underline{\C}^{n+N}]^+_B)$ is an isomorphism, and   $\varinjlim\hskip-12pt\raisebox{-10pt}{${\scriptstyle N}\ $}\  {_{S^1}\alpha}^k_B(B_{+B},[E\oplus \underline{\C}^N]^+_B)\to \varinjlim\hskip-12pt\raisebox{-10pt}{${\scriptstyle N}\ $}\  {_{S^1}\alpha}^k_B(B_{+B},[\underline{\C}^{n+N}]^+_B)$ is injective. Moreover, by  Corollary \ref{indlimN}, the map  $r^k_{\underline{\C}^n}$ is an isomorphism. Now the same arguments as above show that    $r^k_E$ is an isomorphism.
 The limit of (\ref{gysin}) becomes
$${_{S^1}\alpha}^{-1}(B)\stackrel{\rho^{-1}}{\to}\omega^{-1}(B)\to \varinjlim\hskip-12pt\raisebox{-10pt}{${\scriptstyle N}\ $}\ {_{S^1}\alpha}^{-1}_B(S(E\oplus \underline{\C}^N)_{+B},[E\oplus \underline{\C}^N]^+_B)\stackrel{\partial}{\to} {_{S^1}\alpha}^{0}(B)\stackrel{\rho}{\to}\omega^0(B)
$$
But the map 
$$\rho^{-1}:{_{S^1}\alpha}^{-1}(B)={_{S^1}\alpha}^{0}(B_+\wedge S^1)\to \omega^0(B_+\wedge S^1)=\omega^{-1}(B)
$$
is also induced by restriction to the fixed point set, so it is surjective by Remark \ref{ideal} applied to the basis $B_+\wedge S^1$. Therefore $\partial$ induces an isomorphism
$$\varinjlim\hskip-12pt\raisebox{-10pt}{${\scriptstyle N}\ $}\ {_{S^1}\alpha}^{-1}_B(S(E\oplus \underline{\C}^N)_{+B},[E\oplus \underline{\C}^N]^+_B)\textmap{\cong}\ker(\rho)=\omega^0(B_+,\P^\infty_+\wedge S^1)\ . 
$$
\qed

Taking the inductive limit with respect to $N$ of the diagram (\ref{diagram}) written for $E\oplus \underline{\C}^N$, we obtain the commutative diagram
 \begin{equation}\label{diagramnew}
\begin{array}{l}
\unitlength=1mm
\begin{picture}(105,30)(25,0)
\put(25,25){$K^{-1}(B)$}
\put(40,26){\vector(2,0){12}}
\put(45,27){$\Theta$}
\put(55,25){$\varinjlim\hskip-12pt\raisebox{-10pt}{${\scriptstyle N}\ $}\ {_{S^1}\alpha}^{-1}_B(S(E\oplus \underline{\C}^N)_{+B},[E\oplus \underline{\C}^N]^+_B)=\hat\alpha^{-1}(E,E)$}
\put(33,22){\vector(0,-2){10}}
\put(80,22){\vector(0,-2){10}}
\put(81,17){$\partial$}
\put(76,17){$\simeq$}
\put(30,15){$J$}
\put(50,9){$\cdot-1$}
\put(25,7){$ {_{S^1}\alpha}^0(B)^\times
$}
\put(104,7){$ {_{S^1}\alpha}^0(B)\ .
$}
\put(42,8){\vector(2,0){25}}
\put(98,7){$\stackrel{\iota}{\hookrightarrow}$}
\put(70,7){$\omega^0(B_+,\P^\infty_+\wedge S^1)$}
\end{picture} 
\end{array}
\end{equation} 
\vspace{-10mm}\\
\begin{re} The map $\iota\circ\partial\circ\Theta:K^{-1}(B)\to  {_{S^1}\alpha}^0(B)$ satisfies the identity 
$$[\iota\circ\partial\circ\Theta](a+b)=[\iota\circ\partial\circ\Theta](a)[\iota\circ\partial\circ\Theta](b)+[\iota\circ\partial\circ\Theta](a)+[\iota\circ\partial\circ\Theta](b)\ .$$
It is the ``free J-map"  in the terminology of Crabb-Knapp (\cite{CK}, p. 88, p.93).  
\end{re}
\begin{co} The map $J:K^{-1}(B)\to  {_{S^1}\alpha}^0(B)^\times$ is injective.
\end{co}
\pf It suffices to note that $\partial\circ\Theta$ is injective by Corollary 2.5 in \cite{CK}. 
\qed
\\ \\
The group morphism $J$ extends to a ring morphism $\tilde J:\Z[K^{-1}(B)]\to  {_{S^1}\alpha}^0(B)$.
\\ \\
{\bf Question:} {\it Does the subgroup 
$$\tilde J(I[K^{-1}(B)])=\langle \{J(u)-1|\ u\in K^{-1}(B)\}\rangle =\langle \im(\partial\circ\Theta) \rangle\subset \omega^0(B_+,\P^\infty_+\wedge S^1)$$
coincide with the free summand $\omega^0(B_+,\P^\infty_+\wedge S^1)$?}
\\ 
\\

We come back to   the description of  $\alpha^*(x)$: Using Remarks  \ref{viaunits}  and   \ref{firstdesc} one gets the following descriptions of $\alpha^*(x)$.
\begin{pr} For every presentation $(E,F)\in x$ there exist    canonical isomorphisms
$$\alpha^*(x)\cong \qmod{\hat\alpha^*(E,F)}{\tilde J(I[K^{-1}(B)])\hat\alpha^*(E,F)} \ .
$$
Since $\tilde J(I[K^{-1}(B)])$ is contained in $\omega^0(B_+,\P^\infty_+\wedge S^1)$, which is an ideal of $ {_{S^1}\alpha}^0(B)$, we  get  epimorphisms
$$\alpha^*(x)\map \qmod{\hat\alpha^*(E,F)}{\omega^0(B_+,\P^\infty_+\wedge S^1)\cdot \hat\alpha^*(E,F)}\ .
$$
\end{pr}

\subsection{Stabilization}

In this section we will show that the morphism 
\begin{equation}\label{stab}
\tau_*:{_{S^1}\alpha}^k_B(S(E)_{+B}, F^+_B)\to{_{S^1}\alpha}^k_B(S(E')_{+B}, [F']^+_B) 
\end{equation}
 associated with a morphism $\tau:(E,F)\to (E',F')$ in the category ${\cal T}(x)$  is an isomorphism as soon as the rank $f$ of $F$ is sufficiently large. In other words, for fixed $k$, the groups $\alpha^k(x)$ can be computed using only presentations $(E,F)$ with a priori bounded ranks.

\begin{pr} Suppose that $B$ is a finite CW complex. The stabilization morphism (\ref{stab}) is an isomorphism for $2f\geq \dim(B)-k$.
\end{pr}
\pf
 A morphism $\tau$ defines a bundle $U$ and  isomorphisms $E'\cong E\oplus U$, $F'\cong F\oplus U$.  The  long exact sequence  associated with the cofiber sequence over $B$
 $$S(U)_{+B}\map SE'_{+B}\textmap{c} S(E)_{+B}\wedge_B U^+_B\ ,
 $$
and the target space $[F']^+_B$  contains the segment
 $$\to{_{S^1}\alpha}_B^{k-1}(S(U)_{+B}, [F']^+_B)\textmap{\partial}{_{S^1}\alpha}_B^{k}( S(E)_{+B}\wedge_B U^+_B, [F']^+_B)\textmap{c^*}   {_{S^1}\alpha}_B^{k}(SE'_{+B}, [F']^+_B).
 $$
The morphism $\tau_*$ is defined by $c^*$ via the identification ${_{S^1}\alpha}_B^{k}( S(E)_{+B}, F^+_B)={_{S^1}\alpha}_B^{k}( S(E)_{+B}\wedge_B U^+_B, [F']^+_B)$, so it is an isomorphism as soon as
$${_{S^1}\alpha}_B^{k-1}(S(U)_{+B}, [F']^+_B)={_{S^1}\alpha}_B^{k}(S(U)_{+B}, [F']^+_B)=0\ .
$$
Suppose for simplicity $k\geq 0$.  A class $u\in{_{S^1}\alpha}_B^{k}(S(U)_{+B}, [F']^+_B)$ is represented by an $S^1$-equivariant  pointed map over $B$
 $$\varphi:S(U)_{+B}\wedge_B\xi^+_B=\qmod{S(U)\times_B\xi^+}{_B\ S(U)\times_B\infty_\xi} \map [F'\oplus\underline{\R}^k\oplus\xi]^+_B \ ,
 $$
where $\xi=\eta\oplus\xi_0$ is the sum of a complex and a real vector bundle.  We may suppose that $\xi_0$ is an oriented bundle, so that all our bundles become oriented bundles. We will prove that any such map is homotopic to  the  map $\varphi_\infty$ which maps the left hand space fiberwise onto the section $\infty_{F'\oplus\underline{\R}^k\oplus \xi}$. Denote by $q:\P(U)\to B$ the bundle projection and put
$$\tilde F':=q^*(F')(1)\ ,\ \tilde \xi:=q^*(\eta)(1)\oplus q^*(\xi_0)\ .
$$
A  map $\varphi$ as above induces   a pointed bundle map $\tilde \varphi:\tilde \xi^+_{\P(U)}\map [\tilde F'\oplus\underline{\R}^k\oplus \tilde \xi]^+_{\P(U)}$ over $\P(U)$, and the assignment $\varphi\mapsto \tilde \varphi$ is a   bijection. But by Corollary \ref{obsbundle} in section \ref{spheres}, any such pointed bundle map is homotopic to the fiberwise constant bundle map as soon as $\dim_\R(\P(U))+\rk(\tilde \xi)<\rk_\R(\tilde F')+k+\rk(\tilde\xi)$. This condition is equivalent to  $2f>\dim(B)-k-2$. Similarly, we will have  ${_{S^1}\alpha}_B^{k-1}(S(U)_{+B}, [F']^+_B)=0$ as soon as  $2f>\dim(B)-k-1$.
\qed

\subsection{The cohomotopy Euler class of an element in $K(B)$}\label{euler}

Let $x\in K(B)$ and consider a presentation $(E,F)\in x$. The map $o_{(E,F)}:S(E)_{+B}\to F^+_B$ which sends the section $+_B$ of $S(E)_{+B}$ to the infinity section of $F^+_B$ and maps any point $e_b\in S(E_b)$ to $0_b$ is an $S^1$-equivariant map of pointed spaces  over $B$, hence it defines an element 
$\{o_{(E,F)}\}\in {_{S^1}\alpha}^0_B(S(E)_{+B}, F^+_B)$.

One has a canonical isomorphism  (see \cite{CJ} Proposition 12.40)
$${_{S^1}\alpha}^0_B(S(E)_{+B}, F^+_B)\cong {_{S^1}\alpha}^0_{S(E)}(S(E)_{+S(E)}, \pi^*(F)^+_{S(E)})\ ,$$
where $\pi:S(E)\to B$ is the obvious projection. Under this isomorphism the class $\{o_{(E,F)}]\}$ maps to   the equivariant Euler class of the bundle $\pi^*(F)$ over $S(E)$. This class is the pull-back of the equivariant Euler class $\gamma(F)\in {_{S^1}\alpha}^0(B_{+B},F^+_B)$ of the bundle $F$ under the projection $S(E)_{+S(E)}\to B_{+B}$.

For any morphism $\tau=(i,j,E_1,F_1,k):(E,F)\to (E',F')$ in the category 
${\cal T}(x)$ one has $\tau_*(\{o_{(E,F)}\})=\{o_{(E',F')}\}$. 
Therefore the assignment $(E,F)\mapsto -\{o_{(E,F)}\}$ defines a {\it tautological  element}  $\gamma(x)\in \alpha^*(x)$. This element  will be called the equivariant cohomotopy  Euler class of $x$.

 \section{Cohomotopy invariants associated with certain non-linear maps between Hilbert bundles}

\subsection{The cylinder construction}
 \label{cylinder}

  Let $(E,F)$ be a pair of Hermitian vector bundles over a compact basis $B$. Let $V$, $W$ be  Euclidean vector spaces, and let $\mu: E\times V\to [F\times W]^+_B$  be an $S^1$-equivariant   map over $B$. We suppose that $\mu$ is fiberwise differentiable and its fiberwise differential is continuous on $E\times V$. The equivariance property implies that
\begin{equation}
\mu(0^E\times V)\subset \left[0^F\times W
\right]^+_B .
\end{equation}

We assume that $\mu$ has the following  properties:\\
\begin{description}
\item [P1] (properness) There exist positive constants $c$, $C$ such that
$\|\mu(e,v)\| > c$ for all pairs $(e,v)\in E\times V$ with
$\|(e,v)\|\geq C$.
\\  
\item[P2] (restriction to the $S^1$-fixed point set) 
\begin{enumerate} 
\item  There exists a direct sum decomposition $W=H\oplus W_0$ such that 
$$\mu(0^E_y,v)=h(y)+l(v)\ ,\ \forall y\in B\ ,\ \forall v\in V\ ,
$$
where $l:V\textmap{\simeq} W_0\subset W$ is a linear isomorphism, which does
not depend on $y$, and $h:B\to H$ is a continuous map. 
\item   There exists $\varepsilon_0>0$ such that
\begin{equation}\label{restriction}
\|h(y)\|=\|\p_{H}(\mu(0^E_y,v))\| \geq \varepsilon_0\ \  \forall (y,v)\in B\times V\ .
\end{equation}%
\end{enumerate}
\end{description}

We fix an orientation $\oo$ of $H$, and set $b:=\dim(H)$. Choose numbers $R\geq C$ and  $\varepsilon\leq \min(c,\varepsilon_0)$.  The restriction $\mu_{R}$ of $\mu$ to  $D_R(E)\times D_R(V)$ satisfies
$$\|\mu(e,v)\|\geq \varepsilon\ \ \forall (e,v)\in \partial \left[D_R(E) \times D_R(V))\right]\cup \left[0^E\times D_R(V)  \right]\ .
$$
  Therefore,  $\mu_R$ defines an $S^1$-equivariant morphism of pairs over $B$
$$\mu_{R,\varepsilon}: \left(D_R(E)\times D_R(V), \partial[D_{R}(E)\times D_R(V)]\cup [0^E\times D_R(V)] \right) \map\ \ \ \ \ \ \ \ \  \ \ \ \ \ \ \ \ \  \ \ \ \ \ \ \ \ \  $$
$$\ \ \ \ \ \ \ \ \  \ \ \ \ \ \ \ \ \  \ \ \ \ \ \ \ \ \  \ \ \ \ \ \ \ \ \  \ \ \ \ \ \ \ \ \  \ \ \ \ \ \ \ \ \  \map \left(  [F\times W]^+_B, [F\times W]^+_B\setminus \cringle{D}_{\varepsilon}(F\times W)\right)\ .
$$
The first space $D_R(E)\times D_R(V)$ of the pair on which $\mu_{R,\varepsilon}$ is defined can be regarded as  a ``cylinder bundle" over $B$, whose base is the complex disk bundle $D(E)$; the second space of this pair is the union of the boundary of this cylinder bundle with the core  $0^E\times D_R(V)$.
Using polar coordinates in $D_R(E)$ we obtain a map $S(E)\times [0,R]\to D_R(E)$, hence a map
$$\rho:S(E)\times [0,R]\times D_R(V)=S(E)\times D_R(\R\oplus V) \to D_R(E)\times D_R(V)\ ,$$
which maps 
$$[ S(E)\times \{0,R\}\times D_R(V)]\cup [S(E)\times [0,R]\times S_R(V)]=S(E)\times S_R(\R\oplus V) $$
onto the the second component of the pair on which $\mu_{R,\varepsilon}$ is defined. Here we used suitable models
$D(\R\oplus V)$, $S(\R\oplus V)$ for the disc and the sphere in $\R\oplus V$.
Therefore, composing $\mu_{R,\varepsilon}$ with $\rho$ we get   an $S^1$-equivariant map   of pairs over $B$
$$\big(S(E)\times [0,R]\times D_R(V), S(E)\times(\{0,R\}\times D_R(V)\cup[0,R]\times S_R(V))\big)=$$
$$ 
\left(S(E)\times D_R(\R\oplus V),
S(E)\times S_R(\R\oplus V) \right)\to \left([F
\times W]^+_B,\left[F\times W\right]^+_B\setminus
\cringle{D}_{\varepsilon}(F
\times W)\right) 
$$
which we denote by the same symbol $\mu_{R,\varepsilon}$. Collapsing fiberwise  over $B$   the second terms of the two  pairs,  and composing with the natural isomorphism
$$\qmod{[F
\times W]^+_B}{_B\left[F\times W\right]^+_B\setminus
\cringle{D}_{\varepsilon}(F
\times W)}\simeq [F\times W]^+_B\ ,
$$
one gets   an $S^1$-equivariant map   of  pointed spaces over
$B$
$$\mu_{R,\varepsilon}: \qmod{S(E)\times [\R\oplus V]^+ }{_B\  S(E)\times \{\infty\}}= S(E)_{+B}\wedge_B[B\times(\R\oplus V)]^+_B  \map  \left[F\times W
 \right]^+_B .
$$
Using the isomorphism $l:V\stackrel{\simeq}{\to} W_0$ and an orientation
preserving isomorphism
$\R^b\simeq H$, we obtain an element
$$\{\mu\}\in {_{S^1}\alpha}_B^{b-1}(S(E)_{+B},F^+_B)\ ,
$$
which is obviously independent of the choice of the pair $(R,\varepsilon)$. This element will be called the {\it cohomotopy invariant} of $\mu$.

\subsection{General properties of the invariant $\{\mu\}$}

\subsubsection{A vanishing property.}

Let $\mu:E\times V\to [F\times W]^+_B$ be a map satisfying {\bf P1}, {\bf P2}.
\begin{pr}\label{vanishing} If  $\resto{\mu}{D_C(E)\times D_C(V)}$ is nowhere vanishing, then $\{\mu\}=0$.
\end{pr}
\pf We take $\varepsilon\leq\inf\{\|\mu(e,v)\|\ |\ \|e\|\leq C,\ \|v\|\leq C\}$, and we note that the $\left\{[F \times W]^+_B\right\}/_B\big\{\left[F\times W\right]^+_B\setminus
\cringle{D}_{\varepsilon}(F \times W)\big\}$-valued pointed map induced by $\mu_{R,\varepsilon}$ is fiberwise constant.
\qed
\subsubsection{Homotopy invariance.}

Let $\mu'$, $\mu'': E\times V\to [F\times W]^+_{B}$ two maps satisfying properties {\bf P1},   {\bf P2} with constants $C'$, $c'$, $\varepsilon_0'$, and $C''$, $c''$, $\varepsilon_0''$. We suppose that the property {\bf P2}  of the two maps holds for the same decomposition $W=H\oplus W_0$  of  $W$ and for the same isomorphism $l:V\to W_0$. 
 We introduce the notations 
$$\tilde B:= B\times[0,1]\ ,\ \tilde E:=E\times [0,1]=\p_B^*(E)\ , \ \tilde F:=F\times [0,1]=\p_B^*(E)\ .$$
\begin{pr}\label{homotopy} Suppose there exists $C\geq\max(C',C'')$ and a continuous  $S^1$-equivariant map $\tilde \mu: D_C(\tilde E)\times D_C(V)\to [\tilde F\times W]^+_{\tilde B}$  over $\tilde B$  whose restriction to 
$$\partial\left[D_C(\tilde E)\times D_C(V)\right]\cup \left[0^{\tilde E}\times D_C(V)\right]
$$
is  nowhere vanishing, and such that $\resto{\tilde \mu}{t=0}=\mu'$, $\resto{\tilde \mu}{t=1}=\mu''$. Then $\{\mu'\}=\{\mu''\}$ in ${_{S^1}\alpha}_B^{b-1}(S(E)_{+B},F^+_B)$.
\end{pr}
\pf The stable classes  $\{\mu'\}$, $\{\mu''\}$ can be computed using the the cylinder $D_C(\tilde E)\times D_C(V)$ and taking
$$\varepsilon\leq \min\left(\varepsilon_0',\varepsilon_0'',c',c'',\inf\left\{\|\tilde\mu(y)\|\ y\in \partial [D_C(\tilde E)\times D_C(V) ]\cup  [0^{\tilde E}\times D_C(V)]\right\}\right)$$
Applying the cylinder construction  with parameters $C$, $\varepsilon$ to the map $\tilde \mu$ we obtain a homotopy between the corresponding representatives of the classes $\{\mu'\}$, $\{\mu''\}$.
\qed
\subsubsection{A product formula.}\label{prodsubsection}

Let  $V_i$, $W_i$ be Euclidean spaces , $E_i$, $F_i$  Hermitian bundles over a compact base $B$ ($i=1,\ 2$) and $\mu_i:E_i\times  V_i\to [F_i\times W_i]^+_{B}$   $S^1$-equivariant maps over $B$ satisfying the properties {\bf P1}, {\bf P2} (1) of section \ref{cylinder} with constants $C$, $c$. Let $W_i=H_i\oplus W_{0,i}$ be the corresponding direct sum decompositions, and $l_i:V_i\stackrel{\simeq}{\to} W_{0,i}$,  $h_i:B\to H_i$ the maps given by {\bf P2} (1). Fix orientations on the $H_i$, and put 
$$V:=V_1\oplus V_2,\ W:=W_1\oplus W_2,\ H:=H_1\oplus H_2,\ W_0:=W_{0,1}\oplus W_{0,2}\ ,\ l:=l_1\oplus l_2\ ,$$
 and  consider the bundles $E:=E_1\oplus  E_2$, $F:=F_1\oplus F_2$. We have a product map 
$$\mu:E\times V=[E_1\times V_1]\oplus [E_2\times V_2]  \map  [F\times W]^+_B=[F_1\times W_1]^+_B\wedge_B [F_2\times W_2]^+_B$$
 over  $B$.  This map  satisfies properties {\bf P1}, {\bf P2} (1) with the map 
 $$h=(h_1,h_2):B\to H\ .$$
 Note that $\mu$ will also satisfy {\bf P2} (2) as soon as one of the two maps $\mu_1$, $\mu_2$ has this property. Suppose that  $\mu_1$ also satisfies property  {\bf P2} (2) with constant $\varepsilon_0$ and denote by 
 $$\{\mu_1\}\in {_{S^1}\alpha}_B^{b_1-1}(S(E_1)_{+B},[F_1]^+_{B})$$
  the corresponding stable class. The map $\mu_2$ defines a map $[E_2\oplus V_2]^+_{B}\map [F_2\oplus W_2]^+_{B}$ hence a class $\{\mu_2^+\}\in  {_{S^1}\alpha}_{B}^{b_2}([E_2]^+_{B},[F_2]^+_{B})$.   One can then form the product
$$\{\mu_1\}\wedge_B \{\mu_2^+\}\in {_{S^1}\alpha}_B^{b-1}\left(S(E_1)_{+B}\wedge_B[E_2]^+_{B},F^+_{B}\right)\ .
$$
Consider now the contraction map ${_1c}:S(E_1\oplus E_2)_{+B}\to S(E_1)_{+B}\wedge_B
[E_2]^+_B$ introduced in section \ref{alpha} (see formula (\ref{contraction})). Using the  identifications 
$$[E_2]^+_B=\qmod{{D_R(E_2)}}{_B\ {S_R(E_2)}}=\qmod{{E_2}}{_B\ {E_2\setminus \cringle{D}_R(E_2)}}\ ,$$
 we can use  as model for the contraction ${_1c}$ any map  of the form ${_1c}_R^{\Rg}$ given by 
$${_1c}_R^\Rg(e_1,e_2):=\left[\frac{1}{\|e_1\|} e_1, \Rg e_2\right]  \ ,\ (\Rg\geq R)\ . $$
\begin{pr} \label{productformula} Under the above assumptions it  holds $\{\mu\}={_1c}^*\left(\{\mu_1\}\wedge_B \{\mu_2^+\}\right)$.
\end{pr}
\pf The class $\{\mu\}$ is represented by the map of pairs
$$\mu_R:\left(S(E)\times [0,R]\times D_{R}(V),S(E)\times\left([0,R]\times S_R(V)\cup \{0,R\}\times D_{R}(V)\right)\right)\map
$$
$$\map ([F\times W]^+_B\ ,\ [F\times W]^+_B\setminus D_\varepsilon(F\times W))
$$
which is defined by 
$$\mu_R(e_1,e_2,\rho,v_1,v_2)=[\mu_1(\rho e_1,v_1),\mu_2(\rho e_2,v_2)]\ .$$
 The class   ${_1c}^*\left(\{\mu_1\}\wedge_B \{\mu_2^+\}\right)$ is represented by the map $\nu_R^\Rg$  between the same pairs defined by
$$\nu_R^\Rg(e_1,e_2,\rho,v_1,v_2)=\left[\mu_1(\rho \frac{1}{\|e_1\|} e_1,v_1),\mu_2(\Rg e_2,v_2)\right]\ .
$$
Composing $\mu_R$,  $\nu_R^\Rg$ with the projection 
$$p:[F\times W]^+_B\map \qmod{[F\times W]^+_B}{_B\ [F\times W]^+_B\setminus D_\varepsilon(F\times W)}
$$
we obtain two maps
$$m_0\ ,\ m_1:S(E)\times [0,R]\times D_{R}(V)\map \qmod{[F\times W]^+_B}{_B\ [F\times W]^+_B\setminus D_\varepsilon(F\times W)}\simeq [F\times W]^+_B
$$
which map $S(E)\times\left([0,R]\times S_R(V)\cup \{0,R\}\times D_{R}(V)\right)$ onto the infinity section   in the right hand bundle. The natural homotopy between these maps is the map
$$m:[0,1]\times S(E)\times [0,R]\times D_{R}(V)\map \qmod{[F\times W]^+_B}{_B\ [F\times W]^+_B\setminus D_\varepsilon(F\times W)}
$$
given by 
$$m(t,e_1,e_2,\rho,v_1,v_2)=\left[\mu_1\left(\rho\left[1-t+t\frac{1}{\|e_1\|}\right]e_1,v_1\right),\mu_2\left([(1-t)\rho +t\Rg]e_2,v_2\right)\right]
$$
{\bf Claim:} For any $R\geq\sqrt{2} C$ and sufficiently large $\Rg\geq R$ it holds
\begin{enumerate}
\item  the map $m$ is well defined and continuous  at the points  $(t,e_1,e_2,\rho,v_1,v_2)$ with $e_1=0$.
\item the map $m$ maps $[0,1]\times S(E)\times\left([0,R]\times S_R(V)\cup \{0,R\}\times D_{R}(V)\right)$ to the infinity section   in the right hand bundle.
\end{enumerate}

In fact we show that for $e_2\in [E_2]_{y}$, one has 
$$\lim_{u\to (t,0^{E_1}_{b_1},e_2,\rho,v_1,v_2)} m(u)=\infty_{y}\ ,$$
so $m$ maps the locus $e_1=0$ to the infinity section. Let $\eta_R>0$ be sufficiently small, such that $\|\mu_1(e_1,v_1)\|>\varepsilon_0$ for every $(e_1,v_1)\in D_{\eta_R}(E_1)\times D_R(V_1)$. One has
$$\lim_{e_1\to 0}\left\|\rho \left[1-t+t\frac{1}{\|e_1\|}\right]e_1\right\|= \rho t\ .
$$
When $\rho t <\eta_R$, the first component of $m(t,e_1,e_2,\rho,v_1,v_2)$ will already have a norm larger that $\varepsilon_0$.   When $\rho t \geq \eta_R$, we  obtain (using $\|e_1\|^2+\|e_2\|^2=1$):
$$\lim_{e_1\to 0}\left\|[(1-t)\rho +t\Rg]e_2\right\|=(1-t)\rho +t\Rg\geq \eta_R (\frac{1}{t}-1)+t\Rg\geq 2\sqrt{\eta_R \Rg}-\eta_R\ ,
$$
which will be larger than $R$ when $\Rg$ is sufficiently large. The second part of the claim is obvious for the spaces $[0,1]\times S(E)\times[0,R]\times S_R(V)$, $[0,1]\times S(E)\times  \{0\}\times D_{R}(V)$. For $\rho=R$ we obtain
$$\left\|\rho\left[1-t+t\frac{1}{\|e_1\|}\right]e_1\right\|^2+ \|[(1-t)\rho +t\Rg]e_2\|^2\geq R^2(\|e_1\|^2+\|e_2\|^2)= R^2\geq 2C^2\ ,
$$
so at least one of the two norms is $\geq C$.\\

Using the claim, it follows that $m$ descend to an homotopy between two representatives  of  the classes $\{\mu\}$ and ${_1c}^*\left(\{\mu_1\}\wedge_B \{\mu_2^+\}\right)$.
\qed
\vspace{2mm}\\
 An interesting case is the one when also  $\mu_2$ satisfies property {\bf P2} (2). In this case the cylinder construction applies to $\mu_2$ and one can write
 $$\{\mu_2^+\}=\partial_2(\{\mu_2\})\ ,
 $$
 where $\{\mu_2\}\in {_{S^1}\alpha}_B^{b_2-1}(S(E_2)_{+B},[F_2]^+_{B})$ is the invariant associated with $\mu_2$, and $\partial_2$ is the connecting morphism in the long exact sequence associated with the cofiber sequence
 $$S(E_2)_{+B}\map D(E_2)_{+B}\map [E_2]^+_{B}\ .
 $$
Let ${_{2}c}:S(E_1\oplus E_2)_{+B}\to [E_1]^+_{B}\wedge_B S(E_2)_{+B}$ be the standard contraction. In this case, our multiplication formula  becomes
 \begin{co} \label{torsion} Suppose that both maps $\mu_1$, $\mu_2$ satisfy properties {\bf P1}, {\bf P2}. Then 
 $$\{\mu\}={_1c}^*\left(\{\mu_1\}\wedge_B \partial_2(\{\mu_2\})\right)={_2c}^*\left(\partial_1(\{\mu_1\})\wedge_B \{\mu_2\}\right)\ .$$
 \end{co}

Another  corollary is obtained  when $\mu_2$ is defined by a pair of linear isomorphisms $E_2\to F_2$, $V_2\to W_2$. The corresponding formula will play an important role in the proof of the coherence  Lemma \ref{compat} comparing the invariants associated to two finite dimensional approximations of an admissible bundle map between Hilbert bundles.

\begin{pr} \label{smashidnew}  Let  $\mu:E\times V\to F\times W$ be a map satisfying the properties {\bf P1}, {\bf P2} with constants $C$, $c$, $\varepsilon_0$ and maps  $l:V\to W_0$, $h:B\to H$.  Let $a:E'\to F'$ be an isomorphism of complex vector bundles over $B$, and let $b:V'\to W'$  be an isomorphism of real vector spaces. Put $\tilde E:=E\oplus E'$, $\tilde F:=F\oplus F'$, $\tilde V:=V\oplus V'$, $\tilde W:=W\oplus W'$, and define 
$$\tilde \mu(e,e',v,v')=\iota[\mu(e,v)\wedge_B (a(e'),b(v'))]\ ,
$$
where  $\iota$ is the obvious identification 
$$\iota:[F\times W]^+_B\wedge_B (F'\times W')^+_B \to [(F\oplus F')\times (W\oplus W')]^+_B\ . $$
Then
\begin{enumerate}
\item $\tilde \mu$ satisfies  {\bf P1}  with constants $C$, $\gamma$ (for sufficiently small $0<\gamma< c$), and   {\bf P2}  with constant $\varepsilon_0$ and maps $\tilde l:=l\oplus b$, $\tilde h:=h$. 
\item   $\{\tilde \mu\}=\tau_*(\{\mu\})$, where $\tau$ denotes  the obvious  morphism $(E,F)\to (\tilde E, \tilde F)$. 
\end{enumerate}
\end{pr}
 The second statement follows directly from Proposition \ref{productformula}. The first statement (which is specific to  the  case when the second factor is a linear isomorphism) is proved as follows: Since the closed set $\mu^{-1}(D_c(F\times W))$ is contained in the open disk $\cringle{D}_C(E\times V)$, there exists   $r>0$  such that $\|\mu(e,v)\|> c$ as soon as $\|(e,v)\|\geq C-r$.  For a point $(e,e',v,v')$  with $\|(e,e',v,v')\|\geq C$ one has either $\|(e,v)\|\geq C-r$, or $\|(e',v')\|\geq r$. In the first case one obtains $\|\mu(e,v)\|> c$, whereas in the second we get $\|(a(e'),b(v')\|\geq c' r$ for a constant $c'$.
 \qed 

\subsection{A class of non-linear maps between Hilbert bundles} \label{general}

Suppose now that ${\cal V}$, ${\cal W}$ are real Hilbert spaces, and that ${\cal E}$, ${\cal F}$ are complex Hilbert bundles over the compact  basis $B$, and let $\mu: {\cal E}\times {\cal V}\to  {\cal F}\times {\cal W}$ be a continuous  $S^1$-equivariant  map over $B$ which is fiberwise ${\cal C}^\infty$,   and whose fiberwise derivatives are continuous on ${\cal E}\times {\cal V}$. We   assume that the fiberwise differentials
$$d_y:=d_{0_y}\mu_y= {\cal E}_y\times{\cal V}\map {\cal F}_y\times{\cal W} \ ,\ y\in B
$$
at the origins of the fibers $ {\cal E}_y\times {\cal V}$ are Fredholm. The linear operator $d_y$ has the form $d_y=(\delta_y,l_y)$, where $\delta_y:{\cal E}_y\to {\cal F}_y$ and $l_y:{\cal V}\to{\cal W}$ are defined by the derivatives of the restrictions $\resto{\mu}{{\cal E}_y\times \{0^{\cal V}\}}$, $\resto{\mu}{\{0^{\cal E}_y\}\times {\cal V}}$. Note that the continuous  family $\delta:=(\delta_y)_{y\in B}$ of complex Fredholm operators defines an element ${\rm ind}(\delta)\in K(B)$. Let $d:{\cal E}\times {\cal V}\to{\cal F}\times  {\cal W}$ the fiberwise  linear map defined by the family of Fredholm operators $(d_y)_{y\in B}$. We suppose that   $\mu$ also has the properties   
\begin{description}
\item [${\cal P}1$] (properness) There exist positive constants $c$, $C$ such that
$\|\mu(e,v)\| > c$ for all pairs $(e,v)\in {\cal E}\times{\cal V} $ with
$\|(e,v)\|\geq C$.
\vspace{1mm}  
\item[${\cal P}2$] (behavior near   the $S^1$-fixed point set)
\begin{enumerate} 

\item ${\cal W}$ splits orthogonally as ${\cal W}=H\oplus {\cal W}_0$, where 
$H$ is a finite dimensional subspace, and for every 
$y\in B$ one has
$$\mu(0^E_y,v)=h(y)+l(v)\ \ \forall y\in B\ ,\ \forall v\in {\cal V}\ ,
$$
where $l:{\cal V}\textmap{\simeq}{\cal W}_0\subset{\cal W}$ is a linear
 isometry, and  $h$ is a map from $B$ to $H$.

In particular the operator $l_y$ coincides with $l$, so is independent of $y$.

\item There exists $\varepsilon_0>0$ such that for every $y\in B$ one has
$$\|h(y)\|=\|\p_H(\mu(0_y^{\cal E},v))\|\geq \varepsilon_0\ .
$$
\end{enumerate}
\item [${\cal P}3$] (linear+compactness)  The difference $k:=\mu-d$  is  globally compact,  in the sense that for every $R>0$ the image $k(D_R({\cal E}\times{\cal V}))$ of the disk bundle $D_R({\cal E}\times{\cal V})$ is relatively compact in the total space ${\cal F}\times{\cal W}$.
\end{description}
Note that one has the identity
\begin{equation}\label{kk}
k(0^{\cal E}_y,v)=h(y)\in H\ ,\ \forall y\in B\ .
\end{equation}
In the next section we will see that the left hand of the Seiberg-Witten equations on a 4-manifold $M$ defines a map satisfying  properties ${\cal P}_1$ -- ${\cal P}_3$. A different construction of such a map  can be found in \cite{BF}.  
\subsection{The Seiberg-Witten map in dimension 4} 
\label{swmap}

Let $M$ be closed oriented   4-manifold, and let  $L$ be a Hermitian  line bundle on $M$. We fix the following data:
\begin{enumerate}
\item A closed complement ${\cal S}$ of the closed subspace
$iB^1_{\rm DR}(M)=d(iA^0(M))$ of
$iA^1(M)$.
\item A closed complement ${\cal V}$ of the finite
dimensional space
$$i\H^1:=S\cap\ker(d:iA^1(M)\to iA^2(M))\simeq iH^1(M,\R)$$ in ${\cal S}$
\item A complement $i\H^2$ of $d(iA^1(M))$ in $\ker(d:iA^2(M)\to
iA^3(M))$. This complement will come with an isomorphism $i\H^2\simeq
iH^2(M,\R)$.
\item  An affine subspace ${\cal A}$ of the space of connections  ${\cal
A}(L)$ modeled after ${\cal S}$.
\end{enumerate}

Therefore, ${\cal A}$ is a slice to the orbits of the right action of the gauge
group ${\cal G}$ on the space of connections:
$$a\cdot g:=a+ 2 g^{-1}dg
$$

The quotient $\bar {\cal A}:={\cal A}/{\cal V}$ is an affine space modeled after
$iH^1(M,\R)$. Consider the finite dimensional Lie group
$$G:=\{u\in{\cal C}^\infty(M,S^1)|\  u^{-1}du\in {\cal S}\}\ .
$$

One has an obvious short exact sequence
$$\{1\}\map S^1\map G\textmap{\lambda} 2\pi i H^1 (M;\Z)\map\{1\}\ ,
$$
where   $\lambda$ is defined by
$u\mapsto [u^{-1}du]_{\rm DR}$.  The choice of a point $x_0\in M$ defines a
left splitting
${\rm ev}_{x_0}:G\to S^1$ whose kernel is isomorphic to $2\pi i
H^1 (M;\Z)$ and which will be denoted by
$G_{x_0}$. In the  affine space ${\cal A}$ we have a natural $i\H^1$-invariant (hence
$G_{x_0}$-invariant) subset ${\cal A}_0$ defined by
$${\cal A}_0:=\{a\in {\cal A}\vert\ F_a\in i\H^2\}\ .
$$
The curvature $F_{a_0}$ of a connection $a_0\in {\cal A}_0$ is independent of $a_0$, because it coincides with the representative in $i\H^2$ of the de Rham class $-2\pi ic_1^{\rm DR}(L)$; this 2-form will be denoted by $F_0$.
Note that ${\cal A}_0$ is a $G_{x_0}$-invariant complete system of representatives for
the quotient $\bar {\cal A}={\cal A}/{\cal V}$. The space ${\cal A}/{G_{x_0}}$ can be regarded as   an affine bundle over
the torus
$$\Pic(L):=\qmod{\bar {\cal A}}{G_{x_0}} \ ,
$$
  which is naturally a $iH^1(X;\R)/4\pi i
H^1(X;\Z)$-torsor. The fibers of the affine bundle 
$$\pi:{\cal A}/{G_{x_0}}\map \Pic(L)$$
 are affine
${\cal V}$-spaces. Since the quotient
${\cal A}_0/G_{x_0}$ is a section of this affine bundle,  we can regard it as a
${\cal V}$-vector bundle over $\Pic(L)$ with ${\cal A}_0/G_{x_0}$ as zero section. This vector
bundle is  actually trivial: indeed, the map $(a_0,v)\mapsto a_0+v\in {\cal A}$
is $G_{x_0}$-equivariant, and it descends to a trivialization 
$\Pic(L)\times {\cal V}\to {\cal A}/{G_{x_0}}$.
\begin{re} Choosing a Riemannian metric $g$ on $M$ gives canonical choices for the
three objects $S$, $T$, $i\H^2$ above, namely
$${\cal S}= \ker(d^*: iA^1(M)\map iA^0(M))\ ,\ {\cal V}:= d^*(iA^2(M))\ ,\ i\H^2=
i\H^2_g\ ,
$$
where the subscript $_g$ on the right denotes the respective $g$-harmonic
space. With these choices, ${\cal A}_0$ is just the the set of $g$-Yang-Mills
connections in the slice ${\cal A}$.
\end{re}

Let  $g$ be a Riemannian metric on $M$, $\cg\in Spin^c(M)$ an equivalence class of $Spin^c$-structures, and let $\tau:Q\to P_g$ be a $Spin^c$-structure on $M$ representing  the class $\cg$ . Denote by $\Sigma^\pm$,
$\Sigma:=\Sigma^+\oplus\Sigma^-$ the spinor bundles of
$\tau$, $L=\det(\Sigma^\pm)$ the determinant line bundle,  and
$\gamma:\Lambda^1\to\End_0(\Sigma)$ the Clifford map [OT]. Note that the gauge group $\Aut(Q)$ of $Q$ acts on the space of $Spin^c$-structures $\tau:Q\to P_g$ representing $\cg$ (or, equivalently, on the space of Clifford maps $\gamma:\Lambda^1\to\End_0(\Sigma)$ which are compatible with $\cg$). Therefore, the space of $Spin^c$-Dirac operators which are compatible with the pair $(g,\cg)$ has a very complicated topology. Note that, for the construction of a Dirac operator one needs a concrete $Spin^c$-structure $\tau$ (or, equivalently, a concrete Clifford map $\gamma$), not only an equivalence class $\cg$.  \\

The gauge group ${\cal G}$ and its subgroup $G_{x_0}$ act from the left on
the vector spaces of sections $A^0(\Sigma^\pm)$
by the formula
$$(g,\Psi)\mapsto g^{-1}\Psi\ .
$$

Since $G_{x_0}$ acts freely on the affine quotient space $\bar {\cal A}$ we 
get two flat
vector bundles $\bar {\cal A}\times_{G_{x_0}}A^0(\Sigma^\pm)$ over $\Pic(L)$ with standard fibers $A^0(\Sigma^\pm)$.   
In order to use our general formalism  we make the following definitions:
$$B:=\Pic(L)\ ,\  {\cal E}:=\bar {\cal A}\times_{G_{x_0}}A^0(\Sigma^+)\ ,\ {\cal F}:=\bar {\cal A}\times_{G_{x_0}}A^0(\Sigma^-)\ ,\   {\cal W}:=  i A^2_+(M)\ .
$$

Let $\kappa:B\to i\H^+_{g}$ be a smooth map. The   $\kappa$-twisted Seiberg-Witten
map is the map from
$A^0(\Sigma^+)\times {\cal A}$ to $A^0(\Sigma^-)\times iA^2_+$
given by
$$(\Psi,a)\mapsto ( \Dr_a\Psi, (F_a - F_0+\kappa(\pi(a)))^+-\gamma^{-1}((\Psi\bar\Psi)_0)\ .
$$
Via the identification $B\times {\cal V}={\cal A}/{G_{x_0}}$ this map descends to an $S^1$-equivariant  map
$$sw_\kappa:{\cal E} \times {\cal V}\map {\cal F} \times {\cal W}\ .
$$
The restriction of $sw_\kappa$ to the fiber over
$y=[a_0]\in B$ is given by the formula
$$sw_\kappa(\Psi,v)=\left(\Dr_{a_0}\Psi+\frac{1}{2}\gamma(v)\Psi\ ,\ d^+v+\kappa(y)-\gamma^{-1}((\Psi\bar\Psi)_0\right)\ .
$$

The linearization of this map at the zero section in the bundle ${\cal E} \times {\cal V}$ over $B$ is a fiberwise  linear bundle map given by
$$d(\Psi,v)=(\Dr_{a_0}\Psi,d^+v) \ .
$$
over the fiber $[a_0]\in B$. Hence $sw_\kappa$  decomposes as
$$sw_\kappa=d+c_\kappa\ ,
$$
where $c_\kappa$ is the sum of a  quadratic map $c$ and the fiberwise constant map defined by $\kappa$. Denote by $w_\tau$ the expected dimension of the Seiberg-Witten moduli space
corresponding to $\tau$:
$$w_\tau:=\frac{1}{4}(c_1(L)^2-  3\sigma(M)-2e(M ))
$$
We define  Sobolev $L^2_k$-completions of the spaces ${\cal V}$, ${\cal W}$ in the usual way. The construction of Sobolev norms on the bundles ${\cal E}$, ${\cal F}$  is more delicate, because these bundles are quotients  with respect to the group $G_{x_0}$, which does not operate by $L^2_k$-isometries\footnote{We are grateful to Markus Bader for pointing out this subtility to us.}. For a point $y=[a_0]\in B$ (with $a_0\in{\cal A}_0$) one identifies the fibers ${\cal E}_y$, ${\cal F}_y$ with $\{a_0\}\times A^0(\Sigma^\pm)$ and uses the covariant derivatives associated with $\nabla_{a_0}$ to define the $L^2_k$-norm on ${\cal E}_y$. A gauge transformation $g\in G_{x_0}$ defines an isometry $
\{a_0\}\times A^0(\Sigma^+)\to \{a_0\cdot g\}\times A^0(\Sigma^\pm)$, so in this way one obtains a well defined Sobolev norm on the fiber ${\cal E}_y$.

\begin{lm} With respect to suitable Sobolev completions, the following holds:
\begin{enumerate}
\item $sw_\kappa$ is smooth.
\item The fiberwise  linear map $d$ is fiberwise Fredholm of index $w_\tau-b_1+1$, and $c_\beta$ is a compact map.
\item There exists
positive constants $c$, $C$ such that
$$\|(\Psi,v)\|\geq C\Rightarrow \| sw_\kappa(\Psi,v)\|> c\ .
$$
\item The map $c_\kappa=sw_\kappa-d$ is compact.
 \end{enumerate}
\end{lm}
\vspace{2mm}
Therefore   the Seiberg-Witten map $sw_\kappa$ satisfies always the properties  ${\cal P1}$, ${\cal P2}$ (1) and ${\cal P3}$ in section \ref{general}. It also satisfies  ${\cal P2}$ (2) for all maps $\kappa:B\to i\H^+_g\setminus\{0\}$.

The first  and the third statements in the lemma are easy to see.  The crucial properness assertion (2) is stated in  \cite{Fu1}, \cite{Fu2}. A proof of the analogue statement for another version of the Seiberg-Witten map can   be found in \cite{BF}. A detailed proof for our version,  and  an analogue properness property  in a different gauge theoretic  context can be found in \cite{B}.  Similar methods can be also used to treat the 3-dimensional Casson-Seiberg-Witten theory.\\
\\
{\bf The universal Seiberg-Witten map:} With the notations introduced above we fix the parameters $(g,\kappa,\cg, Q)$ on $M$. As we explained above, the family of Dirac operators $\delta:=(\Dr_{a_0})_{[a_0]\in B}$ (and  implicitly the Seiberg-Witten map $sw_\kappa$) still depends on the choice of a $Spin^c$-structure $\tau:Q\to P_g$ in the class $\cg$. This parameter varies in the space  $\Gamma:=\Hom_M(Q,P_g)$ of equivariant bundle morphisms $Q\to P_g$ covering $\id_M$. Since this space has a complicated topology, and our purpose is to construct an invariant which is intrinsically and canonically associated with the base manifold, it is important to understand how the objects $({\cal E}, {\cal F}, \delta, sw_\kappa)$ associated with different bundle morphisms $\tau$ should  be identified. A   construction which has been presented by Furuta in his talk at the Postnikov Memorial Conference \cite{Fu3}  solves this problem in an elegant way\footnote{Furuta explained the details of this construction in an e-mail to the second author, and informed us that similar ideas have been used before.}:

One has a universal family $(\Dr_{a}^\tau)_{(\tau,a)\in \Gamma\times{\cal A}(L)}$ of Dirac operators, which is intrinsically  associated with the system $(g,\kappa,\cg, Q)$. 
An automorphism $f\in\Aut(Q)$ defines automorphisms $f_{\pm}$ of $\Sigma^\pm$ and an automorphism $\det(f)\in\Aut(L)$; the relation between the Dirac operators associated with $\tau$ and $\tau':=\tau\circ f$ is
 
$$
\Dr^{\tau'}_{\det(f)^*(a)}=f_-^{-1}\circ\Dr^\tau_{a}\circ f_+\ .
$$

The group $\Aut(Q)$ acts transitively with constant stabilizer ${\cal G}\subset \Aut(Q)$ on $\Gamma$, and acts  with constant stabilizer $S^1$ on the product $\Gamma\times{\cal A}(L)$. 
Now fix a $Spin^c(4)$-equivariant map $\theta:Q_{x_0}\to P_{g,x_0}$, and  put 
$$\Gamma_0:=\{\tau\in\Gamma,\ \tau_{x_0}=\theta\}\ ,\ \Aut(Q)_{\theta}:=\{f\in\Aut(Q)|\ \theta\circ f_{x_0}=\theta\}\ ,$$
$$ \Aut(Q)_{0}:=\{f\in\Aut(Q)|\   f_{x_0}=\id_{Q_{x_0}}\} \ .$$
The quotient $\Aut(Q)_{\theta}/\Aut(Q)_0$  can be identified with $S^1$.

  The universal family $(\Dr_{a}^\tau)_{(\tau,a)\in \Gamma\times{\cal A}(L)}$ of Dirac operators descends to a a family   $\Dr: \E\to \F$   on the  free quotient
$$\B:=\qmod{\Gamma_0\times{\cal A}(L)}{\Aut(Q)_0}\ .
$$
By choosing an element $\tau\in\Gamma_0$ one obtains an identification $\B\simeq{\cal A}(L)/{\cal G}_{x_0}$, where ${\cal G}_{x_0}$ acts by the formula $(a\cdot g)=a+2g^{-1}dg$, but the identification is not canonical. The free action of  ${\cal V}:= d^*(iA^2(M))$ on the second factor ${\cal A}(L)$ by translations induces a free action on $\B$, and the quotient $B$ with respect to this action is a $b_1$-dimensional torus. The space of pairs $(\tau,a)$ with $a$ Yang-Mills defines a section $\B_0$ of the ${\cal V}$-bundle $\B\to B$, which therefore becomes a trivial vector bundle with fiber ${\cal V}$.  One can  construct a  ``universal" Seiberg-Witten map $sw_\kappa$ over $B$ using the space $\Gamma_0\times A^0(\Sigma^+)\times{\cal A}(L)$ as space of configurations, and $\Aut(Q)_{0}$ as gauge group; this map is intrinsically associated with the system $(g,\kappa,\cg, Q,\theta)$. The important point in this construction is that restricting the universal family $\Dr$ to the torus $\B_0\simeq B$ one obtains a ``universal Segal cocycle" $\Dr:{\cal E}\to{\cal F}$ representing the K-theory element $x={\rm ind}({\Dr})$.   Furuta showed that the corresponding spectrum is independent of the choice of $\theta$, up to homotopy.  On can apply the construction in \cite{BF} and get -- for manifolds with $b_+\geq 2$ and arbitrary $b_1$ -- a well defined Bauer-Furuta invariant  belonging  to a homotopy group  which is functorial with respect to diffeomorphisms.

As we explained in the introduction (see section \ref{mot}), we believe that for some applications  it is useful to have invariants belonging to groups  which are {\it topologically functorial}, as it is the case in classical Seiberg-Witten and Donaldson theories.

 \subsection{Finite dimensional approximation }\label{finite}
 
 We will need the following simple geometric construction. Let ${\cal A}$ be a (real or complex) Hilbert space, and $A\subset {\cal A}$ a finite dimensional subspace. Following \cite{BF} we introduce, for every $\varepsilon>0$ the retraction
$$\rho_{\varepsilon,A}:{\cal A}^+\setminus S_\varepsilon(A^\bot)\to A^+
$$
in the following way. For every $a\in A\setminus\{0\}$ put 
$$s_{\varepsilon,a}:=\frac{\|a\|^2-\varepsilon^2}{2\|a\|^2}\  , \ c_{\varepsilon,a}=s_{\varepsilon,a} a\ ,\ r_{\varepsilon,a}:=\frac{\|a\|^2+\varepsilon^2}{2\|a\|}\ .
$$
Let $S_{\varepsilon,a}\subset \R a+ A^\bot$  be the hypersphere of $\R a + A^\bot$  defined by the equation
$$ \|b -c_{\varepsilon,a}\|^2+ \| a'\|^2= r_{\varepsilon,a}^2\ .
$$
The hypersphere $S_{\varepsilon,a}$ has the properties
$$a\in S_{\varepsilon,a}\ ,\ S_\varepsilon(A^\bot)\subset S_{\varepsilon,a}\ .
$$
Consider also the spherical calotte:
$$C_{\varepsilon,a}:=\{ta+a'\in S_{\varepsilon,a}|\ t>0\}\subset S_{\varepsilon,a}\ .
$$
Denote  by $C_{\varepsilon,\infty}\subset [A^\bot]^+$ the exterior   of the sphere $S_\varepsilon(A^\bot)\subset A^\bot$ (including $\infty$), and by $C_{\varepsilon,0}$ its interior.
Now note that
$${\cal F}_{\varepsilon, A}:=\{C_{\varepsilon, a}|\ a\in A^+\}$$
is a foliation of  ${\cal A}^+\setminus S_\varepsilon(A^\bot)$ with closed leaves; the leaves are  all diffeomorphic to the standard disk of $A^\bot$.  The retraction $\rho_{\varepsilon, A}$ assigns the point $a\in A^+$ to any point of the leaf $C_{\varepsilon,a}\subset {\cal A}^+$.  Note that for any $z\in{\cal A}$ one has the implication
\begin{equation}\label{far}
\left(z\in {\cal A}^+\setminus S_\varepsilon(A^\bot),\ \|z\|\geq \varepsilon\right)\Rightarrow\|\rho_{\varepsilon,A} (z)\|\geq \|z\|
\end{equation}
(equality is obtained when $\|z\|=\varepsilon$ or   $z\in A$). A second important property of the retraction $\rho_{\varepsilon,A}$ is
\begin{equation}\label{proj}
z\in {\cal A} \setminus A^\bot \Rightarrow \big(\rho_{\varepsilon,A}(z)=\lambda_{\varepsilon,z}\p_A(z)\hbox{ with } \lambda_{\varepsilon,z}\geq 1\big)\ .
\end{equation}

Any  $\R$-linear isometry $u$ of ${\cal A}$ which leaves the subspace $A$ invariant will also leave invariant the foliation ${\cal F}_{\varepsilon,A}$. Therefore
\begin{re} \label{equirho} $\rho_{\varepsilon, A}$ is equivariant with respect to any $\R$-linear isometry   of ${\cal A}$ which leaves the subspace $A$ invariant.
\end{re}

These retractions play a fundamental role in the following construction of finite dimensional approximations. This construction is a refinement of the one developed  in \cite{BF}. The main difference is that we have to work over a base $B$, and that we treat the real and complex summands separately.

Consider again an $S^1$-equivariant  map $\mu:{\cal E}\times {\cal V}\to {\cal F}\times {\cal W}$ over $B$ satisfying the properties ${\cal P}1$, ${\cal P}2$, ${\cal P}3$ of section \ref{general}. Recall from section \ref{general} that we denoted by $d$ the linearization of $\mu$ at the 0-section and by $\delta$ and $l$ the complex and the real components of $d$.  We have assumed that the $\R$-linear operator $l$ induces an isometry ${\cal V}\to {\cal W}_0$.  A finite rank subbundle $F\subset{\cal F}$ will be called {\it admissible} if it is mapped surjectively onto the linear space defined by the family of cokernels $(\coker(\delta_y))_{y\in B}$.  A finite dimensional  subspace $W\subset {\cal W}$ will be called admissible if it contains $H$. A pair $(F,W)$ will be called admissible if $F$ and $W$ are both admissible; in this case, for every $y\in B$ the product  $F_y\times W$ is mapped surjectively onto $\coker(d_y)$.

For every admissible pair $\pi=(F,W)$  the preimage $d^{-1}(F\times W)$ is a finite rank subbundle of ${\cal E}\times {\cal V}$ which splits as 
$$d^{-1}(F\times W)=\delta^{-1}(F)\times l^{-1}(W)\ .
$$

We denote by $W_0$ the orthogonal complement of $H$ in $W$, and  put  $V:=l^{-1}(W)=l^{-1}(W_0)$,  $E:=\delta^{-1}(F)\subset {\cal E}$. The pair $(E,F)$ represents ${\rm ind}(\delta)\in K(B)$.   We get {\it topological} orthogonal direct sum decompositions
$$ \ {\cal F}=F\oplus F^\bot\ ,\ {\cal E}=E\oplus E^\bot\ ,\ {\cal W}=W\oplus W^\bot=H\oplus W_0\oplus W^\bot\ ,\ {\cal V}=V\oplus V^\bot\ .
$$
The product  $F\times W$ is a finite dimensional Hilbert subbundle of ${\cal F}\times {\cal W}$ whose orthogonal complement is $F^\bot\times W^\bot$. The retraction 
$$\rho_{\varepsilon,F\times W}: [{\cal F}\times {\cal W}]^+_B\setminus S_\varepsilon(F^\bot\times W^\bot)\map [F\times W]^+_B$$
is defined fiberwise.   We will see that, for sufficiently small $\varepsilon>0$ and sufficiently large admissible pairs $\pi=(F,W)$,  the image of the restriction $\resto{\mu}{E\times V}$ does not intersect $S_\varepsilon(F^\bot\times W^\bot)$. Therefore we can define a map
$$\mu_{\varepsilon,\pi}:=\resto{\{\rho_{\varepsilon,F\times W}\circ \mu\}}{E\times V}:E\times V\map [F\times W]^+_B \ ,$$
which belongs to the class studied in section \ref{cylinder}. Such a map will be called a {\it finite dimensional approximation} of $\mu$. The result we need is very much similar to the first part of Lemma 2.3 in \cite{BF}.  We know that the preimage $\mu^{-1}(D_c({\cal F}\times {\cal W}))$ is contained in the disk bundle $D_C({\cal E}\times{\cal V})\subset D_C({\cal E})\times D_C({\cal V})$. The image $k(D_C({\cal E})\times D_C({\cal V}))$ is   relatively  compact in the total space ${\cal F}\times {\cal W}$, because $k$ is   compact  by property ${\cal P}3$.  Now fix $\eta>0$. 
\begin{dt} A pair $\pi:=(F,W)$ is called $\eta$-{\it admissible} if it is admissible, and  any element of  the compact  set $\overline{k(D_C({\cal E})\times D_C({\cal V}))}$ is $\eta$-close to an element in $F\times W$ belonging to the same fiber. 
\end{dt}
\begin{lm}\label{admiss} Let $K\subset  {\cal F}\times{\cal W}$ a compact set, $F$ a finite rank subbundle of  ${\cal F}$, and $W$ a finite dimensional subspace of ${\cal W}$. The following conditions are equivalent:
\begin{enumerate}
\item Any point $k\in K$ is $\eta$-close to a point of $F\times W$ belonging to the same fiber.
\item There exists a finite system $(\phi_1,\dots,\phi_k)$ of sections of $F$ and a finite system $(w_1,\dots,w_k)$ of vectors of $W$ such that 
$$K\subset \union_{y\in B,\ 1\leq i\leq k} B((\phi_i(y),w_i),\eta)\ .$$
\end{enumerate}
\end{lm}
\pf The implication $(2)\Rightarrow (1)$ is obvious.  For the second it is convenient to introduce the notation
$$B((\phi,w),\eta):=\union_{y\in B} B((\phi(y),w),\eta)\ ,
$$
for a section $\phi\in \Gamma({\cal F})$ and a vector $w\in W$. If $K$ satisfies (1) then it is contained in the union of open sets $\cup_{\phi\in\Gamma(F),w\in W} B((\phi,w),\eta)$. It suffices now to use the compactness of $K$.
\qed
\begin{co}\label{admissco} The set of pairs $(F,W)$ satisfying the $\eta$-admissibility condition is non-empty, open and cofinal.  
\end{co}
\pf  Since $K:=\overline{k(D_C({\cal E})\times D_C({\cal V}))}$ is compact in ${\cal F}\times W$ there exists finite  systems $\underline{\phi}=(\phi_1,\dots,\phi_k)\in \Gamma({\cal F})^k$, $\underline{w}=(w_1,\dots w_k)\in {\cal W}^k$ such that 
\begin{equation}\label{incl}
K\subset \union_{i}B((\phi_i,w_i),\eta).
\end{equation} 
  Now fix an admissible pair $(F_0,W_0)$. Since ${\cal F}$ has infinite rank, it is easy to see that any neighborhood of $\underline{\phi}$ contains a system $\underline{\phi}'$ which is in general position  with respect to $F_0$ in the following sense:  for every point $y\in B$ the system  $\underline{\phi}'(y)$ is linearly independent  in ${\cal F}_y$, and $\langle\underline{\phi}'(y)\rangle\cap F_{0,y}=\{0_y\}$.  Since the condition (\ref{incl}) is obviously open with respect to the pair $(\underline{\phi},\underline{w})$, we can choose such a system $\underline{\phi}'$ which still satisfies (\ref{incl}) and is in general position with respect to $F_0$. We  denote by  $F'$ the rank $k$-subbundle generated by $\underline{\phi}'$, and we put $F:=F_0\oplus F'$ and $W:=W_0+\langle w_1,\dots w_k\rangle$. 

To prove that $\eta$-admissibility is open, note that admissibility is open, and use Lemma \ref{admiss} to prove that the second condition in the definition of $\eta$-admissibility is also open. Finally, to see that the set of $\eta$-admissible pairs is cofinal, we fix an $\eta$-admissible pair $(F_0,W_0)$. For
an  arbitary pair $(F,W)$ consider a small deformation $F_0'$ of $F_0$ for which $(F_0',W)$ is still $\eta$-admissible and such that $F_0'$ is fiberwise transversal to $F$. Then $(F\oplus F_0',W+W)$ will be an $\eta$-admissible pair which contains $(F,W)$.
\qed
\begin{lm}\label{finapprox}  (Finite dimensional approximations)  Let $0<\eta<\frac{c}{4}$. Then
\begin{enumerate}
\item For any $\eta$-admissible pair $\pi=(F,W)$ one has  
$$\im\left(\resto{\mu}{E\times V}\right)\cap  S_c(F^\bot\times W^\bot)=\emptyset\ ,$$
so the finite dimensional approximation 
$$\mu_{c,\pi}:=\resto{\{(\rho_{c,F\times W})\circ \mu\}}{E\times V}:E\times V\map (F\times W)^+_B $$
 is defined.
 \item The restriction $\resto{\mu_{c,\pi}}{D_C(E)\times D_C(V)}$ takes values in $F\times W$.
\item For any $\eta$-admissible pair $\pi=(F,W)$ the finite dimensional approximation  $\mu_{c,\pi}$ satisfies the conditions {\bf P1},  {\bf P2} (see section \ref{cylinder}) with the same  constants $C$, $c$, $\varepsilon_0$, isometry $l:{\cal V}\to{\cal W}_0\subset {\cal W}$ and the same map $h:B\to H$ as $\mu$.
\end{enumerate}
\end{lm}
\pf  
1.  If the intersection  $\im\left(\resto{\mu}{E\times V}\right)\cap  S_c(F^\bot\oplus W^\bot)$ was not empty, there would exist a point $(e,v)\in E\times V$ such that $\mu(e,v)\in S_c(F^\bot\times W^\bot)$. Since $S_c(F^\bot\times W^\bot)\subset D_c({\cal F}\times{\cal W})$, it follows $(e,v)\in D_C({\cal E})\times D_C({\cal V})$. Therefore 
$$\mu(e,v)=d(e,v)+k(e,v)\in F\times W_0+k(D_C({\cal E})\times D_C({\cal V}))\ .
$$
But any element in the second set $k(D_C({\cal E})\times D_C({\cal V}))$ is $\eta$-close to an element in $F\times W$ by assumption, so $\mu(e,v)$ is $\eta$-close to $F\times W$. Since $\eta<\frac{c}{4}$, this contradicts $\mu(e,v)\in  S_c(F^\bot\oplus W^\bot)$.
\\ \\
2.  The same argument shows that $\mu(D_C(E)\times D_C(V))$ does not intersect the complement of $D_c(F^\bot\oplus W^\bot)$ in $F^\bot\oplus W^\bot$.
\\ \\
3. We have to check that, for an $\eta$-admissible pair $\pi=(F,W)$, the finite dimensional approximation $\mu_{c,\pi}$ has the two properties {\bf P1}, {\bf P2} in section \ref{cylinder}.  For a point $(e,v)\in E\times V$ with $\|(e,v)\|\geq C$ it holds $\|\mu(e,v)\| > c$  so, by (\ref{far}), we have
\begin{equation}\label{bignorm}
\|\rho_{c,F\times W}(\mu(e,v))\|\geq \|\mu(e,v)\| > c\ .
\end{equation}
On the other hand, for any $y\in B$, $v\in V$ one has 
$\mu(0_y^E,v)=h(y)+l(y)\in \{0^F_y\}\times W$, 
hence
$$\mu_{c,\pi}(0_y^E,v)=\rho_{c,F\times W}(\mu(0_y^E,v))=\mu(0_y^E,v)=h(y)+l(v)
\ .
$$
\qed

 \subsection{Compatibility properties}
 \begin{lm}\label{compat} (Coherence Lemma) Let $0<\eta<\frac{c}{4}$,   let $\pi=(F,W)$, $\tilde \pi=(\tilde F,\tilde W)$ be two $\eta$-admissible pairs with $\pi\subset \tilde \pi$, and let $F'$, $W'$ be the orthogonal complements of $F$, $W$ in $\tilde F$, $\tilde W$ respectively. 
 The map
 $$\mu_{c,\pi,\tilde\pi}:=\iota\circ\left\{\left[\mu_{c,\pi}\circ (\p_E,\p_V)\right]\wedge_B \left[(\p_{F'},\p_{W'})\circ(\delta, l)\right]^+_B\right\}:\tilde E\times\tilde V\to [\tilde F\times \tilde W]^+_B
 $$
satisfies properties {\bf P1}, {\bf P2} with constants $C$, $\gamma$ (for a sufficiently small $\gamma$ with $0<\gamma<c$), $\varepsilon_0$, and one has   $\{\mu_{c,\tilde \pi}\}=\{\mu_{c,\pi,\tilde\pi}\}$.
\end{lm}
\pf The first statement follows from Proposition \ref{smashidnew}.  We use the same method as in the proof of Lemma 2.3 in \cite{BF} to construct  a homotopy between the restriction of the   two maps to the product $D_C(\tilde E)\times D_C(V)$ and we will apply the homotopy  invariance property of our invariant (see Proposition \ref{homotopy}).  The main difference compared to  \cite{BF} is that we have to control the restriction to the $S^1$-fixed point set, but we do not need an extension of the homotopy to the whole   
$\tilde E\times\tilde V$. For completeness  we include detailed arguments adapted to our  situation. \vspace{2mm}\\ 
\pf    Denote by $E'$, $V'$ the orthogonal complements of $E$, $V$ in $\tilde E$, $\tilde V$. We define the map
\begin{equation}\label{H}
H:[0,4]\times [D_C(\tilde E)\times D_C(\tilde V)]   {\map} \left[{\cal F}\times{\cal W}\right]\setminus \left[\tilde F^\bot\times \tilde W^\bot\setminus \cringle{D}_c(\tilde F^\bot\times \tilde W^\bot)\right]
\end{equation}
 by the formula \footnote{The third branch of the homotopy was omitted in \cite{BF}.}
 \vspace{1mm}\\
$$H_t=\left\{\begin{array}{ccc}
d+\left[(1-t)\ \id_{{\cal F}\times{\cal W}}+ t\ {\rm p}_{F\times W}\right]\circ k&\rm for& 0\leq t\leq 1\ ,\\
d+\p_{F\times W}\circ k\circ \left[(2-t)\ \id_{\tilde E\times \tilde V}+(t-1)\ \p_{E\times  V}\right]&\rm for& 1\leq t\leq 2\ ,\\
\p_{F\times W}\circ k\circ \p_{E\times V}+\left[d-(t-2)\ \p_{F\times W}\circ d\circ \p_{E'\times V'}\right]&\rm for& 2\leq t\leq 3\ ,\\
\p_{F'\times W'}\circ d+\left[(4-t)\ \p_{F\times W}+(t-3)\ \rho_{c,F\times W}\right]\circ\mu\circ\p_{E\times V}&\rm for& 3\leq t\leq 4\ .
\end{array}\right.$$
\vspace{1mm}\\
{\bf Claim:} {\it $H$ is a well defined, continuous, $S^1$-equivariant map over $B$.}
\vspace{2mm} \\ 
This follows from:\vspace{1mm} \\ 
a) For a point $(t,\tilde e,\tilde v)\in [0,4]\times D_C(\tilde E)\times D_C(\tilde V)$, the term $\rho_{c,F\times W}(\mu(\p_{E\times V}(\tilde e,\tilde v)))$ is finite, so the convex combination in the fourth branch is defined and finite. \vspace{1mm}

Indeed,  recall that the retraction $\rho_{c,F\times W}$ is finite on the complement of the  leaf $\left[F^\bot\times W^\bot\right]\setminus  {D}_c(\ F^\bot\times W^\bot)$. Therefore it suffices to note that $k(D_C({\cal E})\times D_C({\cal V}))$ is $\eta$-close to $F\times W$ and $d(E\times V)\subset   F\times  W$,   so the point $\mu(\p_{E\times V}(\tilde e,\tilde v))$ is $\eta$-close to $F\times W$ for   $(\tilde e,\tilde v)\in D_C(\tilde E)\times D_C(\tilde V)$. Therefore 
$$\mu(\p_{E\times V}(\tilde e,\tilde v))\not\in \left[F^\bot\times W^\bot\right]\setminus \cringle{D}_c(\ F^\bot\times W^\bot)\  .$$
b) The formulae given for  the four components of $H$ agree on the 
intersections of their domains. \vspace{1mm}\\ 
c)  $H$ takes values  in $\left[{\cal F}\times{\cal W}\right]\setminus \left[\tilde F^\bot\times \tilde W^\bot\setminus \cringle{D}_c(\tilde F^\bot\times \tilde W^\bot)\right]$.\vspace{1mm}

Indeed,  for $(t,\tilde e,\tilde v)\in [0,4]\times D_C(\tilde E)\times D_C(\tilde V)$ we see as in the proof of a)   that the right hand term of $H_t$ must be   $\eta$-close to $\tilde F\times\tilde W$, so   $H([0,4]\times D_C(\tilde E)\times D_C(\tilde V))$ avoids $\left[F^\bot\times W^\bot\right]\setminus \cringle{D}_c(\ F^\bot\times W^\bot)$. \\  

The map $H$ has the following properties:
\begin{enumerate}
 \item     $H_0$ coincides with the restriction $\resto{\mu}{D_C(\tilde E)\times D_C(\tilde V)}$. 
 \vspace{1mm} 
 \item   $H_4$ coincides with  the map $\mu_{c,\pi,\tilde\pi}$ composed with the inclusion $\tilde F\times\tilde W\hookrightarrow [{\cal F}\times{\cal W}]^+_B\setminus S_c(\tilde F^\bot\times \tilde W^\bot)$.
 \vspace{1mm} 
 \item  One has \begin{equation}\label{fourthnew}
H_t(0^{\tilde E}_y,\tilde v)=h(y)+l(\tilde v)\ ,\ \forall t\in[0,4]\ \ \forall y\in B\  \  \forall 
\tilde v\in D_C(\tilde V)\ .
\end{equation}
 Formula (\ref{fourthnew}) follows from (\ref{kk}) and the fact that $l$ is an
 isometry, so it commutes with orthogonal projections.
\vspace{1mm} 
\item  $H([0,4]\times \partial (D_C(\tilde E)\times D_C(\tilde V))\cap [\tilde F^\bot\times \tilde W^\bot]=\emptyset$.
\vspace{1mm} 

Indeed, for $(\tilde e,\tilde v)\in \partial (D_C(\tilde E)\times D_C(\tilde V))$ we get $\|H_0(\tilde e,\tilde v)\|=\|\mu(\tilde e,\tilde v)\|\geq c$, whereas $\|\mu(\tilde e,\tilde v)\|$ is $\eta$-close to $F\times W\subset \tilde F\times\tilde W$. Moreover, for $t\in[0,1]$ it holds $\|H_t(\tilde e,\tilde v)-H_0(\tilde e,\tilde v)\|=t\|(\p_{F^\bot\times W^\bot}\circ k)(\tilde e,\tilde v)\|\leq \eta$. For $t\geq 2$ we have
$$\p_{F'\times W'}\circ h_t=\p_{F'\times W'}\circ d\ ,
$$
so $H_t(\tilde e,\tilde v)$ can belong to $\tilde F^\bot\times \tilde W^\bot$ only when $\p_{F'\times W'}\circ d(\tilde e,\tilde v)=0$, i.e. when $(\tilde e,\tilde v)\in E\times V$. For such a pair we find
$$\ \ \ \ \ \ \ \ \ \ \ \ \ H_t(\tilde e,\tilde v)=(d+\p_{F\times W}\circ k)(\tilde e,\tilde v)=\mu(\tilde e,\tilde v)-(\p_{F^\bot\times W^\bot}\circ k)(\tilde e,\tilde v)\ \forall t\in[1,3]\ ,
$$
$$H_t(\tilde e,\tilde v)\in \left[\p_{F\times W}(\mu(\tilde e,\tilde v)),\rho_{c,F}(\mu(\tilde e,\tilde v))\right]\ \forall t\in[3,4]\ ,
$$
so $H_t(\tilde e,\tilde v)$ is a non-vanishing vector of $F\times W$ (more precisely a positive multiple of $\p_{F\times W}(\mu(\tilde e,\tilde v))=\mu(\tilde e,\tilde v)-(\p_{F^\bot\times W^\bot}\circ k)(\tilde e,\tilde v)$) for any $t\in[1,4]$.
\end{enumerate}

These properties have the following important consequence:
\\ \\
{\bf Remark:} {\it The composition $\rho_{c,\tilde\pi}\circ H$  is nowhere vanishing on the space}
$$[0,4]\times \left\{ \partial \left[D_C(\tilde E)\times D_C(\tilde V)\right]\cup \left[0^{\tilde E}\times\tilde V\right]\right\}\ .$$
This follows from the fact that the vanishing locus of the retraction $\rho_{c,\tilde\pi}$ is  the leaf $\cringle{D}_c(\tilde F^\bot\times \tilde W^\bot)\subset \tilde F^\bot\times \tilde W^\bot$. 
On the other hand  we have
$$\rho_{c,\tilde\pi}\circ H_0=\resto{\mu_{c,\tilde \pi}}{D_C(\tilde E)\times D_C(\tilde V)}\ ,\ 
\rho_{c,\tilde\pi}\circ H_4=\resto{\mu_{c,\pi,\tilde \pi}}{D_C(\tilde E)\times D_C(\tilde V)}
$$
 It suffices now to apply Proposition \ref{homotopy}.
 \qed

 Using Proposition \ref{smashidnew} and Lemma \ref{compat} we obtain
 \begin{co} Let $\mu:{\cal E}\times{\cal V}\to{\cal F}\times{\cal W}$ be an $S^1$-equivariant map over a compact CW complex $B$  satisfying ${\cal P}1$, ${\cal P}2$, ${\cal P}3$, and let $0<\eta<\frac{c}{4}$.  Fix an orientation $\oo$ of the finite dimensional summand $H$ of ${\cal W}$. The elements
 $$\{\mu_{c,\pi}\}\in  {_{S^1}\alpha}_B^{b-1}(S(E)_{+B},F^+_B)
 $$
 associated with $\eta$-admissible pairs $\pi=(F,W)$ define a unique class 
 $$\{\mu\}\in  \alpha^{b-1}({\rm ind}(\delta))$$
  which depends only on the map $\mu$ and the orientation $\oo$. 
  \end{co}
  
   In particular, using finite dimensional  approximations associated with constants $C'\geq C$ and $0< c'\leq c$ (and parameter $0<\eta<\frac{c'}{4}$), one obtains the same class.\vspace{2mm}\\
\pf Let $\pi=(F,W)$, $\pi_1=(F_1,W_1)$ two $\eta$-admissible pairs. By Lemma
\ref{compat} we can identify the images of the classes $\mu_{c,\pi}$, $\mu_{c,\pi_1}$ in $\alpha^{b-1}({\rm ind}(\delta))$ under the assumption $\pi\subset \pi_1$. The problem is to reduce the general case to this situation. 

By Corollary \ref{admissco} we know that $\eta$-admissibility of $\pi$ is an open condition, i.e. it is stable under small deformations.  On the other hand, by the homotopy property Proposition \ref{homotopy}, the image of the class $\mu_{c,\pi}$ in the group $\alpha^{b-1}({\rm ind}(\delta))$ is  stable under small deormations of $\pi$. Hence it suffices to consider   a generic small deformation $F'$ of the subbundle $F\subset {\cal F}$  which is fiberwise transversal to $F_1$, such that $(F',W)$ is still $\eta$-admissible.  Then we can put $\tilde F:=F'\oplus F_1$, $\tilde W:=W+W_1$ and   apply twice the compatibility  Lemma
\ref{compat}.

\qed
 \begin{pr}\label{vanishingth} Suppose that the restriction $\resto{\mu}{D_C({\cal E})\times D_C(V)}$ is nowhere vanishing. Then $\{\mu\}=0$.
 \end{pr}  
\pf Since $\resto{\mu}{D_C({\cal E})\times D_C(V)}$ is nowhere vanishing, it is easy to see that there exists $\gamma>0$ such that $\|\mu(e,v)\|> \gamma$ for every $(e,v)\in D_C({\cal E})\times D_C(V)$. Indeed, if not there would exist a sequence $(e_n,v_n)\in D_C({\cal E})\times D_C(V)$ such that  $\|\mu(e_n,v_n)\|\to 0$.  Let $K\subset {\cal F}\times{\cal W}$ be a compact subspace which contains $k(D_C({\cal E})\times D_C(V))$.    Since $d=(\delta,l)$ is a continuous family of Fredholm operators, it follows that $d^{-1}(K)\cap [D_C({\cal E})\times D_C(V)]$ is compact.  Therefore $(e_n,v_n)_n$ admits a subsequence which converges in this intersection. The limit will be a vanishing point  of $\mu$, which contradicts the assumption. 

Use now the constant $c':=\min(\gamma,c)$ (instead of $c$)  in the construction of the finite dimensional approximations of $\mu$. The  obtained maps $\mu_{c',\pi}$ are nowhere vanishing  on $D_C(E)\times D_C(V)$, and our assertion follows from the vanishing property Proposition \ref{vanishing} proved in the finite dimensional case.
\qed
\section{Fundamental properties of the cohomotopy invariants}
\subsection{The Hurewicz image of the cohomotopy invariant}\label{hurewicz}
 \subsubsection{The relative Hurewicz morphism}\label{hurewiczm}
 
Let  $B$ be a  compact space, and let $E$, $F$ be Hermitian bundles of ranks $e$, $f$ over $B$. Let $k$ be an integer and $u\in {_{S^1}\alpha}_B^k(S(E)_{+B},F^+_B)$ a stable class. Suppose for simplicity $k\geq 0$. Consider a representative 
$$\varphi:S(E)_{+B}\wedge_B\xi^+_B\to F^+_B\wedge_B[\underline{\R}^k]^+_B\wedge\xi^+_B$$
 of this stable class, where $\xi=\eta\oplus\xi_0$ is the direct sum of a complex vector bundle  $\eta$ and a real vector bundle $\xi_0$. We may suppose that the real summand $\xi_0$ of $\xi$ is orientable. We  choose an orientation of  $\xi_0$; in this way all our bundles become oriented bundles. The space $S(E)_{+B}\wedge_B\xi^+_B$ can be identified with the fiberwise quotient   ${\{S(E)\times_B  \xi^+_B\}}{/_B 
 \{S(E)\times_B \infty_\xi\}}$. Composing $\varphi$ with the canonical projection one obtains a map of pairs over $B$
 $$\tilde \varphi:(S(E)\times_B \xi^+_B, S(E)\times_B \infty_\xi)\to ([F\oplus \underline{\R}^k\oplus \xi]^+_B,\infty_{F\oplus \underline{\R}^k\oplus \xi})\ .
 $$

Consider now the projection $\pi:\P(E)\to B$ and the following bundles over $\P(E)$:
$$\tilde F:=\pi^*(F)(1)\ ,\ \tilde\xi:=\pi^*(\eta)(1)\oplus\pi^*(\xi_0)\ .
$$
The map $\tilde \varphi$ descends to a morphism of pointed sphere bundles over $\P(E)$ 
$$\bar \varphi:\tilde \xi^+_{\P(E)}\map   [\tilde F\oplus \underline{\R}^k\oplus \tilde\xi ]^+_{\P(E)}\ .
$$
Denote by $s$ the real rank of $\xi$.  Let 
$${\mathrm t}_{\tilde\xi}\in H^{s}(\tilde \xi^+_{\P(E)},\infty_{\tilde\xi};\Z)\ ,\ \t_{\tilde F\oplus \underline{\R}^k\oplus \tilde\xi}\in H^{2f+k+s}([\tilde F\oplus \underline{\R}^k\oplus \tilde\xi ]^+_{\P(E)}, \infty_{\tilde F\oplus \underline{\R}^k\oplus \tilde \xi};\Z)$$
be  the Thom classes of the oriented bundles $\tilde\xi$, $\tilde F\oplus \underline{\R}^k\oplus \tilde\xi$.   The formula
$$\bar \varphi^*(\t_{\tilde F\oplus \underline{\R}^k\oplus \tilde\xi})=\p_{\P(E)}^*(h_{\bar\varphi})\cup \t_{\tilde\xi}
$$
defines a cohomology class $h_{\bar \varphi}\in H^{2f+k}(\P(E);\Z)$ which is independent of the chosen orientation  of $\xi_0$  and of the representative $\varphi$ of the stable class $u$.  For $k\leq 0$ one has a similar construction, but uses  a $[\R^{-k}]^+_B$ factor on the left side.  

The assignment $u=[\varphi]\mapsto h_{\bar\varphi}$ defines a morphism
$$h: {_{S^1}\alpha}^k_B(S(E)_{+B},F^+_B)\to H^{2f+k}(\P(E);\Z)\ ,
$$
which we call  the {\it relative Hurewicz morphism} over $B$. 

Denote by $q:\tilde \xi\to \P(E)$ the bundle projection,  and by $\cringle{\varphi}$ the section  in the pull-back $[q^*(\tilde F\oplus \underline{\R}^k\oplus \tilde\xi)]^+_{\tilde\xi}$ over $\tilde\xi$ defined by $\bar \varphi$. Since the vanishing locus $Z(\cringle{\varphi})$ of this  section is compact, one can define its {\it localized Euler class} class $[\cringle{\varphi}]\in H_{d+2e-2-2f-k}(\tilde\xi;\Z)$, which coincides with the fundamental class $[Z(\cringle{\varphi})]$ of the compact oriented submanifold $[Z(\cringle{\varphi})]$ when $\cringle{\varphi}$ is smooth and transversal to the zero section \cite{Br}.
\begin{re} \label{geomint} (The geometric interpretation of the Hurewicz morphism) Suppose that $B$ is an oriented $n$-dimensional compact manifold.   Then 
$$PD_{\P(E)}(h(u))=[\iota_*]^{-1}([\cringle{\varphi}])\ ,$$ 
where  
$$\iota_*:H_{n+2e-2-2f-k}(\P(E);\Z)\to H_{n+2e-2-2f-k}(\tilde\xi;\Z)\ .$$
is the isomorphism induced by the zero section of $\tilde\xi$. 
  If $\cringle \varphi$ is smooth and transversal to the zero section, then   
$$PD_{\P(E)}(h(u))=[\iota_*]^{-1}([Z(\cringle{\varphi})])\ .
$$
\end{re}
\pf The localized Euler class $[\cringle{\varphi}]\in H_{n+2e-2-2f-k}(\tilde\xi;\Z)$ is defined as the cap product $\cringle{\varphi}^*(\t_{q^*(\tilde F\oplus \underline{\R}^k\oplus \tilde\xi)})\cap [\tilde\xi]$, where $[\tilde\xi]$ stands for the  fundamental class of $\tilde\xi$ in cohomology with compact supports \cite{Br}. We get
$$[\cringle{\varphi}]:=\cringle{\varphi}^*(\t_{q^*(\tilde F\oplus \underline{\R}^k\oplus \tilde\xi)})\cap [\tilde\xi]={\bar \varphi}^*(\t_{\tilde F\oplus \underline{\R}^k\oplus \tilde\xi})\cap [\tilde\xi]=[\p_{\P(E)}^*(h(u))\cup\t_{\tilde\xi}]\cap [\tilde\xi]=$$
$$=\p_{\P(E)}^*(h(u))\cap \iota_*([\P(E)])=
\iota_*(h(u)\cap [\P(E)])=\iota_*(PD_{\P(E)}(h(u))\ .
$$
\qed

Let $\nu=(i,E_1):E\to E'$ be a morphism in the category ${\cal U}_B$ of complex vector bundles over $B$ (see section \ref{alpha}).  Such a morphism induces an isomorphism $E'\cong E\oplus E_1$.   The complement $\P(E')\setminus \P(E_1)$ can be identified with the total space of the complex vector bundle $\pi^*(E_1)(1)\to \P(E)$.   Multiplication with the Thom class $\t_{\pi^*(E_1)(1)}$   defines a morphism
$$H^*(\P(E);\Z)\map H^{*+2e_1}({\pi^*(E_1)(1)}^+_{\P(E)},\infty_{\pi^*(E_1)(1)};\Z)\cong$$
$$\cong H^{*+2e_1}(\P(E'),\P(E_1);\Z)\map H^{*+2e_1}(\P(E');\Z)\ ,
$$
which will be denoted by $a_\nu$.  \\ \\
Now fix an element $x\in K(B)$. A morphism $\tau=(i,j;E_1,F_1,l):(E,F)\to (E',F')$
 in the category ${\cal T}(x)$ defines morphisms
$$a_{(i,E_1)}:H^{2f+k}(\P(E);\Z)\to H^{2f'+k}(\P(E');\Z)\ , $$
$$\P(i)_*:H_k(\P(E);\Z)\to H_k(\P(E');\Z) \ .
$$
For an integer $k\in\Z$ we define
$$H^k(x;\Z):=\varinjlim\hskip-24pt\raisebox{-10pt}{$\scriptstyle (E,F)\in x$} H^{2f+k}(\P(E);\Z)\ ,\ H_k(x;\Z):=\varinjlim\hskip-24pt\raisebox{-10pt}{$\scriptstyle (E,F)\in x$} H_{k}(\P(E);\Z)\ .
$$
Using the same methods as in sections \ref{definition}, \ref{alpha} (stabilizing first with respect to trivial bundle enlargements) we see that these  inductive limits exist in  ${\cal A}b$.
\begin{re}
\begin{enumerate}
\item 
One has $H_*(x;\Z)=H_*(B;\Z)\otimes \Z[t]$. 
\item For a compact $n$-dimensional CW complex $B$ there exist isomorphisms
$$H^k(x;\Z)\simeq\bigoplus_{\begin{array}{c}
\scriptstyle s-k\in 2\Z \vspace{-4pt}\\ 
\scriptstyle\max(0,k-2\iota(x)+2)\leq s\leq n
\end{array}
} H^s(B;\Z)\ ,$$
where $\iota(x)\in\Z$ is the index of $x$. In particular, putting $n(x):=2\iota(x)-2+n$, one has $H^{n(x)}(x;\Z)=H^n(B;\Z)$.
\end{enumerate}
\end{re}

The integer  $n(x):=2\iota(x)-2+n$ will be called  {\it  the dimension of the  formal 
projectivization} of $x$.

 \begin{re} Suppose that $B$ is a compact  connected oriented manifold of dimension $n$. The system of Poincar\'e duality isomorphisms $PD_{\P(E)}$  defines isomorphisms
 $$PD_x:H^k(x;\Z)\textmap{\simeq} H_{n(x)-k}(x;\Z)\ . $$
 \end{re}
 \begin{re} The system of Hurewicz morphisms 
 $$h:{_{S^1}\alpha}_B^k(S(E)_{+B},F^+_B)\to H^{2f+k}(\P(E);\Z)$$
  defines  a morphisms of graded groups $h_x:\alpha^*(x)\to H^*(x;\Z)$. If $B$ is  a compact  connected oriented manifold, one also gets a morphism   $PD_x\circ  h_x:\alpha^*(x)\to H_*(x;\Z)$, which we call the homological Hurewicz morphism.
 \end{re}
 The result below has the following important consequence: for a moduli problem with vanishing ``expected dimension", the cohomotopy invariant yields the same information as the classical (co)homological  invariant.  Recall that our cohomotopy invariant $\{\mu\}$ associated with a map satisfying properties ${\cal P}1$ -- ${\cal P}3$ belongs to $ \alpha^{b-1}(x)$, where $x:={\rm ind}(\delta)$, $b:=\dim(H)$ (see section \ref{general}). The {\it expected dimension} $w(\mu):=2\iota(x)+\dim(B)-b-1$ of the moduli problem associated with $\mu$  vanishes if and only if $b-1=n(x)$. 
 \begin{pr} Suppose that $B$ is a finite CW complex of dimension $n$. Then the Hurewicz morphism 
 $$h_x^{n(x)}: \alpha^{n(x)}(x)\map H^{n(x)}(x;\Z)=H^n(B;\Z)\ .
 $$
 is an isomorphism.
 \end{pr}
\pf  Suppose $n(x)\geq 0$ for simplicity. Fix a stabilizing bundle $\xi$. Using the same method and the same notations as in section \ref{hurewiczm} we see that the set 
$${_{S^1}\pi}^0(S(E)_{+B}\wedge_B\xi^+_B, F^+_B\wedge_B[\underline{\R}^{n(x)}]^+_B\wedge\xi^+_B)
$$
can be identified with the set of pointed bundle maps
$$\bar \varphi:\tilde \xi^+_{\P(E)}\map   [\tilde F\oplus \underline{\R}^{n(x)}\oplus \tilde\xi ]^+_{\P(E)}
$$
 over $\P(E)$. The latter set can be identified with  $H^{\dim_\R(\P(E)}(\P(E);\Z)=H^{n}(B;\Z)$ by Proposition \ref{obsbundle} via the map $\bar \varphi\mapsto h_{\bar\varphi}$. The obtained bijections  
 $${_{S^1}\pi}^0(S(E)_{+B}\wedge_B\xi^+_B, F^+_B\wedge_B[\underline{\R}^{n(x)}]^+_B\wedge\xi^+_B)\simeq H^{n}(B;\Z)$$
are compatible with morphisms $\xi\to \xi'$ in the category ${\cal C}_B$ and  with morphisms $(E,F)\to (E',F')$ in the category ${\cal T}(x)$. Therefore we get a bijection $\alpha^{n(x)}(x)\to H^n(B;\Z)$, which coincides with the Hurewicz map  by the definition.
\qed

\subsubsection{A comparison theorem} The main result of this section states: the virtual fundamental class of the moduli space of solutions associated with a map   $\mu$ satisfying  properties ${\cal P}1$, ${\cal P}2$, ${\cal P}3$ can be identified with the image of the cohomotopy invariant under the homological Hurewicz map.   Applied to   Seiberg-Witten theory,  this implies that the {\it full} Seiberg-Witten type invariant coincides with the Hurewicz image of the cohomotopy Seiberg-Witten  invariant. \\

We begin with the finite dimensional case. Let  $B$ be a compact oriented manifold, $p:E\to B$, $q:F\to B$ Hermitian bundles over $B$, let $V$, $W$ be Euclidean spaces, and  let $\mu:{E}\times{ V}\to [{F}\times{W}]^+_B$  be an $S^1$-equivariant  map over $B$ satisfying properties {\bf  P}1, ${\bf  P2}$ of section \ref{cylinder}. The invariant $\{\mu\}\in {_{S^1}\alpha}_B^{b-1}(S(E)_{+B},F^+_B)$ is defined by a map of pairs
$$ 
\left(S(E)\times D_R(\R\oplus V),
S(E)\times S_R(\R\oplus V) \right){\to}  ([F
\times W]^+_B,\left[F\times W\right]^+_B\setminus
\cringle{D}_{\varepsilon}(F
\times W) ) 
$$
induced by the restriction $\mu_{R,\varepsilon}:D_R(E)\times D_R(V)\to (F\times W)^+_B$ of $\mu$ to a sufficiently large cylinder $D_R(E)\times D_R(V)$.
The vanishing locus of $\mu$ (regarded as section in the bundle  $(p^*(F)\times V)\times W\to E\times V$) is an $S^1$-invariant compact space contained in the open subspace $\cringle{D}_R(E)\times \cringle{D}_R(V)\setminus [0^E\times_B D_R(V)]$ of the  cylinder.  Its $S^1$-quotient can be identified with the vanishing locus of the section $\cringle{\mu}_{R,\varepsilon}$ induced by $\mu_{R,\varepsilon}$ on the $S^1$-quotient $\P(E)\times \cringle{D}_R(\R\oplus V)$ of $S(E)\times \cringle{D}_R(\R\oplus V)$. Using Remark  \ref{geomint} one obtains
\begin{co}\label{firstcomp} Suppose that $B$ is a compact oriented manifold. Via the isomorphism $H_*(\P(E)\times \cringle{D}_R(\R\oplus V);\Z)\simeq H_*(\P(E);\Z)$ the Poincar\'e dual $PD_{\P(E)}(h(\{\mu\}))$ coincides with the virtual fundamental class associated with the section   $\cringle{\mu}_{R,\varepsilon}$. If this section is smooth and transversal to the zero section, then  $PD_{\P(E)}(h(\{\mu\}))$ can be identified with the fundamental class of the vanishing locus $Z(\cringle{\mu}_{R,\varepsilon})\subset \P(E)\times \cringle{D}_R(\R\oplus V)$.
\end{co}
Note that $\mu$ is nowhere vanishing outside the cylinder $D_R(E)\times D_R(V)$, so the vanishing loci of $\mu$ and ${\mu}_{R,\varepsilon}$ can be identified. The vanishing locus $M:=Z(\cringle{\mu}_{R,\varepsilon})\cong  {Z(\mu)}/{S^1}$ will be called the {\it ``moduli space"} associated with the map $\mu$.\\

 Let $p:{\cal E}\to B$, $q:{\cal F}\to B$ be complex Hilbert bundles over  $B$, let ${\cal V}$, ${\cal W}$ be  real Hilbert spaces, and let $\mu:{\cal E}\times{\cal V}\to {\cal F}\times{\cal W}$  be an   $S^1$-equivariant   map over $B$ satisfying properties ${\cal P}1$, ${\cal P}2$, ${\cal P}3$ in section \ref{general}.    Denote by $\pi:\P({\cal E})\to B$   the natural projection.  The map $\mu_{R,\varepsilon}$ descends to a smooth section $\cringle{\mu}_{R,\varepsilon}$ in the bundle 
 $$\pi^*({\cal F})(1)\times \cringle{D}_R(\R\oplus {\cal V})\times {\cal W}\to\P({\cal E})\times \cringle{D}_R(\R\oplus {\cal V})\ ,$$
 and again one can identify the moduli space ${\cal M}:=Z(\mu)/S^1$ of $\mu$ with the vanishing locus $Z( \cringle{\mu}_{R,\varepsilon})$ of this section.  Using the same argument as in the proof of Proposition \ref{vanishingth}, we see that the moduli space ${\cal M}$ is compact. 
Suppose now that 
\begin{enumerate}
\item[${\cal P}_4$:] $B$ is a compact,  smooth,  connected, oriented manifold, $\mu$ is smooth and the fiberwise differential of $k:=\mu-d$ at any point is a compact operator.
\end{enumerate}

This condition is always satisfied in practical gauge theoretical situations; indeed, the map $k$ is usually given by the composition of a smooth map ${\cal E}\times{\cal V}\to {\cal F}_1\times{\cal W}_1$ with a map ${\cal F}_1\times{\cal W}_1\to {\cal F}\times{\cal W}$ over $B$   defined by a smooth family of compact operators. The condition ${\cal P}_4$  implies that  $\cringle\mu_{R,\varepsilon}$ is a smooth Fredholm section on the Banach manifold $\P({\cal E})\times \cringle{D}_R(\R\oplus {\cal V})$. In order to give sense to the virtual fundamental class  of the moduli space ${\cal M}$ we have to trivialize the determinant line bundle $\det({\rm index}(D\cringle\mu_{R,\varepsilon}))$ over ${\cal M}$. Equivalently, it suffices to trivialize the line bundle $\det({\rm index}(D \mu))$ over $Z(\mu)$. In these formulae the symbol $D$ stands for the  family of intrinsic derivatives of a section at its zero locus, and $\mu$ is regarded as a section of the bundle $[p^*({\cal F})\times {\cal V}]\times {\cal W}\to {\cal E}\times {\cal V}$. For a point $(e,v)\in Z(\mu)$ with $p(e)=y$ one has a natural identification
$$\det({\rm index}(D_{(e,v)} \mu))=\Lambda^n(T_y(B))\otimes\det({\rm index} (d_{(e,v)}\resto{\mu}{{\cal E}_y\times {\cal V}}))\ ,
$$
where $n:=\dim(B)$ and $\resto{\mu}{{\cal E}_y\times {\cal V}}:{\cal E}_y\times {\cal V}\to {\cal F}_y\times {\cal W}$ is  the restriction of $\mu$ to the fiber over $y$. By the condition ${\cal P}4$, the differential of this restriction is congruent with the operator $d_y=(\delta_y,l)$ modulo a compact operator. Therefore (since the family $\delta=(\delta_y)_{y\in B}$ has a canonical complex orientation, and $B$ is oriented) one obtains a trivialization of  $\det({\rm index}(D \mu))$ for every orientation $\oo$ of $\coker(l)=H$.  This is precisely the orientation parameter involved  in the definition of the cohomotopy invariant $\{\mu\}$. Fix such an orientation $\oo$. Using the results in \cite{Br}, we obtain a virtual fundamental class in C\v{e}ch homology $[{\cal M}]^{\rm vir}\in \check{H}_{w}({\cal M};\Z)$, where
 $$w=w(\mu):=n+2\iota({\rm ind}(\delta))-b-1=n({\rm ind}(\delta))-(b-1)$$
  is the expected dimension of our moduli problem (the  index of the section $\cringle{\mu}_{R,\varepsilon}$).\\
  
Put $x:={\rm ind}(\delta)$, and note that the group
 $$H_{w}(x;\Z)=\bigoplus\limits_{0\leq 2i\leq w} H_{w-2i}(B;\Z)\otimes t^{i}\ .
 $$
 can be identified with  $H_{w}(\P({\cal E})\times \cringle{D}_R(\R\oplus V);\Z)=H_{w}(\P({\cal E});\Z)$.
\begin{dt} The {\it full   homological invariant} of $\mu$ is the image $\{\mu\}_{\mathrm H}$ of the   class $[{\cal M}]^{\rm vir}$ in the group $H_{w}({\rm ind}(\delta);\Z)$.
\end{dt}
\begin{thry} Suppose that conditions ${\cal P}1$ -- ${\cal P}4$ hold.  Then
$$\{\mu\}_{\mathrm H}=PD_x\circ h_x(\{\mu\})\ .$$
\end{thry}
\pf As in section \ref{finite} choose a finite dimensional approximation $\mu_{c,\pi}$ of $\pi$, associated with an $\eta$-admissible pair $(F,W)$. Define   $\mu_{c,\pi,\infty}:D_C({\cal E})\times D_C({\cal V})\to {\cal F}\times {\cal W}$ by
$$\mu_{c,\pi,\infty}(e,v)=\mu_{c,\pi}(\p_E(e),\p_V(v))+ \p_{F^\bot\times W^\bot}\circ d\circ \p_{E^\bot\times V^\bot}\ .
$$
This map takes finite values by Lemma \ref{finapprox}. We claim that there exists a smooth homotopy 
$${\cal H}:[0,4]\times D_C({\cal E})\times D_C({\cal V})\to {\cal F}\times {\cal W}$$
 between $\resto{\mu}{D_C({\cal E})\times D_C({\cal V})}$ and $\mu_{c,\pi,\infty}$   in the space of $S^1$-equivariant Fredholm maps over $B$, such that for $0\leq t\leq 4$ the map ${\cal H}_t$ has no zeroes  in $\partial [D_C({\cal E})\times D_C({\cal V})]\cup 0^{\cal E}\times D_C(V)$. To obtain such a  homotopy it suffices to replace $\tilde E$, $\tilde V$, $\tilde F $, $\tilde W$  in the definition of the homotopy $H$ used in the proof of  Lemma \ref{compat} by ${\cal E}$, ${\cal V}$, ${\cal F}$, ${\cal W}$,  and to compose the resulting map from the right with a smooth homeomorphisms $\theta:[0,4]\to [0,4]$ having the properties
 $$\theta(i)=i\ ,\ \theta^{(k)}(i)=0\hbox{  for }i\in\{0,1,2,3,4\}\ , \ k\geq 1$$
  (to assure differentiability).  Using the homotopy  invariance  of the virtual class \cite{Br}, we can identify $\{\mu\}_{\mathrm H}$ with the image of the virtual class $[\mu_{c,\pi,\infty}]^{\rm vir}$ in $H_w(\P({\cal E});\Z)$. On the other hand, by the ``associativity property" of the virtual class  (see Proposition 14  (4)  in \cite{Br}) and Corollary \ref{firstcomp}, the latter is just the image of $PD_{\P(E)}(h(\{\mu_{c,\pi}\})$ via the embedding $\P(E)\to \P({\cal E})$. But $PD_{\P(E)}(h(\{\mu_{c,\pi}\})$ is a representative of $PD_x\circ h_x(\{\mu\})$.

\qed

\subsection{Cohomotopy invariant jump formulae} \label{cwc}  
\subsubsection{General results}
 Let 
$$M\map N\map P
$$
be  a cofiber  sequence of pointed $S^1$-spaces over a compact basis $B$. 
For every pointed $S^1$-space $Y$ over $B$ there is an associated 
long exact sequence of cohomotopy groups
\begin{equation}\label{longexact}
\dots \to {_{S^1}\alpha}^k_B(P,Y)\to {_{S^1}\alpha}^k_B(N,Y)\to 
{_{S^1}\alpha}^k_B(M,Y)\textmap{\partial} 
{_{S^1}\alpha}^{k+1}_B(P,Y)\to \dots\ \ .
\end{equation}
The connecting morphism 
$$\partial:{_{S^1}\alpha}^k_B(M,Y)={_{S^1}\alpha}^{k+1}(M\wedge_B S^1,Y)\map 
{_{S^1}\alpha}^{k+1}(P,Y)
$$
is given by composition with the contraction map $\cg:P\to M\wedge_B\underline{S}^1$ induced 
by a fixed homotopy equivalence between $P$  and the mapping cone of the 
map $M\to N$.  For the cofiber sequence 
$$S(\xi)_{+B}\map  {D}(\xi)_{+B}\map \xi^+_B
$$
associated with a vector bundle $\xi$ over a compact basis $B$, the morphism $\partial$ can be described in the following way. The obvious isomorphisms
$$S(\xi)_{+B}\wedge_B\underline{S}^1\cong \qmod{S(\xi)\times [0,1]}{S(\xi)\times\{0,1\}}\ ,\ 
\xi^+_B\cong \qmod{S(\xi)\times [0,1]}{\sim}
$$
(where $\sim$ is the equivalence relation generated by  $(v,0)\sim (v',0)$, $(v,1)\sim (v',1)$) allow us to use ${S(\xi)\times [0,1]}/{S(\xi)\times\{0,1\}}$,  $ {S(\xi)\times [0,1]}/{\sim}$ as models for 
$S(\xi)_{+B}\wedge_B\underline{S}^1$ and $\xi^+_B$. Using these models, the morphism $\partial$ is given by composition with the contraction map 
\begin{equation}\label{cgxi}
\cg_\xi: \qmod{S(\xi)\times [0,1]}{\sim}\map  \qmod{S(\xi)\times [0,1]}{S(\xi)\times\{0,1\}}
\end{equation}
induced by the identity of $S(\xi)\times [0,1]$.  

Consider now an oriented $b$-dimensional real vector space $H$ and the 
 cofiber sequence over $B$ associated with the trivial bundle $\underline{H}=B\times H$ over $B$:
$$S(\underline{H})_{+B})\map  {D}(\underline{H})_{+B}\map \underline{H}^+_B\ .
$$
Let $E$ be a Hermitian  vector bundle over $B$. Taking  smash product with $S(E)_{+B}$ over $B$ yields the following cofiber sequence 
over $B$
$$S(E)_{+B}\wedge_B{S}(\underline{H})_{+B}\to
S(E)_{+B}\to S(E)_{+B}\wedge_B \underline{H}^+_B
$$
Since $S(E)_{+B}\wedge_B {S}(\underline{H})_{+B}=[S(E)\times S(H)]_{+B}$, the associated long exact cohomotopy   sequence is
\begin{equation}\label{longexactnew}
\dots\to {_{S^1}\alpha}^{-1}_{B}(S(E)_{+B}\wedge_B\underline H^+_B,
[F\oplus\underline{H}]^+_{B})\to{_{S^1}\alpha}^{-1}_{B}(S(E)_{+B}, 
[F\oplus\underline{H}]^+_{B})\to$$
$$\to{_{S^1}\alpha}^{-1}_{B}([S(E)\times S(H)]_{+B}, 
[F\oplus\underline{H}]^+_{B})\stackrel{\partial}{\to}$$
$$\to{_{S^1}\alpha}^0_{B}(S(E)_{+B},
[F]^+_{B})\to {_{S^1}\alpha}^{0}_{B}(S(E)_{+B}, [F\oplus\underline{H}]^+_{B})\to 
\dots\ .
\end{equation}

Note that one has canonical base change isomorphisms 
\begin{equation}\label{duality}
{_{S^1}\alpha}^{k}_{B}([S(E)\times S(H)]_{+B}, 
[F\oplus\underline{H}]^+_{B})\simeq 
{_{S^1}\alpha}^{k}_{\tilde B}(S(\tilde E)_{+\tilde B}, 
[\tilde F\oplus\underline{H}]^+_{\tilde B})\ .
\end{equation}
associated with the projection 
$$p:\tilde B=B\times S(H)\to B$$
(see \cite{CJ} Proposition 5.37, Proposition 12.40 for the non-equivariant case).
\vspace{5mm}

A map $\kappa:B\to S(H)$ defines a section $j_\kappa^E: S(E)_{+B}\to [S(E)\times S(H)]_{+B}$ over $B$ of the projection $[S(E)\times S(H)]_{+B}\to S(E)_{+B}$, so it defines a splitting of the exact sequence (\ref{longexactnew}).  
\begin{lm}\label{differenceclass} Let $m\in {_{S^1}\alpha}^{-1}_{B}([S(E)\times S(H)]_{+B}, 
[F\oplus\underline{H}]^+_{B})$, and  let $\kappa_0,\ \kappa_1:B\to S(H)$ be two maps. One has the identity
$$(j_{\kappa_1}^E)^*(m)-(j_{\kappa_0}^E)^*(m)=d(\kappa_0,\kappa_1)\cdot \partial(m)\ ,
$$
where where $d(\kappa_0,\kappa_1)\in {_{S^1}\alpha}^{-1}_B(B_{+B}, \underline{H}^+_B)={_{S^1}\alpha}^{b-1}_B(B_{+B}, B_{+B})$ is the difference class of the maps $\kappa_0$, $\kappa_1$ regarded as sections in the sphere bundle $S(\underline{H})$.
\end{lm}
\pf The difference class $d(\kappa_0,\kappa_1)$ is defined by the map
$$\Delta:B_{+B}\wedge_B\underline{S}^1=\qmod{B\times[0,1]}{_B B\times\{0,1\}}\map \qmod{D_\epsilon(\underline{H})}{S_\epsilon(\underline{H})}=\underline{H}^+_B
$$
induced    by
$$(b,t)\mapsto\left\{
\begin{array}{ccc}
{\ } [(1-2t)\kappa_0(b)]&{\rm for}& 0\leq t\leq \frac{1}{2}\ \ \\
{\ }[(2t-1)\kappa_1(b)] &{\rm for}& \frac{1}{2}\leq t\leq 1\ .
\end{array}\right.$$ 
The connecting morphism $\partial_H$ in the long exact sequence
$${_{S^1}\alpha}^{-1}_B(B_{+B}, \underline{H}^+_B)\stackrel{\partial_H}{\to} {_{S^1}\alpha}^{0}_B(B_{+B},  {S}(\underline{H})_{+B}) \to {_{S^1}\alpha}^{0}_B(B_{+B},  B_{+B})\to {_{S^1}\alpha}^{0}_B(B_{+B}, \underline{H}^+_B)
$$
is defined via the identifications
 $$ {_{S^1}\alpha}^{-1}_B(B_{+B}, \underline{H}^+_B)={_{S^1}\alpha}^{0}_B(B_{+B}\wedge_B\underline{S}^1, \underline{H}^+_B)$$
 $${_{S^1}\alpha}^{0}_B(B_{+B},  {S}(\underline{H})_{+B})={_{S^1}\alpha}^{0}_B(B_{+B}\wedge_B\underline{S}^1,  {S}(\underline{H})_{+B}\wedge_B\underline{S}^1)\ ,$$
 by left composition with the contraction $\cg_H:\underline{H}^+_B\to  {S}(\underline{H})_{+B}\wedge_B\underline{S}^1$. The   image of $d(\kappa_0,\kappa_1)$ under $\partial_H$ is just the difference 
 $\{\kappa_1\}-\{\kappa_0\}\in {_{S^1}\alpha}^{0}_B(B_{+B},  {S}(\underline{H})_{+B})$.

One has obviously 
$$(j_{\kappa_1}^E)^*(m)-(j_{\kappa_0}^E)^*(m)=m\circ \left( \{\kappa_1\}-\{\kappa_0\}\right)=
m\circ \partial_H (d(\kappa_0,\kappa_1))\ .$$
We know that $\partial_H (d(\kappa_0,\kappa_1))$ is represented by $\cg_H\circ \Delta$ and the connecting operator $\partial$ in the exact sequence (\ref{longexactnew}) acts by right composition with the same contraction $\cg_H$. Therefore
$$(j_{\kappa_1}^E)^*(m)-(j_{\kappa_0}^E)^*(m)=m\circ(\cg_H\circ \Delta)=\partial(m)\circ d(\kappa_0,\kappa_1)=\partial(m)\circ (d(\kappa_0,\kappa_1)\cdot \{\id_{B_{+B}}\})$$
$$=\left(d(\kappa_0,\kappa_1)\cdot\partial(m)\right)\circ \{\id_{B_{+B}}\}=d(\kappa_0,\kappa_1)\cdot\partial(m) \ .
$$
Here we have used the fact that the composition multiplication $\circ$ is ${_{S^1}\alpha}^*(B)$-bilinear.
\qed

This lemma has an important analogue for the groups $\alpha^*(x)$ associated with a K-theory element $x$.  For a compact space $P$  we put 
$$\alpha^*(P;x)=\varinjlim\hskip-24pt\raisebox{-10pt}{$\scriptstyle(E,F)\in x$} {_{S^1}\alpha}_B^*(S(E)_{+B}\wedge_B \underline{P}_{+B},F^+_B)\ .  
$$
where the inductive limit is taken with respect to the category ${\cal T}(x)$.   Using the methods used in section  \ref{alpha} for the definition of the groups $\alpha^*(x)$, and the results in section \ref{limits}, we see that this inductive limit exists; it can be constructed by taking first the limit of ${_{S^1}\alpha}_B^*(S(E\oplus\underline{\C}^n)_{+B}\wedge_B \underline{P}_{+B},[F\oplus\underline{\C}^n]^+_B)$ over $n$,  and  factorizing the result by the action of $\tilde J(I[K^{-1}(B)])\subset {_{S^1}\alpha}^0(B)$.
The graded group $\alpha^*(P;x)$ comes with an obvious  homomorphism $\alpha^*(P;x)\to \alpha^*(\p_B^*(x))$, where $\p_B:B\times P\to B$ is the projection on the first summand.  
 
 Taking the  inductive limit  of the connection morphisms $\partial=\partial_{E,F}$ in  (\ref{longexactnew}) with respect to the category ${\cal T}(x)$,  one gets a morphism
\begin{equation}\label{deltax}
\partial_x:=\varinjlim\hskip-24pt\raisebox{-10pt}{$\scriptstyle(E,F)\in x$}\partial_{E,F}:\alpha^{b-1}(S(H);x){\map}  \alpha^0(x) \ .
\end{equation}
which is intrinsically associated with $x$. 

Let $\kappa:B\to S(H)$ be a fixed map. The system of morphisms  
$$(j_\kappa^E)^*:{_{S^1}\alpha}^{*}_{B}([S(E)\times S(H)]_{+B},  F^+_{B})\to {_{S^1}\alpha}^{*}_{B}(S(E)_{+B},  F^+_{B})$$
 induces a morphism $j_\kappa^*:\alpha^*(S(H);x)\to \alpha^*(x)$.
 \begin{co}\label{differenceclassnew} Let $m\in  \alpha^{b-1}(S(H);x)$, and let $\kappa_0,\ \kappa_1:B\to S(H)$ be two maps. One has the identity
$$(j_{\kappa_1})^*(m)-(j_{\kappa_0})^*(m)=d(\kappa_0,\kappa_1)\cdot \partial_x(m) \ .
$$
 \end{co}

\subsubsection{The universal perturbation and the invariant jump formulae}

 Let $E$, $F$ be Hermitian vector bundles over a compact basis $B$,   let  $V$, 
$W$ be Euclidean vector spaces, and let $\mu:  E\times V\to[F\times W]^+_B$ be 
an $S^1$-equivariant map over $B$ satisfying the properties {\bf P1} and 
{\bf P2} (1)  with $h=0$. In other words, 
$$\mu(0^E_y,v)=l(v)\ ,\ \forall y\in B\ \forall v\in V\ ,
$$
where $l:V\textmap{\simeq}W_0\subset W$ is a linear embedding. The cylinder 
construction cannot be applied to such a map, because $\mu$ has vanishing
points on the core $0^E\times D^R(V)$ of any cylinder 
$D_R(E)\times D_R(V)$. We   orient the orthogonal 
complement $H$ of $W_0$  in $W$, and we denote by $b$ its dimension.
Let $\epsilon>0$. For every map $\kappa:B\to   S_\epsilon(H)$ we define the perturbation 
$$\mu_\kappa: E\times V\to  [F\times W]^+_B$$
 by putting  $\mu_\kappa(e,v):=T_{\kappa(y)}( \mu(e,v))$ for $e\in E_y$. Here  $T_{\kappa(y)}$ denotes
    the automorphism of $[F\times(H\oplus W_0)]^+_B$ which 
extends the translation 
$$(f,w)\mapsto (f,w+\kappa(y))\ .$$
\begin{re} If $\epsilon>0$ is sufficiently small, the map $\mu_\kappa$ satisfies 
the properties {\bf P1}, {\bf P2} of section \ref{cylinder}, so the cylinder construction applies and 
yields a stable class $\{\mu_\kappa\}\in {_{S^1}\alpha}_{B}^{b-1}(S(E)_{+B},F^+_B)$.
\end{re}
\pf Suppose that $\mu$ satisfies  the property {\bf P1} with constants $C$, $c$. Choose 
$\epsilon<\frac{c}{2}$. The map $\mu_\kappa$ 
 satisfies {\bf P1} with constants $C$, $c':=\frac{c}{2}$, and {\bf P2} 
with constant $\varepsilon_0=\epsilon$.
\qed

Another way to construct a map satisfying properties   {\bf P1} , {\bf P2} is to let  $\kappa$  vary  in the sphere $S_\epsilon (H)$ and consider the  {\it universal 
perturbation}
$$\tilde \mu:\tilde E\times V\map \tilde F\times W
$$
over the basis $\tilde B:= B\times S_\epsilon(H)$ (where $\tilde E:=\p_{B}^*(E)$, 
$\tilde F:=\p_{B}^*(F)$) which acts as $\mu_{\kappa}$ over $
B\times\{\kappa\}$. This map also satisfies properties {\bf P1}, {\bf P2} 
with the same constants as any $\mu_\kappa$, so that
the cylinder  construction applies and yields a   class 
$\{\tilde \mu\}\in
{_{S^1} \alpha}_{\tilde B}^{b-1}(S(\tilde E)_{\tilde B},\tilde F^+_{\tilde B})$. Our next goal
 is to understand this class 
$\{\tilde\mu\}$. The essential point is to identify the image of     $\{\tilde \mu\}
\in {_{S^1}\alpha}_{\tilde B}^{b-1}(S(\tilde E)_{+\tilde B},\tilde F^+_{\tilde B})\}$ 
under the connecting morphism $\partial$.   

Recall from section \ref{euler}  that   $\{o_{(E,F)}\}\in {_{S^1}\alpha}^0_{B}(S(E)_{+B},
F^+_{B})$ is the class of the obvious pointed map $S(E)_{+B}\to F^+_{B}$ over $B$  which maps $+_B$ to the infinity section, and  $S(E)$ to the trivial section. 
\begin{pr} \label{cswc} (The $\partial$-image of the  invariant of the universal perturbation) Via the identification 
$${_{S^1}\alpha}^0_{B}(S(E)_{+B}\wedge_B\underline H^+_B,
[F\oplus\underline{H}]^+_{B})={_{S^1}\alpha}^0_{B}(S(E)_{+B},
F^+_{B})
$$
one has 
$$\partial(\{\tilde \mu\})=-\{o_{(E,F)}\}\ .$$
\end{pr}

\pf    As in section \ref{cylinder}  fix $R>C$ and $\varepsilon <
\min(\varepsilon_0,c')=\min(\epsilon,\frac{c}{2})$. Let $\tau_0<R$ 
be sufficiently small such that $\mu(e,v)$ remains 
finite for every $(e,v)\in D_{\tau_0}(E)\times D_R(V)$.  \\
\\
Step 1. We replace $\resto{\tilde\mu}{D_R(\tilde E)\times D_R(V)}$ by a map $\tilde\mu_\tau$ which represents the  same  class $\{\mu\}$ and coincides with the $\kappa$-independent map $\mu$ outside 
the smaller cylinder $D_\tau(E)\times D_R(V)$.\\

Define  $\tilde\mu_{\tau}:D_R(\tilde E)\times D_R(V)\map 
[\tilde F\times W]^+_{\tilde B}$ by the formula
$$\tilde\mu_{\tau}(e,\kappa,v):=\left\{
\begin{array}{ccc}
\left(1-\frac{1}{\tau} \|e\|\right)(\kappa+l(v))+ \frac{1}{\tau}\|e\|\mu(e,v)&\rm for&
0\leq\|e\|\leq\tau\\
\mu(e,v)&\rm for&\|e\|\geq\tau\ .
\end{array}\right.
$$

The maps $\tilde \mu_{\tau}$ and $\tilde \mu$ coincide on the core 
$0^{\tilde E}\times D_R(V)$ of the cylinder $D_R(\tilde E)\times D_R(V)$ and 
they differ by the translation $T_\kappa$ outside
$D_\tau(\tilde E)\times D_R(V)$. We define a homotopy between $\tilde\mu_{\tau}$
 and $\resto{\tilde \mu}{D_R(\tilde E)\times D_R(V)}$ by putting
$$\tilde\mu_{\tau}^t(e,\kappa,v):=\left\{
\begin{array}{ccc}
(1-t)\tilde\mu_{\tau}(e,\kappa,v)+t\tilde\mu(e,\kappa,v)&\rm for&
\|e\|\leq\tau\\
T_{t\kappa}\circ  \mu(e,v)&\rm for&\|e\|\geq\tau\ .
\end{array}\right.
$$
{\bf Claim:} If $\tau$ is sufficiently small, then 
$\|\tilde\mu^t_\tau\|\geq c'$ on $\partial \left[D_R(\tilde E)\times D_R(V) 
\right]$ for every $t\in [0,1]$.\\

The claim is not obvious only for points $(e,v)\in 
D_\tau(\tilde E)\times S_R(V)$. One has the identity
$$\tilde\mu_{\tau}^t(e,\kappa,v)=(1-t)
\left\{\left(1-\frac{1}{\tau} \|e\|\right)\kappa+l(v)+\frac{1}{\tau}
\|e\|\left[\mu(e,v)-l(v)\right]\right\}+t\mu(e,v)+t\kappa=
$$
$$=l(v)+\left(1-\frac{1-t}{\tau}\|e\|\right)\kappa+
\left[t+\frac{(1-t)}{\tau}
\|e\|\right]\left[\mu(e,v)-l(v)\right]\ .
$$

The first two terms belong to orthogonal complements, so for $e\in 
D_\tau(E)$ one has
$$\|\tilde\mu_{\tau}^t(e,\kappa,v)\|\geq \|l(v)\|-\|\mu(e,v)-l(v)\|\ .
$$
Since $\mu(0^E_y,v)=l(v)$, and $\mu$ is fiberwise differentiable with 
globally continuous derivatives on $E\times V$, it holds
$$\lim_{\tau\to 0}\big\{\sup
\left\{ (\mu(e,v)-l(v))\ \vline\ \ 0\leq\|e\|\leq\tau,\
 \|v\|\leq R\right\}\big\}
=0\ .
$$
On the other hand, for $\|v\|=R$ one has $\|l(v)\|=\|\mu(0^E_y,v)\|>c$. This 
proves the claim.
\\ \\
Using the Claim and 
$\|\tilde\mu_{\tau}^t(e,\kappa,v)\|=\|\kappa\|=\epsilon>0$ we see that 
$(\tilde\mu_{\tau}^t)_{t\in[0,1]}$ defines a homotopy between 
$\tilde \mu_\tau$ and $\resto{\tilde\mu}{D_R(\tilde E)\times D_R(V)}$ in the
space of maps for which the cylinder construction applies. Therefore
\begin{equation}\label{first-step}
\{\tilde \mu\}=\{\tilde\mu_\tau\}\in {_{S^1}\alpha}^{b-1}_{B}([S(E)\times S(H)]_{+B}, 
F^+_{B}) \hbox{ for all sufficiently small }\tau>0.
\end{equation}
\\
Step 2. We  compute the class $-\partial(\{\tilde\mu_\tau\})$.  \\

Regard $\{\tilde \mu_\tau\}$ as an element in the group
$${_{S^1}\alpha}^{-1}_{B}([S(E)\times S(H)]_{+B}, 
[F\oplus\underline{H}]^+_{B})={_{S^1}\alpha}^{0}_{B}([S(E)_{+B}\wedge_B \underline{S}(H)_{+B}
\wedge_B \underline{S}^1, 
[F\oplus\underline{H}]^+_{B})\ .
$$
As explained at the beginning of this section the morphism $\partial$ is given by composition with the contraction map
$$\cg_H:H^+=\qmod{S_{\epsilon}(H)\times [0,R]}{\sim_H}\to S(H)_+\wedge S^1=
\qmod{S_\epsilon(H)\times[0,R]}{S(H)\times\{0,R\}} $$
induced by the identity of $S_{\epsilon}(H)\times [0,R]$. The morphism $-\partial$ is defined by composition with $\cg'$, where $\cg'$  is induced by the map $(\kappa,\rho)\to (\kappa,R-\rho)$.
\\

The class $\{\tilde\mu_\tau\}$ is represented by the map 
$$\tilde m_\tau:S(E)\times S_\epsilon(H)\times [0,R]\times D_R(V)\map \qmod{[F\times 
W]^+_B}{_B\ [F\times W]^+_B\setminus D_\varepsilon(F\times W)}
$$
given by
$$\tilde m_\tau(e,\kappa,\rho,v)=[\tilde \mu_{\tau}(\rho e,\kappa,v)]\ .
$$
As we have seen in section \ref{cylinder}, this map induces a map
$$S(E)_{+B}\wedge_B \underline{S}(H)_{+B}
\wedge_B \underline{S}^1\wedge_B \underline{V}^+_B\map
\qmod{[F\times 
W]^+_B}{_B\ [F\times W]^+_B\setminus D_\varepsilon(F\times W)}$$
because it has the following properties
\begin{enumerate}
\item $\tilde m_\tau(e,\kappa,0,v)$ and  $m_\tau(e,\kappa,R,v)$ belong always 
to the infinity section of the right 
hand space,
\item $\tilde m_\tau(e,\kappa,\rho,v)$ belongs to the infinity section of the 
right  hand space when $\|v\|=R$.
\end{enumerate}
 The class 
$-\partial(\{\tilde\mu_\tau\})$ is defined by the map
$$\tilde m_\tau':S(E)\times S_\epsilon(H)\times [0,R]\times D_R(V)\to \qmod{[F\times 
W]^+_B}{_B\ [F\times W]^+_B\setminus D_\varepsilon(F\times W)}$$
given by 
$$\tilde m_\tau'(e,\kappa,\rho,v)=\tilde m_\tau\left(e,\kappa,R-\rho,
v\right)\ .
$$
This map descends to a map 
$$S(E)_{+B}\wedge_B \underline{H}^+_B
\wedge_B   \underline{V}^+_B\to
\qmod{[F\times 
W]^+_B}{_B\ [F\times W]^+_B\setminus D_\varepsilon(F\times W)}$$
because it has the following properties:
\begin{enumerate}
\item   $\tilde m_\tau'(e,\kappa,0,v)$ and $\tilde m_\tau'(e,\kappa,R,v)$ are independent 
of $\kappa$.
\item $\tilde m_\tau'(e,\kappa,R,v)$ belongs always to the infinity section of the right 
hand space.
\item $\tilde m_\tau'(e,\kappa,\rho,v)$ belongs to the infinity section of the 
right  hand space when $\|v\|=R$.
\end{enumerate}
These three conditions characterize the maps of pointed spaces over $B$ defined  on $S(E)\times S_\epsilon(H)\times [0,R]\times D_R(V)$ which descend to $S(E)_{+B}\wedge_B \underline{H}^+_B
\wedge_B   \underline{V}^+_B$.
\\ \\
Step 2 (a). We deform the map $\tilde m_\tau'$ in the space of maps satisfying the three properties above, by composing it with a 1-parameter family of contractions in the $\rho$-direction.\\

For $t\in [0,1]$ define the map
$$[\tilde m'_\tau]^t:S(E)\times S_\epsilon(H)\times [0,R]\times D_R(V)\to \qmod{[F\times 
W]^+_B}{_B\ [F\times W]^+_B\setminus D_\varepsilon(F\times W)}
$$
by
$$[\tilde m'_\tau]^t(e,\kappa,\rho,v)=\tilde m_\tau\left(e,\kappa,
(1-t+t\frac{\tau}{R})(R-\rho),v\right)\ .
$$
The family $([\tilde m'_\tau]^t)_{t\in[0,1]}$ defines a homotopy in the space of 
maps satisfying properties (1), (2), (3) above. The main point in checking (1) is the fact that the map $\tilde m_\tau$ is constant with respect to $\kappa$ for $\rho\in [\tau,R]$. Therefore
it holds
$$-\partial(\{\tilde\mu_\tau\})=\{[\tilde m'_\tau]^0\}=\{[\tilde m'_\tau]^1\}\ .
$$
Putting $\tilde m''_\tau:=[\tilde m'_\tau]^1$, one has 
$$\tilde m''_{\tau}(e,\kappa,\rho,v)=\tilde m_\tau (e,\kappa,
\frac{\tau}{R}(R-\rho),v )=
\left[1-\frac{R-\rho}{R}\right](\kappa+l(v))+ \frac{R-\rho}{R}
\mu(\tau\frac{R-\rho}{R}e,v) 
$$
$$=
\frac{\rho}{R}\kappa+l(v)+\frac{R-\rho}{R}
\left(\mu(\tau\frac{R-\rho}{R}e,v)-l(v)\right)\ .
$$
Step 2 (b). We remark that the family of maps $\tilde m''_\tau$ has a uniform  limit as $\tau\to 0$ and we compute this limit explicitly .\\

Using arguments as in the proof of the claim above, we see that 
$$\lim_{\tau\to 0}\frac{R-\rho}{R}
\left(\mu(\tau\frac{R-\rho}{R}e,v)-l(v)\right)=0
$$
uniformly. Therefore $\tilde m'':=\lim_{\tau\to 0} \tilde m''_\tau$ operates by  $\tilde m''(e,\kappa,\rho,v)=\frac{\rho}{R}\kappa+l(v)$. 
It is now easy to see that the map 
$$S(E)_{+B}\wedge_B \underline{H}^+_B
\wedge_B   \underline{V}^+_B\to \qmod{[F\times 
W]^+_B}{_B\ [F\times 
W]^+_B\setminus D_\varepsilon(F\times W)}=F^+_B\wedge_B\underline H^+_B\wedge_B[\underline W_0]^+_B
$$
induced by $\tilde m''$ is homotopic to the smash product over $B$ of the obvious map  $S(E)_{+B}\to F^+_B$ (which represents $o(E,F)$) with  $l^+_B:\underline V^+_B\to [\underline W_0]^+_B$, and 
$\id:\underline H^+_B\to \underline H^+_B$.
\qed
\\

For a map $\kappa:B\to S_\epsilon(H)$ one has
\begin{equation}\label{individual}
\{\mu_\kappa\}=(j^E_\kappa)^*(\{\tilde \mu\})\ .
\end{equation}
This formula shows that the individual invariant $\{\mu_\kappa\}$ associated with a map $\kappa:B\to S_\epsilon(H)$ is determined by the invariant associated with the universal perturbation $\tilde \mu$ and the homotopy class of $\kappa$. Using Corollary \ref{differenceclassnew} we obtain
\begin{co} \label{jump} (Cohomotopy invariant jump formula) One has
$$\{\mu_{\kappa_0}\}-\{\mu_{\kappa_1}\}=o_{(E,F)}\cdot d(\kappa_0,\kappa_1)\ ,
$$
where $d(\kappa_0,\kappa_1)\in {_{S^1}\alpha}^{-1}_B(B_{+B}, \underline{H}^+_B)$ is the difference class of the maps $\kappa_0$, $\kappa_1$ regarded as sections in the sphere bundle $S_\epsilon(\underline{H})$.
\end{co}
 
Suppose now that $b=1$. In this case $S_\epsilon(H)$ has two elements $\kappa_0$, $\kappa_1$, and the difference class $d(\kappa_0,\kappa_1)$ is just the unit element of  ${_{S^1}\alpha}^{0}_B(B_{+B}, B_{+B})$. Therefore, in this case, our result gives 
\begin{co}\label{wc} (Cohomotopy wall crossing) Suppose $b=1$.  Then the two classes   $\{\mu_{\kappa_0}\}$, $\{\mu_{\kappa_1}\}$ associated with the two perturbations $\mu_{\kappa_0}$,  $\mu_{\kappa_1}$ of $\mu$ are related by the formula
$$\{\mu_{\kappa_0}\}-\{\mu_{\kappa_1}\}=\{o(E,F)\}\ .
$$
\end{co} 

We can now extend our results to the infinite dimensional case.  Let  $B$ be an oriented compact manifold, ${\cal E}$, ${\cal F}$ complex Hilbert bundles over $B$, ${\cal V}$, ${\cal W}$ real Hilbert spaces, and $\mu:{\cal E}\times{\cal V}\to {\cal F}\times{\cal W}$  an $S^1$-equivariant,  fiberwise differentiable map over $B$ satisfying properties ${\cal P}1$, ${\cal P}3$ and ${\cal P}2$ (1) with $h= 0$.  Then we have an orthogonal decomposition ${\cal W}=H\oplus {\cal W}_0$, and $\mu(0^{\cal E}_y,v)=l(v)$ for every $v\in {\cal V}$, where $l:{\cal V}\to {\cal W}_0$ is a linear isometry. We fix an orientation of the finite dimensional summand $H$. Defining in the same way as in the finite dimensional framework the universal perturbation $\tilde \mu$, one gets a stable class 
$$\{\tilde \mu\}\in \alpha^*({S}_\epsilon(H);x)\ ,$$
where $x\in K(B)$ is the  index of the complex part of the fiberwise linearization of $\mu$ at the zero section. Recall that the Euler class $\gamma(x)\in \alpha^0(x)$ is defined  by the system of stable classes $-\{o_{(E,F)}\}\in {_{S^1}\alpha}^0_B(S(E)_{+B},F^+_B)$ defined by the obvious maps $S(E)_{+B}\to F^+_B$ (see  section \ref{euler}).  Using the results obtained above  and taking inductive limit over ${\cal T}(x)$, we obtain
\begin{co}\label{finalresult}
  \begin{enumerate}
\item  The image of $\{\tilde \mu\}$ under the morphism  
$$\partial_x:\alpha^{b-1}(S_\epsilon(H);x) \to  \alpha^0(x)$$
 is given by
 $$\partial_x(\{\tilde\mu\})=\gamma(x)\ .$$
\item Let $\kappa_0,\ \kappa_1:B\to S(H)$  two maps. Then
$$\{\mu_{\kappa_1}\}-\{\mu_{\kappa_0}\}=d(\kappa_0,\kappa_1)\cdot\gamma(x)\ .
$$
\item Suppose $b=1$ and write $S_\epsilon(H)=\{\kappa_0,\kappa_1\}$. Then 
$$\{\mu_{\kappa_1}\}-\{\mu_{\kappa_0}\}=\gamma(x)\ .
$$
\end{enumerate}

\end{co}

\subsection{A product formula and a vanishing theorem}

In this section we give the infinite dimensional analogue of the product formula proven in section \ref{prodsubsection}.

Let  ${\cal V}_i$, ${\cal W}_i$ be real Hilbert spaces , ${\cal E}_i$, ${\cal F}_i$  complex Hilbert bundles over a compact base $B$ ($i=1,\ 2$), and let $\mu_i:E_i\times  V_i\to [F_i\times W_i]^+_{B}$  be $S^1$-equivariant maps over $B$, satisfying the properties ${\cal P}1$, ${\cal P}2$ (1) and ${\cal P}3$ of section \ref{general} with constants $C$, $c$. Let ${\cal W}_i=H_i\oplus {\cal W}_{0,i}$ be the corresponding orthogonal sum decompositions,   $l_i:{\cal V}_i\stackrel{\simeq}{\to} {\cal W}_{0,i}$ isometries, $x_i\in K(B)$ the K-theory elements defined by the corresponding families $\delta_i$ of Fredholm operators, and $h_i:B\to H_i$ the maps given by ${\cal P}2$ (1). We introduce the notations: 
$${\cal V}:={\cal V}_1\oplus {\cal V}_2,\ {\cal W}:={\cal W}_1\oplus {\cal W}_2,\ H:=H_1\oplus H_2,\ {\cal W}_0:={\cal W}_{0,1}\oplus {\cal W}_{0,2}\ ,\ l:=l_1\oplus l_2\ ,$$
 and  consider the Hilbert bundles ${\cal E}:={\cal E}_1\oplus  {\cal E}_2$, ${\cal F}:={\cal F}_1\oplus {\cal F}_2$. The   product map 
$$\mu:{\cal E}\times {\cal V}=[{\cal E}_1\times {\cal V}_1]\oplus [{\cal E}_2\times {\cal V}_2]  \map  [{\cal F}\times {\cal W}]^+_B=[{\cal F}_1\times {\cal W}_1]^+_B\wedge_B [{\cal F}_2\times {\cal W}_2]^+_B$$
also   satisfies  properties ${\cal P}1$ ${\cal P}2$ (1)  (with associated map $h=(h_1,h_2):B\to H$)  and  ${\cal P}3$; it satisfies ${\cal P}2$ (2) as soon as one of the two maps $\mu_1$, $\mu_2$ does.

Suppose that  $\mu_1$   satisfies property  ${\cal P}2$  (2). In this case the construction of section \ref{general} applies and yields an invariant   
 $$\{\mu_1\}\in  \alpha^{b_1-1}(x_1)\ .$$
The finite dimensional approximations  of the map $\mu_2$ define classes 
$$\{(\mu_2)_{c,\pi_2}^+\}\in  {_{S^1}\alpha}_{B}^{b_2}([E_2]^+_{B},[F_2]^+_{B})\ .$$
It can be shown that a compatibility result similar to Proposition \ref{compat}  holds, so that one obtains an invariant
$$\{\mu_2^+\}\in \alpha^{b_2}(x_2^+):=\varinjlim\hskip-28pt\raisebox{-10pt}{$\scriptstyle (E_2,F_2)\in x_2$\ \ }{_{S^1}\alpha}_B^{b_2}([E_2]^+_{B},[F_2]^+_{B})\ .
$$
Here the inductive limit on the right is taken over the category ${\cal T}(x_2)$ and is constructed using the same methods as in the definition of the groups $\alpha^*(x)$ (see section \ref{alpha}). 
The direct limit of the obvious products 
$${_{S^1}\alpha}_{B}^{b_1-1}(S(E_1)_{+B},[F_1]^+_B)\times {_{S^1}\alpha}_{B}^{b_2}([E_2]^+_{B},[F_2]^+_{B})\to $$
$$ {_{S^1}\alpha}_{B}^{b_1+b_2-1}(S(E_1)_{+B}\wedge_B [E_2]^+_{B_2},[F_1\oplus F_2]^+_B)\stackrel{{_1c}^*}{\to } {_{S^1}\alpha}_{B}^{b_1+b_2-1}(S(E_1\oplus E_2)_{+B},[F_1\oplus F_2]^+_B)
$$
gives a well defined product
$$\cdot: {_{S^1}\alpha}^{b_1-1}(x_1)\times \alpha^{b_2}(x_2^+)\to {_{S^1}\alpha}^{b_1+b_2-1}(x_1+x_2)\ .
$$
Using finite dimensional approximations of $\mu$ of the form 
$$\mu_{c,\pi_1\times\pi_2}=(\mu_1)_{c,\pi_1}\times(\mu_2)_{c,\pi_2}$$
 and applying  Proposition \ref{productformula} we obtain
\begin{re}   Under the assumptions and with the notations above, the invariant of the product map $\mu=\mu_1\times\mu_2$ is given by the formula
$$\{\mu_1\times\mu_2\}=\{\mu_1\}\cdot \{\mu_2^+\}\ .
$$
\end{re}

Note that in this formula the map $\mu_2$ is allowed  to have $S^1$-invariant zeroes.
In the case when both maps $\mu_i$ satisfy  ${\cal P}2$ (2) (so they are nowhere zero on their $S^1$-fixed point loci) one has the following important vanishing result for the Hurewicz image of the invariant associated with a product map:
\begin{pr} Put $x:=x_1+x_2\in K(B)$ and let $h_x: \alpha^*(x)\to H^*(x;\Z)$ be the Hurewicz morphism associated with $x$. Suppose that both maps $\mu_i$ satisfy properties ${\cal P}1$, ${\cal P}2$ (1), ${\cal P}2$ (2) and ${\cal P}3$, and that $B$ is a finite CW complex. Then 
$$h_x(\{\mu_1\times\mu_2\})=0\ .
$$
\end{pr}
\pf Let $m_i:=(\mu_i)_{c,\pi_i}$ be   finite dimensional approximations of $\mu_i$ and put $m:=m_1\times m_2$. Applying the cylinder construction to this maps we get a representative
$$
m_R:S(E_1\oplus E_2)_{+B}\wedge_B[\underline{\R}\oplus \underline{V_1}\oplus\underline{V_2}]^+_B\to [F_1\oplus F_2\oplus\underline{W_1}\oplus\underline{W_2}]^+_B
$$
of the class $\{\mu_1\times\mu_2\}$. Put $E:=E_1\oplus E_2$, $F:=F_1\oplus F_2$, $V:=V_1\oplus V_2$, $W:=W_1\oplus W_2$, and $b=b_1+b_2$. Let
$$\bar m_R:[\underline{\R}\oplus \underline{V}]^+_{\P(E)}\to [\tilde F\oplus \underline{W}]^+_{\P(E)}
$$
be the associated sphere bundle map, constructed as in section \ref{hurewiczm}. We denote by 
$$\p:[\underline{\R}\oplus \underline{V}]^+_{\P(E)}\to\P(E)\ ,\ {\mathrm q}:[\tilde F\oplus \underline{W}]^+_{\P(E)}\to\P(E) $$
 the two bundle projections, 
and  by $h:=h_{\bar m_R}\in H^{2f+b_1+b_2-1}(\P(E);\Z)$ the corresponding Hurewicz class, which is defined by the equality
\begin{equation}\label{thhu}
(\bar m_R)^*(\t_{\tilde F\oplus\underline{W}})=\p^*(h)\cup \t_{\underline{\R}\oplus\underline{V}}
\end{equation}
in $H^*([\underline{\R}\oplus \underline{V}]^+_{\P(E)},\infty_{\underline{\R}\oplus \underline{V}};\Z)$.  Since both maps $\mu_i$ satisfy property {\bf P2}, it follows that,  for a sufficiently small neighborhood ${\cal P}$ of $\P(E_1)\cup\P(E_2)$ in $\P(E)$, the map   $\bar m_R$ maps $\p^{-1}({\cal P})$  to the infinity section of the right hand bundle.  We can suppose that ${\cal P}$ is a standard compact neighborhood of this union, i.e. it has the form
$${\cal P}=\P(E)\setminus \left\{[e_1,e_2]\in \P(E)|\ e_i\ne 0,\ \ln \frac{\|e_1\|}{\|e_2\|} \in (-s\,s)\right\}
$$
for sufficiently large $s>0$. The pull-back class $(\bar m_R)^*(\t_{\tilde F\oplus\underline{W}})$  can be regarded as an element in $H^*([\underline{\R}\oplus \underline{V}]^+_{\P(E)},\infty_{\underline{\R}\oplus \underline{V}}\cup \p^{-1}({\cal P});\Z)$, which can be identified  with $H^{*-(\dim(V)+1)}(\P(E),{\cal P};\Z)$ via the relative Thom isomorphism over the pair $(\P(E),{\cal P})$. Therefore, the equality 
\begin{equation}\label{thhunew}
(\bar m_R)^*(\t_{\tilde F\oplus\underline{W}})=\p^*(h')\cup \t_{\underline{\R}\oplus\underline{V}}
\end{equation}
in $H^*(([\underline{\R}\oplus \underline{V}]^+_{\P(E)},\infty_{\underline{\R}\oplus \underline{V}}\cup \p^{-1}({\cal P});\Z)$ defines a class $h'\in H^*(\P(E),{\cal P};\Z)$, and $h$ is just the image of $h'$ via the  morphism $C^*:H^*(\P(E),{\cal P};\Z)\to H^*(\P(E);\Z)$ associated with the map   $C:(\P(E),\emptyset)\to (\P(E),{\cal P})$. Put now 
$$\P_0:=\P(E)\setminus(\P(E_1)\cup\P(E_2))\ ,\ {\cal P}_0:={\cal P}\setminus(\P(E_1)\cup\P(E_2))\ ,$$
 and denote by $h'_0$ the  image  of $h'$ via the morphism $I^*:H^*(\P(E),{\cal P};\Z)\to  H^*(\P_0,{\cal P}_0;\Z)$ defined by the map   $I:(\P_0,{\cal P}_0)\to (\P(E),{\cal P})$.  The main point in the proof of our proposition is that  the  restriction  
$$
\resto{\bar m_R}{\P_0}: \p^{-1}(\P_0)\to {\mathrm q}^{-1}(\P_0) \ .
$$
 is equivariant with respect to the {\it free} $S^1$-action $(\zeta,[e_1,e_2])\mapsto [\zeta e_1,e_2]$ on $\P_0$ and the obvious lift of this action in the  bundle $\resto{\tilde F}{\P_0}$. This is just because  $\mu$ is the product of two $S^1$-equivariant maps $\mu_i$. Therefore, $\resto{\bar m_R}{\P_0}$ descends to a bundle map
$$[\bar n_R]_0:\qmod{\p^{-1}(\P_0)}{S^1}\map \qmod{{\mathrm q}^{-1}(\P_0)}{S^1} 
$$
over $\Q_0:=\P_0/S^1$. The two sphere bundles above coincide with the fibrewise compactifications  $[\underline{\R}\oplus \underline{V}]^+_{\Q_0}$, $[\tilde F_0\oplus \underline{W}]^+_{\Q_0}$, where $\tilde F_0$ is the $S^1$-quotient of $\tilde F$, regarded as a bundle over $\Q_0$. We denote by $\p_0$, ${\mathrm q}_0$ the corresponding bundle  projections on $\Q_0$. Put ${\cal Q}:= {\cal P}/{S_1}$, ${\cal Q}_0:={\cal Q}\cap \Q_0$. Using   the relative Thom isomorphism over the pair $(\Q_0,{\cal Q}_0)$, it follows that the equality 
$$[\bar n_R]_0^*(\t_{\tilde F_0\oplus \underline{W}})={\mathrm p}_0^*(k_0)\cup\t_{\underline{\R}\oplus \underline{V}}
$$
defines a class $k_0\in H^*(\Q_0,{\cal Q}_0;\Z)$. Taking the pull-back of this equality via the projection $\Pi_0:(\P_0,{\cal P}_0)\to (\Q_0,{\cal Q}_0)$,
 (and comparing the obtained formula with a similar equality satisfied by $h'_0$), we see that $\Pi^*_0(k_0)=h'_0$. Therefore
\begin{equation}\label{fromk}
h=C^*\circ {I^*}^{-1}\circ \Pi_0^*(k_0)=C^*\circ \Pi^*\circ [J^*]^{-1}(k_0)\ ,
\end{equation}
where 
$$\Pi: (\P(E),{\cal P})\to \left(\qmod{\P(E)}{S^{1}},{\cal Q}\right)\ ,\ J: (\Q_0,{\cal Q}_0)\to \left(\qmod{\P(E)}{S^{1}},{\cal Q}\right)$$
 denote the obvious maps. In this formula we used the identity $J\circ\Pi_0=\Pi\circ I$,  and that the maps $I$, $J$   induce  isomorphisms in cohomology, by the excision theorem.    The result follows now directly from Lemma \ref{spectral} below.
\qed

\begin{lm}\label{spectral}  The morphism 
$$U^*:H^*\left(\qmod{\P(E)}{S^{1}},{\cal Q};\Z\right)\map  H^*(\P(E);\Z)$$
induced by the map $U:=\Pi\circ C:(\P(E),\emptyset)\to \left(\qmod{\P(E)}{S^{1}},{\cal Q}\right)$, vanishes.
\end{lm}
\pf  By the excision and homotopy invariance theorem   one  has 
$$H^*\left(\qmod{\P(E)}{S^{1}},{\cal Q};\Z\right)=H^*\left(\qmod{\P(E)}{S^{1}}\setminus \cringle{{\cal Q}},{\cal Q}\setminus \cringle{{\cal Q}};\Z\right)\ ,$$
where $\cringle{{\cal Q}}$ is the interior of ${\cal Q}$.  One has a natural homeomorphism 
$$\qmod{\P(E)}{S^{1}}\setminus \cringle{{\cal Q}}\cong [\P(E_1)\times_B\P(E_2)]\times [-s,s]\ ,\ [e_1,e_2]\mapsto\left([e_1],[e_2],\ln \frac{\|e_1\|}{\|e_2\|} \right)\ ,$$
and   this homeomorphism identifies   ${\cal Q}\setminus \cringle{{\cal Q}}$   with $[\P(E_1)\times\P(E_2)]\times\{-s,s\}$. Multiplication with the  Thom class of the trivial bundle   
$$\P(E_1)\times_B\P(E_2)\times(-s,s)\to \P(E_1)\times_B\P(E_2)$$
 defines an isomorphism 
$$H^i(\P(E_1)\times_B\P(E_2);\Z)\stackrel{\cong}{\to} H^{i+1}\left(\qmod{\P(E)}{S^{1}}\setminus \cringle{{\cal Q}},{\cal Q}\setminus \cringle{{\cal Q}};\Z\right)=H^{i+1}\left(\qmod{\P(E)}{S^{1}},{\cal Q};\Z\right).
$$
{\it Step 1.} When $B$ is a point, the statement  of the Lemma is obvious because in this case both  spaces $\P(E_1)\times_B\P(E_2)$ and   $\P(E)$ have trivial cohomology in odd dimensions. 
\\ \\
{\it Step 2.} For a general basis, note that $U$ induces a  morphism of the Leray spectral sequences associated with the projections 
$$\P(E)   \map B\  ,\   \left(\qmod{\P(E)}{S^{1}},{\cal Q}\right)\map B\ .
$$
But the  Leray  spectral sequence for the relative cohomology of the pair  $\left(\qmod{\P(E)}{S^{1}},{\cal Q}\right)$  can be identified with the  spectral sequence for the cohomology with compact supports of $\qmod{\P(E)}{S^{1}}\setminus {\cal Q}$.  It suffices to note that the induced spectral sequence morphism   vanishes at the $E_1^{p,q}$-level, by {\it Step 1}.
\qed

\section{Appendix}

\subsection{Inductive limits of functors}\label{limits}
We recall the following important  
\begin{dt}  (\cite{AM} p. 148) \label{filtering}
A filtering category is category ${\cal C}$ with the properties
\begin{enumerate}
\item[F1.]  For every pair $(O,O')$ of objects, there exists an object $O''$ and morphisms $O\to O''$, $O'\to O''$.
\item[F2.]  For every two morphisms  $u$, $v:O\to O'$ there exists an object $O''$ and a morphism $w:O'\to O''$ such that $w\circ u=w\circ v$.
\end{enumerate}
\end{dt}
For {\it small} filtering categories one has the following basic fact:
\begin{pr}\label{old} (\cite{AM}, p. 149-150) Let ${\cal A}$ be one of the categories ${\cal S}ets$, ${\cal A}b$ or ${\cal G}r$, and let ${\cal C}$ be a filtering  {small} category. Then any functor  $F:{\cal C}\to {\cal A}$ has an inductive limit, which can constructed in the classical way: one  factorizes the disjoint union  
$\coprod_{O\in {\cal O}b({\cal C})} F(O)$ by the equivalence relation 
\begin{equation}\label{relation}
(O,x)\sim (O',x')\hbox{ if } \exists\ u:O\to O'',\ u':O'\to O''\hbox { with } F(u)(x)=F(u')(x')\ .
\end{equation}
When ${\cal A}={\cal A}b$ or ${\cal G}r$, one endows the obtained set of equivalence classes with the operation induced by the group operations on the summands $F(O)$ of the disjoint union.
\end{pr}
We will say that ${\cal C}$ is {\it weakly filtering} if it satisfies F1 and   the following weak form of the axiom $F2$. 
\begin{enumerate}
\item[$\tilde{\rm F}2$.]  For every two morphisms  $u$,  $v:O\to O'$ there exists an object $O''$ and  morphisms $w$,$z:O'\to O''$ such that $w\circ u=z\circ v$.
\end{enumerate}
\begin{lm} Suppose that ${\cal C}$ is weakly filtering and small. Then the relation $\sim$ defined in (\ref{relation}) is still an equivalence relation, and the conclusion of Proposition \ref{old} holds for ${\cal A}={\cal S}ets$.
\end{lm}
\pf  It suffices to check that $\sim$ is transitive. Let $x\in F(O)$, $x'\in F(O')$, $x''\in F(O'')$ with $x\sim x'$, $x'\sim x''$. Therefore there exists morphisms $u:O\to \hat O$, $u':O'\to \hat O$, $v':O'\to \tilde O$, $v'':O''\to \tilde O$ such that $F(u)(x)=F(u')(x')$ and $F(v')(x')=F(v'')(x'')$.  By F1 there exists morphisms $\hat w:\hat O\to O_0$, $\tilde w:\tilde O\to O_0$. We apply $\tilde{\rm F}2$ to the morphisms $\hat w u'$, $\tilde w v':O'\to O_0$. We obtain morphisms $\hat z$, $\tilde z: O_0\to O_1$ such that $\hat z\hat w u'=\tilde z\tilde w v'$. Therefore
$$F( \hat z\hat w u)(x)=F(\hat z\hat w)(F(u)(x))=F(\hat z\hat w)(F(u')(x'))=F(\hat z\hat w u')(x')=$$
$$=F(\tilde z\tilde w v')(x')=F(\tilde z\tilde w)(F(v')(x'))=F(\tilde z\tilde w)(F(v'')(x''))=F(\tilde z\tilde w v'')(x'')\ ,
$$
hence $x\sim x''$.
\qed
\vspace{1mm}\\
For ${\cal A}={\cal A}b$ or ${\cal G}$  one cannot endow the quotient of the disjoint union by this equivalence relation with a coherent group structure using only the weakly filtering condition.

Unfortunately, we will need inductive limits of functors defined on index categories which are not small.  In this case the disjoint union considered in  Remark \ref{old} might not be a set.  However, there exists a simple situation when the existence of an inductive limit is guaranteed:

\begin{lm}\label{epi}  Let ${\cal C}$ be a weakly filtering category, $Q\in{\cal O}b({\cal C})$ a fixed object and $F:{\cal C}\to{\cal A}$ a functor  such that  $F(u)$ is  surjective for every morphism  $u:Q\to O$.%
\begin{enumerate}
\item Suppose ${\cal A}={\cal S}ets$. 
\begin{enumerate}
\item The relation on $F(Q)$ defined by   
 \begin{equation}\label{newrelation}
 y\approx y' \  \hbox{ if }\  \exists\ u,\ v:Q\to O \hbox{ such that } F(u)(y)=F(v)(y')  
 \end{equation}
is an equivalence relation. Put $L:=F(Q)/\approx$.
\item 
For any $O\in {\cal O}b({\cal C})$  there exists a unique map $f_O:F(O)\to L$ defined by $f_O(x)=[y]$ for any pair $(x,y)\in F(O)\times F(Q)$ for which  there exist morphisms $u:O\to \hat O$, $v:Q\to \hat O$ with $F(u)(x)= F(v)(y)$. The system $(f_O)_{O\in {\cal O}b({\cal C})}$ is $F$-compatible (i.e. it holds $f_{O'}\circ F(w)=f_O$ for any morphism $w:O\to O'$).
\item The system $(f_O)_{O\in {\cal O}b({\cal C})}$ satisfies the universal property of the inductive limit, so the inductive limit of $F$ exists and can be identified with $L$.
\end{enumerate} 
\item Suppose ${\cal A}={\cal A}b$ or ${\cal G}r$. 
\begin{enumerate}
\item Let $H$  be a smallest   normal subgroup of $F(Q)$ which contains the elements $x'x^{-1}$ with $x\approx x'$. Put $L:=F(Q)/H$.
\item The  system of morphism  $(f_O:F(O)\to L)_{O\in{\cal O}b({\cal C})}$ defined  in a similar way as in (1) is $F$-compatible and satisfies the universal property of the inductive limit. Therefore the inductive limit of $F$ exists and can be identified with $L$.
\end{enumerate}
\end{enumerate} 
\end{lm} 
\pf (1) (a) is clear. For (b) we have to prove that the map $f_O$ is well defined.  Let $y\in F(Q)$, $y'\in F(Q)$,  $u:O\to\hat O$, $v:Q\to \hat O$,  $u':O\to\hat O'$, and $v':Q\to \hat O'$ such that $F(u)(x)=F(v)(y)$ and $F(u')(x)=F(v')(y)$. We can find an object $\tilde O$ and morphisms $w:\hat O\to\tilde O$, $w':\hat O'\to\tilde O$. Since ${\cal C}$ is weakly filtering, the exist morphisms $z:\tilde O\to O_0$,  $z':\tilde O\to O_0$ such that $zwu=z'w'u'$. This implies 
$$F(zwv)(y)=F(zw)(F(u)(x))=F(z'w')(F(u')(x))=F(zwv')(y')\ ,$$
so $y\approx y'$. 

The $F$-compatibility of the system  $(f_O)_{O\in {\cal O}b({\cal C})}$ and the fact that this system satisfies the universal property of the inductive limit are easily verified.  
\vspace{2mm}\\
(2) Follows easily from (1).
\qed

\begin{dt} (\cite{AM} p. 149) Let ${\cal N}$, ${\cal C}$ be  categories. A functor $\Theta:{\cal N}\to {\cal C}$ is  called
\begin{enumerate}
\item  cofinal, if
\begin{enumerate}
\item[C1.] For any   $O\in{\cal O}b({\cal C})$ there  exists $n\in {\cal O}b({\cal N})$ and  $u:O\to\Theta(n)$.
\item[C2.]  For every $n\in {\cal O}b({\cal N})$, $O\in{\cal O}b({\cal C})$, and  $u:\Theta(n)\to O$, there exists $m\in {\cal O}b({\cal N})$, $\nu:n\to m$  and $v:O\to\Theta(m)$ such that $vu=\Theta(\nu)$.
\end{enumerate}
\item cofinal in the sense of Artin-Mazur (\cite{AM} p. 149) , if 
\begin{enumerate} 
\item[C1.] holds,
\item[$\tilde{\rm C}$2.] For every $O\in{\cal O}b({\cal C})$, $n\in {\cal O}b({\cal N})$  and $u$,  $v:O\to\Theta(n)$, there exists  a morphism $\mu:n\to m$ in ${\cal N}$ such that $\Theta(\mu)u=\Theta(\mu) v$.
\end{enumerate}
\end{enumerate}
\end{dt}
\begin{lm}\label{AMOT}\begin{enumerate}
\item If ${\cal N}$ is filtering and $\Theta$ is cofinal in the sense of Artin-Mazur, then $\Theta$ is cofinal and ${\cal C}$ is filtering. 
\item If ${\cal C}$ is filtering  and $\Theta$ is cofinal, then $\Theta$ is cofinal in the sense of Artin-Mazur.
\item Suppose $\Theta:{\cal N}\to{\cal C}$ is cofinal, and ${\cal N}$, ${\cal C}$ are both small and  filtering. For any functor $F:{\cal C}\to {\cal A}$   (with ${\cal A}={\cal S}ets$, ${\cal A}b$ or ${\cal G}r$)  the canonical morphism
$$
\varinjlim\hskip-24pt\raisebox{-11pt}{$\scriptstyle n\in {\cal O}b({\cal N})$}F(\Theta(n))\to \varinjlim\hskip-24pt\raisebox{-11pt}{$\scriptstyle O\in {\cal O}b({\cal C})$}F(O)\
$$ 
is an isomorphism.
\end{enumerate}
\end{lm}
\pf  1. Let $u:\Theta(n)\to O$ be a morphism. Using C1, we can find a morphism $w:O\to \Theta(m)$; since ${\cal N}$ is filtering, we can find morphisms $\eta:n\to k$, $\kappa:m\to k$. Therefore, we get two morphisms $\Theta(\eta)$, $\Theta(\kappa)wu:\Theta(n)\to\Theta(k)$. By $\tilde{\rm C}$2, there exists $\mu:k\to l$ such that $\Theta(\mu)\Theta(\eta)=\Theta(\mu)\Theta(\kappa)wu$. This shows $[\Theta(\mu\kappa)w]u=\Theta(\mu\eta)$, so C2 holds with $v=\Theta(\mu\kappa)w$ and $\nu=\mu\eta$.  The fact that ${\cal C}$ is filtering is stated in \cite{AM} p. 149.
\vspace{2mm}\\
2.  Let $u$,  $v:O\to\Theta(n)$ be two morphisms. Since  ${\cal C}$ is filtering, there exists $w:\Theta(n)\to O'$ with $w u=wv$. By C2, we can find $m\in{\cal O}b({\cal N})$, $\nu:n\to m$ and $v':O'\to \Theta(n)$, such that $v' w=\Theta(\nu)$. We will have $\Theta(\nu)u=v'wu=v'wv=\Theta(\nu)v$, which proves $\tilde{\rm C}$2.
\vspace{2mm}\\
3. See Proposition 1.8 in \cite{AM} p. 150.
\qed
\hskip-4mm{\bf Example 1.} Let $B$ be a compact space and let ${\cal U}_B$ be the category of complex vector bundles over $B$. A morphism $U\to U'$ is a pair $u=(i,U_1)$ consisting of a bundle embedding $i:U\to U'$ and a complement $U_1$ of $i(U)$ in $U'$ (see section \ref{alpha}).  The category ${\cal U}_B$ satisfies F1 but not F2, so it is {\it not} filtering. Let ${\cal N}$ be category associated with the ordered set $(\N,\leq)$. Then the functor $\Theta:{\cal N}\to{\cal U}_B$ which associates to $n$ the trivial  bundle $\underline{\C}^n$  and to an inequality $n\leq m$ the standard morphism $\underline{\C}^n\to \underline{\C}^m$ is cofinal. This follows from the fact that  any vector bundle on $B$ possesses a complement.   Note however that $\Theta$ is not cofinal in the sense of Artin-Mazur.\vspace{2mm}\\
{\bf Example 2.} For a category ${\cal C}$ and an object $Q\in {\cal O}b({\cal C})$ we will denote by ${\cal C}_{Q}$ the category whose objects are morphism $u:Q\to O$ and whose morphisms are
$$\Hom(Q\stackrel{u}{\to}  O,Q\stackrel{v}{\to} O'):=\{w:O\to O'|\ w\circ u=v\}\ .
$$
A morphism $u:Q\to Q'$ induces in an obvious way a pull-back functor $u^*:{\cal C}_{Q'}\to{\cal C}_Q$.  If ${\cal C}$ is filtering then ${\cal C}_Q$ is filtering and the target functor functor $T:{\cal C}_{Q}\to {\cal C}$ is both cofinal and cofinal in the sense of Artin-Mazur.
 \begin{dt} A   category with automorphism push-forward is a pair $({\cal U},A)$, where ${\cal U}$ is a category and $A:{\cal U}\to{\cal G}r$ a functor, such that  
\begin{enumerate}
\item [F1.] holds in ${\cal U}$.
\item[S1.] $A(O)=\Aut(O)$ for every $O\in{\cal O}b({\cal C})$.
\item[S2.] For any   $u:O\to O'$  and $a\in\Aut(O)$ one has $A(u)(a)\circ u= u\circ a$
\item[S3.] For every two morphisms  $u$, $v:O\to O'$ in ${\cal U}$ there exists an object  $O''$,  a morphism $w:O'\to O''$ and an automorphism $a\in\ A(O'')$ such that $a\circ w\circ u=w\circ v$.
\end{enumerate}
\end{dt}
Note that when $({\cal U},A)$ is a  category with automorphism push-forward, then ${\cal U}$ is weakly filtering (use S3).\\

\hskip -11pt{\bf Example 3.} Defining the automorphism push-forward functors in the obvious way, the categories ${\cal U}_B$, ${\cal C}_B$, ${\cal T}(x)$ introduced in this article become  categories with automorphism push-forward.  
\vspace{2mm}\\
Let $({\cal U},A)$ be a  category with automorphism push-forward, $Q\in{\cal O}b({\cal U})$ a fixed object, and $F:{\cal U}\to{\cal A}b$ a functor such  that $F(u)$ is a isomorphism for any morphism $u:Q\to O$. We know by Lemma \ref{epi} that   the inductive limit of $F$ exists and is a quotient of $F(Q)$. We need an explicit description of this quotient.   For every object $u:Q\to O$ in the category ${\cal U}_Q$   the group $A(T(u))$ acts on $F(Q)$ via the isomorphism $F(u):F(Q)\to F(T(u))$. A morphism $w:T(u)\to T(v)$ can be regarded as an element in $\Hom_{{\cal U}_Q}(u,v)$ and defines  a group morphism $A(w): A(T(u))\to A(T(v))$ which intertwines the actions of these groups on $G(Q)$.
\begin{pr}\label{iso} Let $({\cal U},A)$ be a category with automorphism push-forward, $Q\in{\cal O}b({\cal U})$ a fixed object, and $F:{\cal U}\to{\cal A}b$ a functor such  that $F(u)$ is a isomorphism for any  $u:Q\to O$. Let ${\cal N}$ be a \ub{small} filtering category and $\Theta:{\cal N}\to{\cal U}_Q$ a   functor satisfying the cofinality axiom C1.
Put 
$$\A:= \varinjlim\hskip-22pt\raisebox{-11pt}{$\scriptstyle n\in {\cal O}b({\cal N})$} A(T(\Theta(n)))\ .$$
Then  $\A$ acts on $F(Q)$ in a natural way, the inductive  limit $\varinjlim\hskip-22pt\raisebox{-11pt}{$\scriptstyle O\in {\cal O}b({\cal U})$}F(O)$ exists and can be identified with the quotient   ${F(Q)}/{I[\A]F(Q)}$.
\end{pr}
\pf  
By Lemma \ref{epi} the inductive limit of $F$ exists and can be identified with a quotient $F(Q)/H$. Here $H$ is the group generated by the elements of the form $x-x'$ where $x$, $x'\in F(Q)$ are such that there exists $u$, $u':Q\to O$ with $F(u)(x)=F(u')(x')$.  We claim that the set of such pairs $(x,x')$ coincides with the set of pairs of the form $(\ag x',x')$ with $x'\in F(Q)$, $\ag\in \A$.

Indeed, if $F(u)(x)=F(u')(x')$, choose $v:O\to\hat O$ and $a\in A(\hat O)$ such that $vu'=avu$. The morphism $vu$ can be regarded as an object in the category ${\cal U}_Q$. Since $\Theta$ satisfies the axiom C1, there exists $n\in{\cal O}b({\cal N})$ and a morphism $vu\to\Theta(n)$ in ${\cal U}_Q$, i.e. a morphism $w:\hat O\to T(\Theta(n))$ such that $wvu=\Theta(n)$. We obtain
$$F(\Theta(n))(x)=F(wvu)(x)=F(wvu')(x')=F(wavu)(x')=F(A(w)(a)wvu)(x')=$$
$$=A(w)(a)(F(wvu)(x'))=A(w)(a)(F(\Theta(n))(x'))\ ,
$$
which shows that $x=\ag x'$, where $\ag$ is the class of $A(w)(a)\in A(T(\Theta(n))$ in $\A$.
Conversely let  $\ag=[a]\in\Ag$ be represented by $a\in A(T(\Theta(n))$ and suppose that $x=\ag x'$. This means $F(\Theta(n))(x)=a(F(\Theta(n))(x'))$ so,  putting $u:=\Theta(n)$, $u':=a\Theta(n)$ one has $F(u)(x)=F(u')(x')$.
\qed
\vspace{2mm}

 Let $({\cal U},A)$ be a category with automorphism push-forward, and  let $G:{\cal C}\to  {\cal A}$ be a functor, where ${\cal A}$ is one of the categories ${\cal S}ets$, ${\cal G}r$ or ${\cal A}b$. 
 \begin{dt}
 We say that the stabilized  automorphisms  act trivially on $G$    if
 \begin{enumerate}
\item[TSA.]
For every $O\in{\cal O}b({\cal C})$, $x\in G(O)$ and $a\in A(O)$ there exists a morphism $u:O\to O'$ such that $G(u)(G(a)(x))=G(u)(x)$.
\end{enumerate}
In the presence of functor $\Theta:{\cal N}\to{\cal U}$, we say that   the $\Theta$-stabilized  automorphisms  act trivially on $G$   if
 \begin{enumerate}
\item[$\Theta$SA.]
For every $n\in{\cal O}b({\cal N})$, $x\in G(\Theta(n))$ and $a\in A(\Theta(n))$ there exists a morphism $\nu:n\to m$ such that $G(\Theta(\nu))(G(a)(x))=G(\Theta(\nu))(x)$.
\end{enumerate}
\end{dt}
\begin{re}\label{obvious} If $\Theta$ is cofinal and $G$ satisfies $\Theta$SA, then it also satisfies TSA.  If ${\cal C}$ is filtering, then any functor $G:{\cal C}\to  {\cal A}$ satisfies  TSA. If, moreover, $\Theta$ is cofinal, then $G$ also satisfies $\Theta$SA.
\end{re}
Let   $({\cal U},A)$ be a  category with automorphism push-forward, and  let $G:{\cal U}\to  {\cal A}$ be a functor.  Let ${\cal N}$ be a {\it small   filtering} category and $\Theta:{\cal N}\to {\cal U}$  a cofinal functor such that $\Theta$SA holds.  Consider the classical inductive limit $L_\Theta:=\varinjlim\hskip-24pt\raisebox{-11pt}{$\scriptstyle n\in {\cal O}b({\cal N})$}G(\Theta(n))$. For every $O\in{\cal O}b({\cal U})$ we define a morphism $f_O:G(O)\to L_\Theta$ by $f_O(x):=[G(v)(x)]$ where $v:O\to \Theta(n)$ is a morphism (whose existence is guaranteed by C1). 
\begin{pr} \label{trivact} Under the assumptions and with the notations above it holds
\begin{enumerate}
\item For any $O\in{\cal O}b({\cal U})$ the map $f_O$ is well defined. The system of maps $(f_O)_{O\in{\cal O}b({\cal U})}$ is $G$-compatible i.e. for any  $u:O\to O'$ one has $f_{O'}\circ G(u)=f_O$. When  ${\cal A}={\cal A}b$ or ${\cal G}r$, the map  $f_O$ is a group morphism.
\item The system $(f_O)_{O\in{\cal O}b({\cal U})}$ satisfies the universal property of the inductive limit, therefore the functor $G$ admits an inductive limit in ${\cal A}$ which can be identified with $L_\Theta$.
\end{enumerate}
\end{pr}
We agree to write $u(x)$, $v(x)$ \dots, instead of $G(u)(x)$, $G(v)(x)$ \dots, to save on notations.

\pf 1. Let $v:O\to\Theta(n)$, $v':O\to \Theta(n')$ be two morphisms. Since ${\cal N}$ is filtering, there exists morphisms $\nu:n\to m$, $\nu':n'\to m$. Applying  axiom  S3 to the morphisms $\Theta(\nu)v$, $\Theta(\nu')v'$,   we get a morphism $w:\Theta(m)\to \hat O$ and an automorphism $a\in A(\hat O)$ such that $w\Theta(\nu')v'=aw\Theta(\nu)v$. Now we apply the axiom C2 to $w$ and we get  morphisms $u:\hat O\to \Theta(k)$, $\mu:m\to k$ such that $uw=\Theta(\mu)$. We have  
$$\Theta(\mu\nu')v'=uw\Theta(\nu')v'=uaw\Theta(\nu)v=A(u)(a)uw\Theta(\nu)v=A(u)(a) \Theta(\mu\nu)v\ .
$$
Using the axiom $\Theta$SA we obtain a morphism $\eta:k\to l$ such that 
$$\Theta(\eta)\left[A(u)(a) \Theta(\mu\nu)v(x)\right]=\Theta(\eta)\left[\Theta(\mu\nu)v(x)\right]\ .$$
Therefore $\Theta(\eta\mu\nu')(v'(x))=\Theta(\eta\mu\nu)(v(x))$, which shows that $v(x)$ and $v'(x')$ define the same element in $L_\Theta$. The second and the third claim are   obvious. 
\vspace{2mm}\\
2.  Let $\Lambda\in{\cal O}b({\cal A})$ and $(g_O)_{O\in{\cal O}b({\cal U})}$, $g_O:G(O)\to \Lambda$ a system of $G$-compatible morphisms.  Using the system $(g_{\Theta(n)})_{n\in{\cal O}b({\cal N})}$ (which is $G\circ\Theta$-compatible) we get a unique morphism $g:L_\Theta\to\Lambda$ such that $g\circ c_n=g_{\Theta(n)}$ for every $n\in {\cal O}b({\cal N})$, where $c_n:G(\Theta(n))\to L_\Theta$ is the canonical morphism. It remains to prove that $g\circ f_O=g_O$ for every $O\in{\cal O}b({\cal U})$. Let $x\in G(O)$ and choose $v:G(O)\to \Theta(n)$. One has
$$g\circ f_O(x)=g(c_n(v(x)))=g_{\Theta(n)}(v(x))=g_O(x)\ .
$$
\qed
 \subsection{Bundle maps between pointed sphere bundles}\label{spheres}
 
Let $X$ be a CW complex and $Y\subset X$ a subcomplex. For two sections $s'$, $s''$ in a an oriented $r$-sphere bundle over   $X$ which coincide over $Y$, we denote by $o_Y(s',s'')\in H^r(X,Y;\Z)$ the {\it primary obstruction} to the existence of a homotopy between $s'$ and $s''$ in the space of sections which coincide  with $\resto{s'}{Y}=\resto{s''}{Y}$ on $Y$ \cite{St}. \\
 
Let  $\pi_\zeta:\zeta\to{\cal B}$ be an oriented real bundle  of rank  $r$ over a CW complex  ${\cal B}$.   Denote by $\pi_\zeta^+: \zeta^+_{\cal B}=:\hat{\cal B}\to {\cal B}$ the bundle projection of the associated sphere bundle, and consider the pull-back bundle $\hat \zeta:=[\pi_\zeta^+]^*(\zeta)$ on $\hat{\cal B}$. The sphere bundle  $\hat \zeta^+_{\hat{\cal B}}=[\pi_\zeta^+]^*(\zeta^+_{\cal B})$ comes with a tautological section $\theta_\zeta$ and an ``infinite" section $s^\infty_{\hat\zeta}$. These sections coincide on the subspace $\infty_\zeta\subset\hat{\cal B}$.  We endow the space $\hat{\cal B}$ with a CW structure in the following way: First, on the subspace $\infty_\zeta$ we copy the CW structure  from the base ${\cal B}$ via $s^\infty_\zeta$. Second, for every $k$-cell $e\subset {\cal B}$ we put $\hat e:=\pi_\zeta^{-1}(e)$.   The attaching map corresponding to $\hat e$ is defined in the following way: let $u:D^k\to\bar e\subset {\cal B}$ the attaching map of  $e$. The pullback bundle $u^*(\zeta)$ is trivial, so it can be identified with $D^k\times \R^r=D^k\times \cringle{D}^r$. The induced map $D^k\times \cringle{D}^r\to \pi_\zeta^{-1}(\bar e)\subset\zeta$ can be extended to map $\hat u:D^k\times {D}^r\to [\pi_\zeta^+]^{-1}(\bar e)\subset\zeta^+$ in an obvious way. Let $\t_\zeta$ be the Thom class of the bundle $\zeta$. We claim
\begin{lm}\label{thom} With respect to such a cellular structure on $\hat {\cal B}$ one has $o_{\infty_\zeta}(s^\infty_{\hat\zeta},\theta_\zeta)=\tg_\zeta$ in $H^r(\hat{\cal B},\infty_\zeta;\Z)$.
\end{lm}
\pf Let $P:\E\to \B:={\mathrm BSO}(r)$ be the universal vector bundle  with structure group ${\mathrm SO}(r)$ and  a fixed  $CW$ structure on the classifying space $\B$.  Since $H^r(\hat \B, \infty_\E;\Z)\simeq H^0(\B;\Z)\simeq\Z$, there exists an integer $N$ such that $o_{\infty_\E}(s^\infty_{\hat \E},\theta_{\E})=N {\mathrm t}_\E$.
Let $f:{\cal B}\to\B$ a cellular map which induces the bundle $\zeta$. This map is covered by a bundle map $\hat f:\hat {\cal B}\to\hat\B$, which is obviously cellular and maps the subcomplex $\infty_\zeta$ of  $\hat {\cal B}$ into the subcomplex $\infty_\E$ of $\hat\B$. Using the functorial properties of the relative obstruction class  and of the Thom class, we obtain $o_{\infty_\zeta}(s^\infty_{\hat\zeta},\theta_\zeta)=N\tg_\zeta$. The integer $N$ can be computed using any bundle $\zeta$, so we will choose the bundle $\R^r\to\{*\}$. The tautological section is just the identity of $[\R^r]^+$. It's easy to see that both classes can be identified with the generator of $H^r([\R^r]^+,\infty;\Z)$.
\qed
\begin{co}\label{thomnew} Let $\zeta$ be an oriented $r$-bundle over a CW complex ${\cal B}$, and let $s$ be a section in $\zeta^+_{\cal B}$ which coincides with $s^\infty_\zeta$ on a subcomplex ${\cal A}\subset{\cal B}$. Then $o_{\cal A}(s^\infty_\zeta,s)=s^*({\mathrm t}_ \zeta)$ in $H^r({\cal B},{\cal A};\Z)$.
\end{co}
\pf Note that, with respect to the cellular decomposition of $\hat{\cal B}$ considered above, the section $s:{\cal B}\to\hat {\cal B}$ is a cellular map and maps to subcomplex ${\cal A}$ into the subcomplex $\infty_\zeta$. It suffices to apply the functorial property of the relative obstruction classes with respect to cellular maps.
\qed
\begin{co}\label{obs} Let $\zeta$ be an oriented $r$-bundle over a finite CW complex ${\cal B}$ of dimension $n\leq r$, and let ${\cal A}\subset{\cal B}$ be a subcomplex. Then the map $o_{\cal A}:s\mapsto s^*({\mathrm t}_ \zeta)$ defines a bijection between the set $\Gamma_{\cal A}(\zeta^+_{\cal B})$ of homotopy classes of sections in $\zeta^+_{\cal B}$ which coincide with $s^\infty_\zeta$ on ${\cal A}$, and $H^r({\cal B},{\cal A};\Z)$.
\end{co}
\pf  Injectivity: Since $\dim({\cal B})\leq r$, for a section $s\in \Gamma_{\cal A}(\zeta^+_{\cal B})$ the only  obstruction to the existence of a homotopy between  $s^\infty_\zeta$ and $s$ is the primary obstruction $o_{\cal A}(s^\infty_\zeta,s)$.  To prove surjectivity, consider, for any $r$-cell $e\subset {\cal B}\setminus{\cal A}$, a section $s_e$ which coincides with $s^\infty_\zeta$ on $\bar e\setminus  e$ and has a single vanishing point, which is non-degenerate. The pull-back $s_e^*({\mathrm t}_ \zeta)$ is a generator of $H^r({\cal B},{\cal B}\setminus e;\Z)\cong \Z$.
\qed

\begin{co}\label{obsbundle} Let $\zeta_0$, $\zeta_1$ be two oriented bundles of ranks $r_0$, $r_1$ over an $n$-dimensional complex ${\cal B}$.
\begin{enumerate}
\item If $n+r_0<r_1$,  any pointed bundle map $f:[\zeta_0]^+_{\cal B}\to [\zeta_1]^+_{\cal B}$ over ${\cal B}$ is homotopic (in the space of pointed bundle maps over ${\cal B}$) to the fiberwise constant map $f^\infty$ which maps   $[\zeta_0]^+$ into $\infty_{\zeta_1}$.
\item If $n+r_0=r_1$, then a pointed bundle map $f:[\zeta_0]^+_{\cal B}\to [\zeta_1]^+_{\cal B}$ over ${\cal B}$ is homotopic to $f^\infty$ if and only if the class $h_f\in H^{n}(B;\Z)$, defined by the condition $f^*({\mathrm t}_{\zeta_1})=[\pi_{\zeta_0}^+]^*(h_f)\cup {\mathrm t}_{\zeta_0}$, vanishes. Moreover, the assignment $f\mapsto h_f$ defines a bijection between the set of homotopy classes of pointed bundle maps $[\zeta_0]^+_{\cal B}\to [\zeta_1]^+_{\cal B}$ and $H^n({\cal B};\Z)$.
\end{enumerate}
\end{co}
\pf It suffices to apply Corollary \ref{obs} to the pull-back bundle $\tilde\zeta_1:=[\pi_{\zeta_0}^+]^*(\zeta_1)$ over   $\tilde{\cal B}:=[\zeta_0]^+_{\cal B}$ and to identify the space of pointed bundle maps  $[\zeta_0]^+_{\cal B}\to [\zeta_1]^+_{\cal B}$ with the space of those sections in $[\tilde\zeta_1]^+_{\tilde{\cal B}}$ which coincide with $s^\infty_{\tilde\zeta_1}$ on $\infty_{\zeta_0}\subset \tilde {\cal B}$. Then use the Thom isomorphism $\cdot \cup {\mathrm t}_{\zeta_0}: H^{n}({\cal B};\Z)\to H^{r_1}(\tilde {\cal B},\infty_{\zeta_0};\Z)$.
\qed

\vspace{10mm}
{\small Christian Okonek: \\
Institut f\"ur Mathematik, Universit\"at Z\"urich,
Winterthurerstrasse 150, CH-8057 Z\"urich,\\
e-mail: okonek@math.unizh.ch
\\  \\
Andrei Teleman: \\
CMI,   Universit\'e de Provence,  39  Rue F. Joliot-Curie, F-13453
Marseille Cedex 13,   e-mail: teleman@cmi.univ-mrs.fr
}

  \end{document}